\documentclass[10pt,a4paper]{article}

\usepackage[utf8]{inputenc}
\usepackage[T1]{fontenc}
\usepackage{amsmath,amssymb,amsfonts,amsthm}
\usepackage{mathrsfs}
\usepackage{dsfont}
\usepackage[dvipsnames]{xcolor}
\usepackage{enumitem}
\usepackage[numbers,comma,sort]{natbib}
\usepackage{microtype}
\usepackage{aliascnt}
\usepackage[
  a4paper,
  textwidth=15cm,
  textheight=23cm,
  centering,
  footskip=1cm
]{geometry}

\definecolor{db}{RGB}{0,0,130}
\usepackage[
  colorlinks=true,
  citecolor=db,
  linkcolor=black,
  urlcolor=black,
  pdfstartview=FitH,
  pdftitle={Fractional SDEs with singular drift and superlinear dissipation},
  pdfauthor={Konstantinos Dareiotis, El Mehdi Haress, Khoa Le}
]{hyperref}
\usepackage[nameinlink,capitalize]{cleveref}

\numberwithin{equation}{section}

\newcommand{\R}{\mathbb{R}}
\newcommand{\N}{\mathbb{N}}
\newcommand{\EE}{\mathbb{E}}
\newcommand{\PP}{\mathbb{P}}
\newcommand{\eps}{\varepsilon}
\newcommand{\ind}{\mathds{1}}

\theoremstyle{plain}
\newtheorem{theorem}{Theorem}[section]
\newaliascnt{proposition}{theorem}
\newtheorem{proposition}[proposition]{Proposition}
\aliascntresetthe{proposition}
\crefname{proposition}{Proposition}{Propositions}
\Crefname{proposition}{Proposition}{Propositions}

\newaliascnt{lemma}{theorem}
\newtheorem{lemma}[lemma]{Lemma}
\aliascntresetthe{lemma}
\crefname{lemma}{Lemma}{Lemmas}
\Crefname{lemma}{Lemma}{Lemmas}

\newaliascnt{corollary}{theorem}
\newtheorem{corollary}[corollary]{Corollary}
\aliascntresetthe{corollary}
\crefname{corollary}{Corollary}{Corollaries}
\Crefname{corollary}{Corollary}{Corollaries}

\theoremstyle{definition}
\newaliascnt{definition}{theorem}
\newtheorem{definition}[definition]{Definition}
\aliascntresetthe{definition}
\crefname{definition}{Definition}{Definitions}
\Crefname{definition}{Definition}{Definitions}

\newaliascnt{assumption}{theorem}
\newtheorem{assumption}[assumption]{Assumption}
\aliascntresetthe{assumption}
\crefname{assumption}{Assumption}{Assumptions}
\Crefname{assumption}{Assumption}{Assumptions}

\theoremstyle{remark}
\newaliascnt{remark}{theorem}
\newtheorem{remark}[remark]{Remark}
\aliascntresetthe{remark}
\crefname{remark}{Remark}{Remarks}
\Crefname{remark}{Remark}{Remarks}

\theoremstyle{definition}
\newaliascnt{notation}{theorem}
\newtheorem{notation}[notation]{Notation}
\aliascntresetthe{notation}
\crefname{notation}{Notation}{Notations}
\Crefname{notation}{Notation}{Notations}

\author{Konstantinos Dareiotis\textsuperscript{1}\footnote{\href{mailto:K.Dareiotis@leeds.ac.uk}{\texttt{K.Dareiotis@leeds.ac.uk}}},
El Mehdi Haress\textsuperscript{2}\footnote{\href{mailto:el-mehdi.haress@inria.fr}{\texttt{el-mehdi.haress@inria.fr}}},
Khoa L\^e\textsuperscript{1}\footnote{\href{mailto:K.Le@leeds.ac.uk}{\texttt{K.Le@leeds.ac.uk}}}
\\[4pt] {\small\textsuperscript{1}University of Leeds, School of Mathematics, Leeds, United Kingdom}
\\ {\small\textsuperscript{2}Universit\'e C\^ote d'Azur, Inria, France}}
\date{\today}
\title{Regularising Flows and Singular Fractional SDEs with Superlinear Dissipation}

\begin{document}

\maketitle

\begin{abstract}
We study stochastic differential equations driven by the fractional Riemann--Liouville process
with Hurst parameter $H\in(0,1)$, whose drift is the sum of a possibly
distributional vector field $b$ and a smooth, one-sided Lipschitz,  dissipative vector field $F$
with coercivity of order $k+1$ and polynomially growing derivatives.  If
$b\in\mathcal C^\gamma(\R^d;\R^d)$ for some $\gamma<1$ satisfying
$\gamma>\max\{1-1/(2H),2-1/H\}$, and if
$H+1/(k-1)>(1-\gamma)H/2$, we construct a solution by
smooth approximation and prove uniqueness in a class defined through local
comparison with the flow generated by the dissipative drift. We also
establish quantitative temporal regularity
and stability estimates.  The proof combines a conditional decomposition of
fractional Brownian motion, Malliavin integration by parts, uniform
non-degeneracy estimates for the regularising flow, and the stochastic sewing
lemma.
\end{abstract}


\medskip

\noindent\textbf{MSC2020 subject classification: }
60H10, 60G22, 60H50.

\medskip

\noindent\textbf{Acknowledgement: } Most of this work was completed while the second author was affiliated with the University of Leeds. All authors gratefully acknowledge support from the Engineering and Physical Sciences Research Council (EPSRC) under grant EP/Y016955/1. The authors used AI-based tools for assistance with editing, proofreading, and checking some technical details.

\tableofcontents

\section{Introduction}\label{sec:introduction}

Let $d\geq 1$, let $t_0<T$, and let $H\in(0,1)$. We study the formal
$d$-dimensional stochastic differential equation
\begin{align}\label{eq:sde-main}
\begin{split}
dX_t
=F(X_t)dt+b(X_t)dt+dV_t+dB_t^{t_0},
\qquad t\in[t_0,T],
\end{split}
\end{align}
with an $\mathcal F_{t_0}$-measurable initial condition $X_{t_0}$.
Here $B^{t_0}$ is a Riemann--Liouville process with Hurst parameter $H$
started at $t_0$ (a Wiener process when $H=1/2$), $V$ is an
$\mathcal F_{t_0}$-measurable continuous input process, $b$ is a possibly
distributional vector field, and $F$ is a smooth dissipative vector field
whose derivatives may grow polynomially.

Equation \eqref{eq:sde-main} is singular when $b$ has negative
Besov--H\"older regularity, since the composition $b(X_t)$ is then not
defined pointwise. Accordingly, the drift integral associated with $b$ is
defined as a limit along smooth approximations. Although the dissipativity
of $F$ prevents explosion, $F$ is not assumed to be globally Lipschitz and
its derivatives may be unbounded. The regularisation provided by the
fractional noise must therefore be combined with quantitative
coming-down-from-infinity estimates for the nonlinear flow generated by
$F$.

\paragraph{Regularisation by noise.}
In the absence of the dissipative component $F$, equation
\eqref{eq:sde-main} has been studied extensively. For Brownian SDEs, the
first regularisation-by-noise results go back to Zvonkin and Veretennikov,
who established well-posedness for bounded measurable drifts
\cite{Zvonkin,Veretennikov}. These results were subsequently extended by
Krylov and R\"ockner to time-dependent integrable drifts satisfying a
subcritical Ladyzhenskaya--Prodi--Serrin condition
\cite{KrylovRockner}. Brownian equations with distributional drifts have
also been considered under various notions of solution; see, among others,
\cite{portenko1990generalized,FRW,BassChen,DelarueDiel,FIR,LeGall,Lejay}.

SDEs with irregular drifts driven by fractional Brownian motion have likewise
received considerable attention, with early contributions including
\cite{NualartOuknine}. The development of new methods has since led to a
renewed interest in such equations. Catellier and Gubinelli introduced an
averaging-along-irregular-curves approach, which makes it possible to
quantify the regularising effect of a path through the properties of its
occupation measure \cite{CatellierGubinelli}. A complementary probabilistic
method is provided by the stochastic sewing lemma, introduced in
\cite{le2020stochastic}, which is particularly well suited to conditional
estimates and can be combined with Davie-type arguments \cite{Davie}.
These techniques have led to increasingly general strong well-posedness
results for fractional equations with irregular or distributional drifts;
we refer to \cite{GaleatiGerencser} and the references therein for recent
developments. They have also been extended to distribution-dependent
equations and equations with generalised drifts
\cite{GHM,Banos}. Weak well-posedness and the identification of the law of
solutions have been investigated in \cite{ButkovskyLeMytnik}.
Regularisation estimates of this kind have also become an important tool in
the numerical approximation of SDEs with irregular or distributional
drifts; see, for example,
\cite{neuenkirch2021euler,butkovsky2021approximation,
de2019numerical,GHR2024}.

\paragraph{Dissipative stochastic equations.}
When $b=0$ and the driving process is Brownian, the well-posedness and
long-time behaviour of dissipative stochastic equations are classical; see,
for example, \cite{markovchains}. Under suitable Lyapunov and dissipativity
conditions, one typically obtains existence and uniqueness of an invariant
measure together with convergence of the law of the solution toward its
stationary regime.

When the equation is driven by fractional Brownian motion, the solution
itself is no longer Markovian because of the temporal correlations of the
noise. Hairer introduced a stochastic-dynamical-systems framework in which
the state is enlarged to include the relevant history of the noise and used
it to study invariant measures and convergence to equilibrium
\cite{Hairer}. This approach has subsequently been developed in
\cite{hairerohashi,hairer9,hairerpillai}. Related works study, among other
questions, the regularity of stationary densities and the approximation of
invariant measures; see \cite{li2023non,cohen,panloup2020general} and the
references therein.

\paragraph{The combined regime.}
Singular interactions and dissipative confinement naturally occur together
in many stochastic systems. In interacting-particle models, for instance, a
singular or strongly repulsive interaction is often accompanied by a
confining force that prevents the particles from escaping to infinity.
Classical examples include Coulomb gases and logarithmic interactions
\cite{serfaty2015coulomb,forrester2010log}. A similar competition between
local singular behaviour and global dissipation appears in stochastic
interface models, constrained dynamics, and equations involving reflection
terms or singular potentials; see, among others,
\cite{BounebacheZambotti,athreya2020well,nualartpardoux,
zambottireflection,scarpa}. In these examples, the singular component
describes a local interaction or constraint, whereas the dissipative
component controls the dynamics at large distances.

Despite the natural coexistence of these two mechanisms, general
well-posedness results usually treat only one of them in full strength. The
theory of regularisation by noise allows highly irregular, and even
distributional, drifts, but typically assumes that any additional regular
component is globally Lipschitz or has uniformly bounded derivatives.
Conversely, the theory of dissipative stochastic equations permits
superlinear and non-globally-Lipschitz vector fields, but generally requires
the remaining terms to be classically defined. Other works concern specific
singular interactions, reflection measures, or model-dependent structures
rather than a general class of singular vector fields.
Fractional SDEs combining a singular drift with a Lipschitz dissipative
component are studied in our previous work \cite{DHL25} and in
\cite{Avi}. In the Brownian case $H=1/2$, Bao and Hao
\cite{bao2024limit} consider a positively H\"older drift together with a
locally Lipschitz component satisfying a dissipativity condition that is
slightly different from ours.

The present paper provides a general framework in which the dissipative
component is coercive of order $k+1>2$ and need not be globally Lipschitz,
while the second drift belongs to a Besov--H\"older class in the
subcritical regime \eqref{eq:subcritical} of regularisation by noise. The fractional results of \cite{DHL25} and
\cite{Avi} concern globally Lipschitz dissipative components, which include
linear dissipation but exclude superlinear growth; for superlinear
dissipative components, the results of the present paper are new for every
$H\ne1/2$. When $H=1/2$,
the condition \eqref{eq:subcritical} forces $\gamma>0$, and this case is
close to that of \cite{bao2024limit}, whose dissipativity condition on the
regular part of the drift differs slightly from ours. For $H<1/2$, the
framework allows, in particular, genuinely distributional drifts.

\paragraph{The role of path input in \eqref{eq:sde-main}.} 
Besides well-posedness, our estimates are designed to serve as a foundation
for the study of long-time properties such as invariant measures and
convergence to equilibrium. We do not pursue these questions here, but they
motivate the formulation adopted below. A direct application of the
stochastic-dynamical-systems framework of \cite{Hairer} would require
suitable continuity properties of the solution with respect to the entire
driving noise, which are particularly inconvenient to establish in the
presence of an irregular drift. As in \cite{DHL25}, our intended
approach is instead to couple the solution with the past of the noise and to
evolve the two components jointly.

After fixing an initial time $t_0$, the future of a two-sided fractional
Brownian motion decomposes into a contribution determined by its past and a
fresh Riemann--Liouville process independent of the past. In
\eqref{eq:sde-main}, the fresh part is $B^{t_0}$, and the contribution of
the past corresponds to the particular input $V=\bar B^{t_0}$ (see
\Cref{sec:framework}); we allow a general $\mathcal F_{t_0}$-measurable
input $V$, which is the form needed for restarting. The same decomposition, applied at any later time $u$, splits
the Riemann--Liouville increment into an $\mathcal F_u$-measurable history
and a fresh innovation independent of $\mathcal F_u$. Freezing the state and
the accumulated history at time $u$ therefore yields exact fixed-time
conditioning identities, without requiring a pathwise Markov property; at
$H=1/2$ the history vanishes and these identities reduce to the usual
restarting property of Brownian motion. It is therefore important to
establish well-posedness and stability for a prescribed initial condition
and a prescribed input $V$, with estimates that remain suitable for
restarting and concatenating the dynamics.

\paragraph{Main results.}
We assume that $b\in\mathcal C^\gamma(\R^d;\R^d)$ for some $\gamma<1$ with
$(1-\gamma)H<1/2$ and $(2-\gamma)H<1$; see \eqref{eq:subcritical}. The first
condition is the Catellier--Gubinelli regularisation threshold. The second
one is redundant when $H\le1/2$ and is an additional integrability condition
when $H>1/2$. The drift $F$ has coercivity of order $k+1$, polynomial
derivative bounds, and a global one-sided Lipschitz estimate
(\Cref{ass:dissipative-drift}). Finally, we assume the compatibility 
condition $H+1/(k-1)>(1-\gamma)H/2$, which holds automatically whenever
$\gamma\ge-1$ (in particular whenever $H\ge1/4$) and, under the subcritical
conditions above, whenever $k\le5$.

Under these assumptions, \Cref{thm:main} provides a solution which, for
every finite moment order $m$, satisfies
$X-V-B^{t_0}\in C^{1+(\gamma\wedge0)H}([t_0,T];L_m)$, the temporal
regularity predicted by the regularisation threshold. More importantly, the
solution satisfies a conditional local comparison estimate with the flow
generated by $F$, and it is pathwise unique among solutions satisfying this
comparison estimate, for a loss $\eps$ of H\"older regularity in an explicit
range depending on $H$, $\gamma$, and $k$, and having a marginal moment of
order slightly above four. \Cref{thm:stability} gives a quantitative stability estimate with
respect to the initial condition and the singular drift, for a fixed input
$V$.

\paragraph{Main difficulties and strategy.}
The obstruction to extending the argument of \cite{DHL25} is the following.
For a globally Lipschitz $F$, drift differences satisfy
$|F(x)-F(y)|\le L|x-y|$, so they can be estimated at the same moment
order as the state difference. Here the corresponding bound is only
$|F(x)-F(y)|\lesssim(1+|x|^{k-1}+|y|^{k-1})|x-y|$. The unbounded state
factors prevent the same direct closure: applying H\"older's inequality
requires higher moments of the unknown difference. The one-sided Lipschitz
condition controls the signed energy contribution of $F$, but does not
supply the absolute-value estimates needed in all the sewing terms.
Moreover, when $b$ is distributional, the usual energy calculation must
itself be justified through approximation and regularisation estimates.

We therefore incorporate the entire dissipative drift into the reference
dynamics. The regularising flow $\phi^{s,x,\mathbf v}$ solves the equation
with drift $F$ and the same additive signal, but without $b$; see
\eqref{eq:flow-front}. Comparing $X$ with this flow preserves the one-sided
cancellation and isolates the contribution of the singular drift.
Coming down from infinity controls the flow uniformly in its initial point
after any positive time. Since the reference dynamics is nonlinear and no
longer a Gaussian shift, Gaussian smoothing is replaced by Malliavin
integration by parts for the flow. The sewing defects involve differences of
evaluations of $b$ at two and four points; interpolation between their
arguments produces affine combinations of flows, sometimes with negative
coefficients. This explains the need for spatial derivative bounds and
non-degeneracy estimates for these affine combinations. The Gaussian
regularisation estimates of \cite{DHL25} remain useful for temporal
regularity and passage to the limit; it is the direct closure of the a
priori and stability arguments that must be replaced here. One could
instead truncate $F$ and localise with stopping times, but stopping times
are not compatible with the conditional estimates on which stochastic
sewing relies, and the growth of the truncated solutions would still have
to be controlled. The flow-based approach avoids this: the dissipative
dynamics is handled through results on equations without the singular
drift, and the sewing estimates only involve these flows. Similar
coming-down-from-infinity estimates play an important role in dissipative
stochastic reaction--diffusion equations; see, for instance,
\cite{moinatweber,MourratWeber,pavlosweber}.

The comparison with the regularising flow is imposed only below a local time
scale $\tau_{\mathbf v;I}$, defined in \eqref{eq:local-scale}, at which the
H\"older increment of the input is balanced by the coming-down-from-infinity
scale of the polynomial drift. On good events, where the interval is shorter
than this scale, the solution is close to the restarted regularising flow.
On bad events, the local-scale inequality produces an arbitrarily high power
of the interval length, at the cost of Gaussian moments of the H\"older
seminorm of the history.

The proof proceeds in four steps. First, we develop Malliavin calculus for
the Riemann--Liouville process and prove integration-by-parts estimates at
the critical regularity for non-degenerate random variables. Second, we
estimate the regularising flow: coming down from infinity, Malliavin
derivatives up to the finite order required by the integration by parts,
spatial and four-point differences, and uniform inverse moments of the
Malliavin covariance of affine combinations of flows. Third, we insert these
estimates into conditional stochastic sewing arguments, in which the
singular drift is evaluated along the flow restarted at the left endpoint of
each increment. This yields uniform temporal regularity and comparison
estimates for smooth approximations of $b$, and, after a dissipative
pathwise estimate, stability with respect to the state and the drift. The
resulting local recursion is iterated on a deterministic mesh; this is where
the condition $H+1/(k-1)>(1-\gamma)H/2$ is used. Finally, heat
regularisations of $b$ generate a Cauchy family of smooth solutions; the
uniform estimates pass to the limit, the limiting drift process is shown to
be independent of the approximation, and the stability estimate yields
uniqueness.

\paragraph{Organisation of the paper.}
\Cref{sec:framework} introduces the fractional-noise decomposition, the
Besov--H\"older spaces, the definition of a solution, the structural
assumptions, the local scale, the solution classes, and the main results.
\Cref{sec:preliminaries} develops the Malliavin calculus, proves the
integration-by-parts estimates, records the conditioning identities used
when restarting the equation, collects Gaussian estimates for the
history of the noise, and recalls the stochastic sewing lemma. The regularising-flow estimates are established in
\Cref{sec:flow}. In \Cref{sec:a-priori}, we prove a priori estimates for
smooth approximations and stability estimates for solutions with singular
drifts. The existence, uniqueness, and stability proofs are given in
\Cref{sec:proof-main}.

\section{Framework and main results}\label{sec:framework}

This section fixes the probabilistic and analytic framework for the equation.
We first describe the fractional-noise decomposition and the notation used
throughout, then introduce the solution concept, the structural assumptions,
and the solution classes, and finally state the main results.

\subsection{Notation}

Let $(\Omega,\mathcal F,\mathbb F=(\mathcal F_t)_{t\in\R},\PP)$ be a filtered
probability space satisfying the usual conditions. We write
$\EE^t[\,\cdot\,]:=\EE[\,\cdot\mid\mathcal F_t]$ and denote the law of a random
variable $Z$ by $\mathcal L(Z)$. For $m\in[1,\infty]$, the norm in
$L_m(\Omega)$ is denoted by $\|\cdot\|_{L_m}$; the underlying probability space
will be omitted from the notation.

Let $(W_t)_{t\in\R}$ be a two-sided $d$-dimensional $\mathbb F$-Wiener
process; that is, $W_0=0$, $W$ has continuous adapted paths, and, for
every $s<t$, the increment $W_t-W_s$ is independent of $\mathcal F_s$ and has
law $\mathcal N(0,(t-s)I_d)$.

Fix $H\in(0,1)$.  When $H\ne1/2$, let $c_H>0$ denote the normalising
constant in the Mandelbrot--Van Ness representation and, for $s\in\R$ and
$t\ge s$, define the Riemann--Liouville process started at $s$ and the
history at time $s$ by
\begin{align*}
 B_t^s:=c_H\int_s^t(t-r)^{H-\frac12}\,dW_r,
 \qquad
 \bar B_t^s
 :=c_H\int_{-\infty}^s
 \Big((t-r)^{H-\frac12}-(s-r)^{H-\frac12}\Big)\,dW_r.
\end{align*}
The history $\bar B^s$ is $\mathcal F_s$-measurable, whereas $B^s$ is
independent of $\mathcal F_s$. If $B$ denotes the two-sided fractional
Brownian motion given by the Mandelbrot--Van Ness representation with the
same $W$, then
\begin{align}\label{eq:decomp-fBm}
 B_t-B_s=\bar B_t^s+B_t^s,
 \qquad t\ge s.
\end{align}
The process $B$ itself is not used in our analysis herein; \eqref{eq:decomp-fBm} only serves
to interpret \eqref{eq:sde-main} when the path input is given by $\bar B^{t_0}$, see the discussion in \cref{sec:introduction} and \Cref{rem:full-fbm} below.
Similarly, for $s\le u\le t$, the past accumulated between $s$ and $u$,
viewed from time $u$ onwards, is
\begin{align*}
 \bar B_t^{s,u}
 :=c_H\int_s^u
 \Big((t-r)^{H-\frac12}-(u-r)^{H-\frac12}\Big)\,dW_r,
\end{align*}
and
\begin{align}\label{eq:restart-RL}
\begin{split}
 B_t^s-B_u^s=\bar B_t^{s,u}+B_t^u,
 \qquad t\ge u.
\end{split}
\end{align}
When $H=1/2$, we set $B_t:=W_t$, $B_t^s:=W_t-W_s$, and
$\bar B_t^s=\bar B_t^{s,u}:=0$; with these definitions, both
\eqref{eq:decomp-fBm} and \eqref{eq:restart-RL} remain valid. Since $W$ is
two-sided, the histories $\bar B^{s}$ and $\bar B^{s,u}$ are defined on the
whole half-line $[u,\infty)$ (with $u=s$ for $\bar B^s$). Splitting the
integration range $[s,c]$ at $u$ in the definition gives the cocycle
identity
\begin{align}\label{eq:history-cocycle}
\begin{split}
 \bar B_r^{s,c}
 =\bar B_r^{s,u}-\bar B_c^{s,u}+\bar B_r^{u,c},
 \qquad s\le u\le c\le r,
\end{split}
\end{align}
which is what makes restarting at successive times $u\le c$ consistent: the
input frozen at time $c$ equals the input frozen at time $u$, re-centred at
$c$, plus the history accumulated on $[u,c]$.

For a path $f$ and $s\le t$, we write $f_{s,t}:=f_t-f_s$, and for an interval
$I\subset\R$ we let $\Delta_I:=\{(s,t)\in I^2:s<t\}$.
We denote by $C(I;\R^d)$ the space of continuous paths with the topology of
uniform convergence on compact subsets. For $\alpha\in(0,1]$, we write
$C^\alpha(I;\R^d)$ for the temporal H\"older space with seminorm
\[
 [f]_{C^\alpha(I)}
 :=\sup_{(s,t)\in\Delta_I}\frac{|f_{s,t}|}{|t-s|^\alpha}.
\]
A path defined only on $[u,\infty)\cap I$ for some $u\in I$, such as a
history $\bar B^{s,u}$ restarted at time $u$, is extended to the left of $u$
by the constant value at $u$ whenever its seminorm on $I$ is taken. This
extension does not change the H\"older seminorm, and it never takes place
outside $I$.
For an $L_m$-valued path $f$ on $I$, we similarly write
\[
 [f]_{C^\alpha(I;L_m)}
 :=\sup_{(s,t)\in\Delta_I}\frac{\|f_{s,t}\|_{L_m}}{|t-s|^\alpha}.
\]
For $0<\beta<\alpha\le1$, we shall use the elementary interpolation
estimate
\begin{align}\label{eq:holder-interpolation}
 [f]_{C^\beta(I;L_m)}
 \le2\Bigl(\sup_{r\in I}\|f_r\|_{L_m}\Bigr)^{1-\beta/\alpha}
 [f]_{C^\alpha(I;L_m)}^{\beta/\alpha}.
\end{align}

We use $\mathcal C^\alpha=\mathcal C^\alpha(\R^d)$ for spatial
Besov--H\"older spaces with the convention $\mathcal C^0=L_\infty$ for bounded
measurable functions. For $\alpha\in(0,1]$,
$\|f\|_{\mathcal C^\alpha}:=\|f\|_\infty
+\sup_{x\ne y}|f(x)-f(y)|/|x-y|^\alpha$.
For $\alpha<0$, a Schwartz distribution $f$ belongs to $\mathcal C^\alpha$ if
\begin{align}\label{eq:alpha-norm}
\begin{split}
 \|f\|_{\mathcal C^\alpha}
 :=\sup_{\delta\in(0,1]}\delta^{-\alpha/2}\|P_\delta f\|_\infty<\infty,
\end{split}
\end{align}
where $P_\delta f:=p_\delta*f$ is the heat semigroup, with
$p_\delta(x):=(4\pi\delta)^{-d/2}e^{-|x|^2/(4\delta)}$.
We write $C_b^\infty(\R^d;\R^d)$ for smooth bounded vector fields whose
derivatives of every order are bounded.

Throughout the paper, $C$ denotes a finite positive constant that may change
from line to line. Its dependence is restricted to the parameters specified in
the corresponding statement. We write $A\lesssim B$ when $A\le CB$ for such a
constant.

\subsection{Definition of a solution}

We next formulate the equation and specify how the distributional drift is
interpreted. Let $t_0<T$, let $X_{t_0}$ be an $\mathcal F_{t_0}$-measurable
$\R^d$-valued random variable, and let $V$ be an $\mathcal F_{t_0}$-measurable
random variable with values in $C([t_0,T];\R^d)$. The integral form of
\eqref{eq:sde-main} is formally
\begin{align*}
 X_t
 =X_{t_0}+\int_{t_0}^tF(X_r)\,dr
 +\int_{t_0}^tb(X_r)\,dr
 +V_{t_0,t}+B_t^{t_0}.
\end{align*}
The distributional drift integral is defined by approximation.

\begin{definition}[Solution by approximation]\label{def:solution}
Let $\gamma<1$, let $b\in\mathcal C^\gamma$, and let $F:\R^d\to\R^d$ be continuous. An
adapted continuous process $X=(X_t)_{t\in[t_0,T]}$ is a solution of
\eqref{eq:sde-main}, with respect to
$(\Omega,\mathcal F,\mathbb F,\PP)$ and $W$, if there exists an adapted
continuous process $K=(K_t)_{t\in[t_0,T]}$ with $K_{t_0}=0$ such that:
\begin{enumerate}[label=(\roman*)]
 \item almost surely, for every $t\in[t_0,T]$,
 \begin{align*}
  X_t
  =X_{t_0}+\int_{t_0}^tF(X_r)\,dr+K_t+V_{t_0,t}+B_t^{t_0};
 \end{align*}
 \item for every sequence $(b^N)_{N\ge1}\subset C_b^\infty(\R^d;\R^d)$
 with $\sup_{N\ge1}\|b^N\|_{\mathcal C^\gamma}<\infty$ and
 $b^N\to b$ in $\mathcal C^{\gamma-1}$, one has
 \begin{align*}
  \sup_{t\in[t_0,T]}
  \left|\int_{t_0}^tb^N(X_r)\,dr-K_t\right|
  \longrightarrow0
 \end{align*}
 in probability as $N\to\infty$.
\end{enumerate}
\end{definition}

The process $K$ is the drift integral associated with $b$ along $X$. We call
a sequence $(b^N)$ as in (ii) an admissible approximation of $b$; an example
is $b^N=P_{1/N}b$.  If $\gamma>0$, Besov interpolation between the uniform
$\mathcal C^\gamma$ bound and convergence in $\mathcal C^{\gamma-1}$ shows
that every admissible sequence converges to $b$ uniformly.  Hence, when
$\gamma>0$, the definition agrees with the classical notion of a solution and
$K_t=\int_{t_0}^tb(X_r)\,dr$.

\subsection{Structural assumptions}

The singular drift is measured by Besov--H\"older regularity,
whereas the smooth drift is controlled through coercivity, polynomial
growth, and a global one-sided Lipschitz bound.

Fix $\gamma<1$. We work in the subcritical regime
\begin{align}\label{eq:subcritical}
\begin{split}
 (1-\gamma)H<\frac12,
 \qquad
 (2-\gamma)H<1,
\end{split}
\end{align}
or equivalently $\gamma>\max\{1-1/(2H),\,2-1/H\}$.
The second inequality in \eqref{eq:subcritical} is redundant for
$H\le1/2$, while for $H>1/2$ it is the stronger one.  In particular, when
$H=1/2$ the condition reads $\gamma>0$, and when $H>1/2$ it forces
$\gamma>0$. We set $\alpha_\gamma:=1+(\gamma\wedge0)H$; this is the
temporal regularity exponent of the solution with the additive signal
removed, see \Cref{thm:main}.

\begin{assumption}[Singular drift]\label{ass:singular-drift}
There exists $\Xi>0$ such that $b\in\mathcal C^\gamma(\R^d;\R^d)$ and
$\|b\|_{\mathcal C^\gamma}\le\Xi$.
\end{assumption}

In particular, the heat approximations satisfy
$\sup_{N\ge1}\|P_{1/N}b\|_{\mathcal C^\gamma}\lesssim\Xi$ and
$P_{1/N}b\to b$ in $\mathcal C^{\gamma-1}$, so they form an admissible
approximation of $b$.

The coming-down-from-infinity estimate used in the analysis of the
regularising flow requires coercivity of order $k+1$. Only finitely many
derivative bounds are needed: we use orders up to $\lfloor-\gamma\rfloor+5$.

\begin{assumption}[Dissipative drift]\label{ass:dissipative-drift}
Let $k>1$. The vector field $F\in C^\infty(\R^d;\R^d)$ satisfies the following conditions:
\begin{enumerate}[label=(\roman*)]
 \item There exist $\kappa_1>0$ and $\kappa_2\ge0$ such that
\begin{align}\label{eq:Fdissipative}
\begin{split}
  \langle F(x),x\rangle
  \le-\kappa_1|x|^{k+1}+\kappa_2,
  \qquad x\in\R^d.
\end{split}
\end{align}
 \item There exists $\kappa_3\ge0$ such that, for every
 $j\in\{0,\ldots,\lfloor-\gamma\rfloor+5\}$ and
 $x\in\R^d$,
\begin{align}\label{eq:growth}
\begin{split}
  |\nabla^jF(x)|
  \le\kappa_3\big(1+|x|^{(k-j)\vee0}\big),
\end{split}
\end{align}
 where $\nabla^0F:=F$.
 \item There exists $\kappa_4\ge0$ such that, for every $x,v\in\R^d$,
\begin{align}\label{eq:one-sided-lipschitz}
\begin{split}
  v^\top\nabla F(x)v\le\kappa_4|v|^2.
\end{split}
\end{align}
\end{enumerate}
We write $\boldsymbol\kappa:=(\kappa_1,\kappa_2,\kappa_3,\kappa_4)$.
\end{assumption}

Integrating \eqref{eq:one-sided-lipschitz} along the segment joining $x$
and $y$ gives
\begin{align}\label{eq:F-one-sided-global}
\begin{split}
 \langle F(x)-F(y),x-y\rangle
 \le\kappa_4|x-y|^2,
 \qquad x,y\in\R^d.
\end{split}
\end{align}

\begin{remark}[Derivative orders]\label{rem:R-gamma-choice}
The integration-by-parts estimates at the critical regularity use
$\lfloor-\gamma\rfloor+3>2-\gamma$ differentiations. Two additional
spatial derivative levels suffice for the four-point flow estimates, which
explains the bound $\lfloor-\gamma\rfloor+5$ in the assumption. This is a
convenient sufficient choice, not an optimal one.
\end{remark}

\begin{remark}[Examples]\label{rem:examples}
For $k>1$, the vector field $F(x)=-x(1+|x|^2)^{(k-1)/2}$ satisfies
\Cref{ass:dissipative-drift} with $\kappa_4=0$: it is smooth,
$\langle F(x),x\rangle\le-|x|^{k+1}$, its derivatives of order $j$ grow like
$(1+|x|)^{k-j}$, and
$\nabla F(x)=-(1+|x|^2)^{(k-1)/2}I_d-(k-1)(1+|x|^2)^{(k-3)/2}xx^\top\le0$.
When $k$ is an odd integer, the polynomial field $F(x)=-x|x|^{k-1}$ can be
used as well; for other values of $k$ it is not smooth at the origin. More
generally, one may take $F=-\nabla U+G$, where $U\in C^\infty$ satisfies
$\nabla^2U\ge-a_1I_d$, $\langle\nabla U(x),x\rangle\ge a_2|x|^{k+1}-a_3$ for
some constants $a_1,a_3\ge0$ and $a_2>0$, and
the bounds \eqref{eq:growth} for $\nabla U$ in place of $F$, and $G$ is
smooth with bounded derivatives. Typical singular drifts are distributional derivatives of
bounded H\"older functions, $b=\partial_{x_i}f$ with $f\in\mathcal C^{1+\gamma}$,
and, in a more specific regime, $b=e\,\mu$, where $e\in\R^d$ is fixed and
$\mu$ is the Dirac mass at the origin, which belongs to $\mathcal C^{-d}$; \eqref{eq:subcritical} then
requires $H<1/(2(d+1))$, and \eqref{eq:H-k-iteration-condition} requires,
for $d\ge2$, $H<2/((d-1)(k-1))$, which is an additional restriction only
when $(d-1)(k-1)>4(d+1)$.
\end{remark}

\subsection{Regularising flow, local scale, and solution classes}

The a priori estimates compare the solution with the equation obtained by
removing the singular drift. We therefore introduce the corresponding
regularising flow, the local time scale below which this comparison is
imposed, and the solution classes that quantify the comparison remainder.

For $s<T$, $x\in\R^d$, and a continuous path
$\zeta\in C([s,T];\R^d)$, let $\mathcal S_{F,s}(x,\zeta)$ be the unique
solution of the deterministic integral equation
\begin{align*}
 \mathcal S_{F,s}(x,\zeta)_t
 =x+\int_s^tF\bigl(\mathcal S_{F,s}(x,\zeta)_r\bigr)\,dr
 +\zeta_{s,t},
 \qquad t\in[s,T].
\end{align*}
Smoothness of $F$ gives local existence, uniqueness, and continuous dependence
on $(x,\zeta)$, while \eqref{eq:Fdissipative} prevents finite-time explosion.
For a deterministic path $\mathbf v\in C([s,T];\R^d)$, the regularising flow
is defined by
\begin{align}\label{eq:flow-front}
\begin{split}
 \phi^{s,x,\mathbf v}(\omega)
 :=\mathcal S_{F,s}(x,\mathbf v+B^s(\omega)),
 \quad\text{that is,}\quad
 \phi_t^{s,x,\mathbf v}(\omega)
 =x+\int_s^tF(\phi_r^{s,x,\mathbf v}(\omega))\,dr
 +\mathbf v_{s,t}+B_t^s(\omega).
\end{split}
\end{align}
Furthermore, if $\xi$ and $\mathbf w$ are $\mathcal F_s$-measurable random variables with
values in $\R^d$ and $C([s,T];\R^d)$, respectively, we write
$\phi_t^{s,\xi,\mathbf w}(\omega)
:=\mathcal S_{F,s}\bigl(\xi(\omega),\mathbf w(\omega)+B^s(\omega)\bigr)_t$.
The continuity of the deterministic solution map makes this evaluation
measurable, $\phi^{s,\xi,\mathbf w}$ is an adapted process and it is the unique $\omega$-wise solution of
$Z_t=\xi+\int_s^tF(Z_r)\,dr+\mathbf w_{s,t}+B_t^s$. In particular,
expressions such as $\phi_t^{s,X_s,V+\bar B^{t_0,s}}$ are unambiguous.

For $M\ge1$, $\mathbf x=(x_1,\ldots,x_M)\in(\R^d)^M$, and $\lambda$ in the
affine set
$\mathcal A_M:=\{\lambda\in\R^M:\sum_{i=1}^M\lambda_i=1\}$,
we write
\begin{align*}
 \Phi_t^{s,\mathbf x,\lambda,\mathbf v}
 :=\sum_{i=1}^M\lambda_i\phi_t^{s,x_i,\mathbf v}.
\end{align*}
The coefficients may be negative. Such affine combinations arise when the
arguments of two- and four-point differences of $b$ are interpolated.

For $\eps\in(0,H)$, set $\chi_\eps:=H-\eps+\frac1{k-1}>0$. For an interval $I$ of length at most one and a continuous $\mathbf{v}:I\to\R^d$,
we define the local time scale
\begin{align}\label{eq:local-scale}
 \tau_{\mathbf v;I}
 :=\big(1\vee[\mathbf v]_{C^{H-\eps}(I)}\big)^{-1/\chi_\eps}
 =\big(1\vee[\mathbf v]_{C^{H-\eps}(I)}\big)^{-\frac{k-1}{1+(k-1)(H-\eps)}}
 \in[0,1],
\end{align}
with the convention $\tau_{\mathbf v;I}:=0$ when
$[\mathbf v]_{C^{H-\eps}(I)}=\infty$.
The parameter $\eps$ is always clear from the context and is omitted from the
notation. Below this scale, the H\"older increment of the input is balanced
by the coming-down-from-infinity scale of the polynomial drift; see
\Cref{lem:coming-down-infinity}. For $s\in[t_0,T]$ we write
$I_s:=[s,(s+1)\wedge T]$.

\begin{definition}[Solution classes]\label{def:solution-classes}
Let $X$ be a solution to \eqref{eq:sde-main} on $[t_0,T]$.
\begin{enumerate}[label=(\roman*)]
 \item For $m\in[2,\infty)$, we say that $X\in\mathcal V(m)$ if
 $[X-V-B^{t_0}]_{C^{1/2}([t_0,T];L_m)}<\infty$, and we write
 $\mathcal V(<\infty):=\bigcap_{m\in[2,\infty)}\mathcal V(m)$.
 \item Set $p_k:=8\bigl((k-2)\vee1\bigr)$. For $\eps\in(0,H)$, we
 write $X\in\mathcal V_\eps(\infty)$ if, for every $p\in[2,p_k]$, there
 exists a constant $C$ such that
\begin{align}\label{eq:classVinf}
\begin{split}
 \sup_{\substack{(s,t)\in\Delta_{[t_0,T]}\\ t-s\le1}}
 \frac{
 \left\|
 \ind_{\{t-s\le\tau_{V+\bar B^{t_0,s};I_s}\}}
 \big(\EE^s|X_t-
 \phi_t^{s,X_s,V+\bar B^{t_0,s}}|^p\big)^{1/p}
 \right\|_{L_\infty}}
 {(t-s)^{\alpha_\gamma}}
 \le C.
\end{split}
\end{align}
\end{enumerate}
\end{definition}

The symbol $\infty$ refers to the outer $L_\infty$ norm in
\eqref{eq:classVinf}. Only conditional moments up to the finite order $p_k$
are required: this is the largest conditional moment of the comparison
remainder, and of the polynomial factors produced by $\nabla F$ and
$\nabla^2F$, that appears in the stability argument of \Cref{sec:stab}; the
count is given at the beginning of that subsection.

\subsection{Main results}

We assume uniform moment control of the random input at every positive loss
of H\"older regularity.

\begin{assumption}[Input]\label{ass:input}
For every $\eps\in(0,H)$, there exists $\Lambda_{V,\eps}\ge0$ such that, for
every $q\in[2,\infty)$,
\begin{align}\label{eq:assumpV}
\begin{split}
 \sup_{s\in[t_0,T]}
 \left\|[V]_{C^{H-\eps}(I_s)}\right\|_{L_q}
 \le\Lambda_{V,\eps}\sqrt q.
\end{split}
\end{align}
\end{assumption}

Covering $[t_0,T]$ by finitely many overlapping intervals of length at most
one and splitting increments at the overlap points shows that
\eqref{eq:assumpV} also implies
$\|[V]_{C^{H-\eps}([t_0,T])}\|_{L_q}\le C_{T-t_0}\Lambda_{V,\eps}\sqrt q$
for every $q\ge2$.

\begin{theorem}\label{thm:main}
Let $H\in(0,1)$ and let $\gamma<1$ satisfy \eqref{eq:subcritical}.
Assume that $b$ and $F$ satisfy
\Cref{ass:singular-drift,ass:dissipative-drift}, and that the exponent
$k>1$ in \Cref{ass:dissipative-drift} satisfies
\begin{align}\label{eq:H-k-iteration-condition}
 H+\frac1{k-1}>\frac{(1-\gamma)H}{2}.
\end{align}
Let $t_0<T$, let
$X_{t_0}$ be $\mathcal F_{t_0}$-measurable and have finite moments of every
order, and let $V$ be an
$\mathcal F_{t_0}$-measurable continuous random path satisfying
\Cref{ass:input}.

Then there exists a solution $X=(X_t)_{t\in[t_0,T]}$ of
\eqref{eq:sde-main} in the classes
$\mathcal V(<\infty)$ and
$\mathcal V_\eps(\infty)$, for every $\eps\in(0,H)$.

Pathwise uniqueness holds among all solutions $Y$ of the same
equation, with the same initial condition $X_{t_0}$, input $V$, driving Wiener
process, and Riemann--Liouville process, for which there exist $\delta_0>0$
and
\begin{align*}
 \eps\in \left(0,H \wedge
 \left(H+\frac1{k-1}-\frac{(1-\gamma)H}{2}\right)\right)
\end{align*}
such that $Y\in\mathcal V_\eps(\infty)$ and
$\sup_{t\in[t_0,T]}\|Y_t\|_{L_{4+\delta_0}}<\infty$.

Moreover, for every $m\in[2,\infty)$ and every $\eps\in(0,H)$, there exist
$q\ge m$ and a constant $C=C(\eps,m,\gamma,H,k,d,\boldsymbol\kappa,\Xi,
 \Lambda_{V,\eps},T-t_0)$ such that
\begin{align}\label{eq:thm-reg}
\begin{split}
 [X-V-B^{t_0}]_{C^{\alpha_\gamma}([t_0,T];L_m)}
 \le C\left(
 1+\|X_{t_0}\|_{L_{q}}
 +\left\|[V]_{C^{H-\eps}([t_0,T])}\right\|_{L_{q}}
 \right)^C.
\end{split}
\end{align}
\end{theorem}

The preceding theorem gives well-posedness for each admissible initial
condition and singular drift. The next result compares the solutions
obtained when these two data are perturbed.

\begin{theorem}[Stability]\label{thm:stability}
Let $H\in(0,1)$ and let $\gamma<1$ satisfy \eqref{eq:subcritical}. Let
$F$ satisfy \Cref{ass:dissipative-drift}, and assume that the exponent
$k>1$ satisfies \eqref{eq:H-k-iteration-condition}. Let $t_0<T$, and let
$V$ be an $\mathcal F_{t_0}$-measurable continuous random path satisfying
\Cref{ass:input}.
Let $\Xi>0$, let $b,g\in\mathcal C^\gamma(\R^d;\R^d)$ satisfy
$\|b\|_{\mathcal C^\gamma}+\|g\|_{\mathcal C^\gamma}\le\Xi$,
and let $X_{t_0},Y_{t_0}$ be $\mathcal F_{t_0}$-measurable random variables with
finite moments of every order. Denote by $X$ and $Y$ the respective unique
solutions given by \Cref{thm:main}, driven by the same input $V$ and
the same Riemann--Liouville process.

Then, for every $m\in[2,\infty)$ and $\theta\in(0,1)$, there exist
$q\ge2m$ and a constant $C$ depending only on
$m,\theta,\gamma,H,k,d,\boldsymbol\kappa,\Xi,T-t_0,
\|X_{t_0}\|_{L_q},\|Y_{t_0}\|_{L_q}$,
and on the input moment bounds in \eqref{eq:assumpV},
such that
\begin{align}\label{eq:main-stability}
 \sup_{t\in[t_0,T]}\|X_t-Y_t\|_{L_{2m}}
 \le C\left(
 \|X_{t_0}-Y_{t_0}\|_{L_{2m}}
 +\|b-g\|_{\mathcal C^{\gamma-1}}
 \right)^{1-\theta}.
\end{align}
\end{theorem}

\begin{remark}\label{rem:full-fbm}
Taking $V=\bar B^{t_0}$ in \eqref{eq:sde-main} recovers an equation driven by
the full fractional Brownian increment $B_t-B_{t_0}$, by
\eqref{eq:decomp-fBm}. For this choice, \Cref{ass:input} follows from
\Cref{lem:history-holder-moments} below, applied with $u=t_0$ and $u'=s$,
since $I_s\subset[s,s+1]$.
\end{remark}

\begin{remark}
Condition \eqref{eq:H-k-iteration-condition} expresses the compatibility
between the two mechanisms underlying the proof. On the one hand, the noise
must be sufficiently rough to regularise the singular drift. On the other
hand, the dissipative dynamics must lose its dependence on large initial
conditions sufficiently rapidly. The parameters $H$ and $k$ quantify these
two effects, respectively. Condition \eqref{eq:H-k-iteration-condition}
holds automatically whenever $\gamma\ge-1$, since then
$(1-\gamma)H/2\le H$. Moreover, the first inequality in
\eqref{eq:subcritical} gives $(1-\gamma)H/2<1/4$, so that
\eqref{eq:H-k-iteration-condition} also holds automatically whenever
$H+\frac1{k-1}\ge\frac14$, in particular for every $k\le5$. In general,
\eqref{eq:H-k-iteration-condition} is equivalent to
$\gamma>-1-\frac{2}{(k-1)H}$, and under \eqref{eq:subcritical} it is a
genuine restriction only when $k>5$ and $H<\frac14-\frac1{k-1}$. In the
proof, the condition is used
only at the end of the stability argument, \Cref{prop:stab}, on which the
existence, uniqueness, and stability statements all rely. There, local
estimates on a mesh of size $h$ gain only a factor $h^{1-(1-\gamma)H}$, so
that their iteration produces a factor $\exp(Ch^{-(1-\gamma)H})$, whereas
the contributions of the events on which an input exceeds its local scale
decay like $\exp(-ah^{-2\chi_\eps})$. Choosing $h$ so that the latter matches
the target error, the former is smaller than any negative power of that
error exactly when $2\chi_\eps>(1-\gamma)H$, that is,
$\eps<H+\frac1{k-1}-\frac{(1-\gamma)H}2$; condition
\eqref{eq:H-k-iteration-condition} ensures that such $\eps\in(0,H)$ exist.
We do not know whether the condition is necessary.
\end{remark}

\section{Preliminaries}\label{sec:preliminaries}

The proof relies on the following ingredients. We first develop the Malliavin
calculus and the integration-by-parts estimates at the critical regularity
needed to exploit the smoothing effect of the fresh fractional noise. We then
record the conditioning identities that allow these estimates to be used
after restarting the equation. Finally, we collect Gaussian estimates for the
history of the noise and recall the stochastic sewing lemma.

\subsection{Malliavin calculus for the Riemann--Liouville process}

We recall the Malliavin calculus associated with the Riemann--Liouville
process. The construction below follows the standard isonormal-Gaussian
framework and integration-by-parts formalism of
\cite{nualart2006malliavin}. Since we work with the fresh
Riemann--Liouville component rather than the full fractional Brownian motion,
the corresponding kernel representation is slightly simpler.

Fix $s\in\R$ and $\bar t\in(s,s+1]$, set $I=[s,\bar t]$, and let $B^s$ be
the Riemann--Liouville process from \Cref{sec:framework} (the Brownian
increment when $H=1/2$). We denote by $\mathcal H_{s,\bar t}$, or simply $\mathcal H$ when $s$ and $\bar t$ are fixed, the
completion of the step functions on $I$ with respect to the inner product
$\langle \mathbf 1_{[s,t]},\mathbf 1_{[s,u]}\rangle_{\mathcal H}
=\EE[(B_t^s)^{(1)}(B_u^s)^{(1)}]$, where $(\cdot)^{(1)}$ denotes the first
coordinate; all coordinates have the same covariance. Writing
\begin{align*}
 \mathfrak K(t,r):=
 \begin{cases}
 c_H(t-r)^{H-\frac12}\mathbf 1_{\{r<t\}},&H\ne\frac12,\\
 \mathbf 1_{\{r<t\}},&H=\frac12,
 \end{cases}
\end{align*}
the map $\mathfrak K^*$ defined on step functions by
$\mathfrak K^*\mathbf 1_{[s,t]}=\mathfrak K(t,\cdot)$ extends to an isometry
from $\mathcal H$ into $L_2(I)$. If $H<1/2$, then, up to a deterministic
constant depending only on $H$, $\mathfrak K^*$ agrees with the right-sided
fractional derivative $D_{\bar t-}^{1/2-H}$ of order $\frac12-H\in(0,\frac12)$.
In particular,
\begin{align}\label{eq:H-fractional-norm}
\begin{split}
 \|h\|_{\mathcal H}\asymp_H
 \|D_{\bar t-}^{1/2-H}h\|_{L_2(I)}.
\end{split}
\end{align}
When $H=1/2$, the canonical Hilbert space is simply
\begin{align}\label{eq:brownian-canonical-space}
 \mathcal H_{s,\bar t}=L_2([s,\bar t]),
 \qquad \mathfrak K^*=\mathrm{Id}.
\end{align}
When $H>1/2$, the definition of $\mathfrak K^*$ and an elementary
integration give, for a step function $h$ on $I$,
\begin{align}\label{eq:smooth-K-star}
\begin{split}
 (\mathfrak K^*h)(u)
 =c_H\Bigl(H-\frac12\Bigr)\int_u^{\bar t}
 h(r)(r-u)^{H-\frac32}\,dr,
 \qquad u\in I.
\end{split}
\end{align}
Thus $\mathfrak K^*$ is, up to a deterministic constant, the right-sided
fractional integral of order $H-\frac12\in(0,\frac12)$.  Young's convolution
inequality gives
\begin{align}\label{eq:smooth-H-embedding}
\begin{split}
 \|h\|_{\mathcal H_{s,\bar t}}
 =\|\mathfrak K^*h\|_{L_2(I)}
 \le C_H(\bar t-s)^{H-\frac12}\|h\|_{L_2(I)}.
\end{split}
\end{align}
Since step functions are dense in $L_2(I)$, this proves the continuous
embedding $L_2(I)\hookrightarrow\mathcal H_{s,\bar t}$.  By tensorisation,
for every $\ell\ge1$ and every
$h\in L_2(I^\ell)$,
\begin{align}\label{eq:smooth-H-tensor-embedding}
\begin{split}
 \|h\|_{\mathcal H_{s,\bar t}^{\otimes\ell}}
 \le C_H^\ell(\bar t-s)^{\ell(H-\frac12)}
 \|h\|_{L_2(I^\ell)}.
\end{split}
\end{align}

For a $d$-dimensional Riemann--Liouville process, $\mathcal H$ is understood
componentwise. The Malliavin derivative with respect to this isonormal
Gaussian process is denoted by $D$; in particular,
$D_\sigma B_t^s=\mathbf 1_{[s,t]}(\sigma)I_d$.
Tensor products are Hilbert tensor products over $\mathcal H$ and Euclidean
tensor products in the state variables; all finite-dimensional tensor norms
are the associated Hilbert--Schmidt norms. For a finite-dimensional or
Hilbert-valued random variable $Z$, $\ell\in\N$, and $p\ge1$, we use the norm
\begin{align*}
 \|Z\|_{\mathbb D^{\ell,p}}
 :=\left(
 \EE|Z|^p
 +\sum_{i=1}^\ell
 \EE\|D^iZ\|_{\mathcal H^{\otimes i}}^p
 \right)^{1/p}.
\end{align*}

Throughout this subsection, $F$ satisfies \Cref{ass:dissipative-drift},
$x\in\R^d$, and $\mathbf v\in C(I;\R^d)$ is deterministic. We study the
regularising flow $\phi_t=\phi_t^{s,x,\mathbf v}$, $t\in I$, from
\eqref{eq:flow-front}. Since $\mathbf v$ is deterministic, $D\mathbf v=0$.

The following interpolation estimate will be used when passing to the limit
in the Malliavin derivative. Its purpose is to turn uniform convergence of
Lipschitz functions without jumps into convergence in $\mathcal H_{s,\bar t}$.

\begin{lemma}\label{lem:Hinterp}
There exists a constant $C=C(H)$ such that the following holds. Let
$s\in\R$, $\bar t\in(s,s+1]$, $I=[s,\bar t]$, and $t\in I$. Then, for every
Lipschitz function $\psi:I\to\R^{d\times d}$ satisfying $\psi=0$ on
$[t,\bar t]$,
\begin{align}\label{eq:Hinterp}
\begin{split}
 \|\psi\|_{\mathcal H_{s,\bar t}}
 \le
 \begin{cases}
 C\Big(
 \|\psi\|_\infty
 +\|\psi\|_\infty^{\frac12+H}
 [\psi]_{\operatorname{Lip}(I)}^{\frac12-H}
 \Big),&H<\frac12,\\[3pt]
 C\|\psi\|_\infty,&H\ge\frac12.
 \end{cases}
\end{split}
\end{align}
\end{lemma}

\begin{proof}
If $H=1/2$, \eqref{eq:brownian-canonical-space} gives directly
$\|\psi\|_{\mathcal H}=\|\psi\|_{L_2(I)}\le(\bar t-s)^{1/2}\|\psi\|_\infty$.
If $H>1/2$, \eqref{eq:smooth-H-embedding} gives
$\|\psi\|_{\mathcal H}\le C_H(\bar t-s)^{H-\frac12}\|\psi\|_{L_2(I)}
\le C_H(\bar t-s)^H\|\psi\|_\infty$.
We may therefore assume that $H<1/2$ and write $a_H:=\frac12-H$.
By \eqref{eq:H-fractional-norm}, it is enough to estimate the right-sided
fractional derivative. For $\sigma\in I$,
\begin{align*}
 D_{\bar t-}^{a_H}\psi(\sigma)
 =\frac1{\Gamma(1-a_H)}\left(
 \frac{\psi(\sigma)}{(\bar t-\sigma)^{a_H}}
 +a_H\int_\sigma^{\bar t}
 \frac{\psi(\sigma)-\psi(\upsilon)}
 {(\upsilon-\sigma)^{1+a_H}}\,d\upsilon
 \right),
\end{align*}
where $\Gamma$ is the Gamma function. Both terms vanish for
$\sigma\in(t,\bar t]$. Since $2a_H<1$, the boundary term satisfies
\begin{align*}
 \left\|\frac{\psi(\cdot)}{(\bar t-\cdot)^{a_H}}
 \right\|_{L_2([s,t])}
 &\le \|\psi\|_\infty
 \left(\int_s^t(\bar t-\sigma)^{-2a_H}\,d\sigma\right)^{1/2} \le C\|\psi\|_\infty.
\end{align*}
For the integral term, put $\theta_0:=\|\psi\|_\infty$ and
$\theta_1:=[\psi]_{\operatorname{Lip}(I)}$. If $\theta_0\theta_1=0$, the
conclusion is immediate. Otherwise, let
$r_0:=\min\{\bar t-s,\theta_0/\theta_1\}$.
Using the Lipschitz bound for $\upsilon-\sigma\le r_0$ and the supremum bound
for $\upsilon-\sigma>r_0$ gives, uniformly in $\sigma$,
\begin{align*}
 \left|\int_\sigma^{\bar t}
 \frac{\psi(\sigma)-\psi(\upsilon)}
 {(\upsilon-\sigma)^{1+a_H}}\,d\upsilon\right|
 \le C\left(\theta_1r_0^{1-a_H}
 +\theta_0r_0^{-a_H}\right).
\end{align*}
If $r_0=\theta_0/\theta_1$, the right-hand side is bounded by
$C\theta_0^{1-a_H}\theta_1^{a_H}$. If $r_0=\bar t-s<\theta_0/\theta_1$, then
$\theta_1r_0^{1-a_H}\le\theta_0r_0^{-a_H}$, and the right-hand side is
bounded by $C\theta_0(\bar t-s)^{-a_H}$; its $L_2([s,t])$ norm is then at most
$C\theta_0(\bar t-s)^{1/2-a_H}=C\theta_0(\bar t-s)^H\le C\theta_0$. Combining
the two estimates yields \eqref{eq:Hinterp}, with a constant depending only
on $H$, since $\bar t-s\le1$.
\end{proof}

The proofs of \Cref{lem:DY} and \Cref{prop:higher-DY} below approximate $F$
by globally Lipschitz vector fields. To pass to the limit in derivatives of
all orders, the approximation must keep the one-sided bound
\eqref{eq:one-sided-lipschitz} and control every derivative that enters.
The Yosida approximation does both.

\begin{lemma}[Regularisation of the drift]\label{lem:yosida}
Set $\mathfrak m:=\kappa_4\,\mathrm{id}-F$, which is monotone by
\eqref{eq:F-one-sided-global}, and, for $N\ge1$,
$\mathcal R_N:=(\mathrm{id}+N^{-1}\mathfrak m)^{-1}$ and
$F_N:=\kappa_4\,\mathrm{id}-\mathfrak m\circ\mathcal R_N$. Then $F_N$ is
smooth, globally Lipschitz with $\|\nabla F_N\|_\infty\le\kappa_4+2N$, and
satisfies \eqref{eq:one-sided-lipschitz} with the same constant $\kappa_4$.
Moreover, there exist constants $C,\mu$ depending only on $k$,
$\boldsymbol\kappa$, $\gamma$, and $d$ such that, for every
$j\in\{0,\ldots,\lfloor-\gamma\rfloor+5\}$,
\begin{align}\label{eq:yosida-derivatives}
\begin{split}
 \nabla^jF_N\longrightarrow\nabla^jF\quad\text{locally uniformly},
 \qquad
 \sup_{N\ge1}|\nabla^jF_N(x)|\le C(1+|x|)^{\mu},
 \quad x\in\R^d.
\end{split}
\end{align}
\end{lemma}

\begin{proof}
The map $\mathrm{id}+N^{-1}\mathfrak m$ is smooth, strongly monotone, and its
Jacobian has a positive definite symmetric part, so it is a smooth
diffeomorphism of $\R^d$ and $\mathcal R_N$ is smooth. If $A+A^\top\ge0$,
then $|(I+A)w|\,|w|\ge\langle(I+A)w,w\rangle\ge|w|^2$, so
$\|(I+A)^{-1}\|\le1$; hence
$\nabla\mathcal R_N=(I+N^{-1}\nabla\mathfrak m(\mathcal R_N))^{-1}$ has norm
at most $1$. Testing $y+N^{-1}\mathfrak m(y)=0$, $y:=\mathcal R_N0$, against
$y$ and using monotonicity gives $|y|\le|F(0)|$, so
$|\mathcal R_Nx|\le|x|+|F(0)|$. Since
$\mathfrak m\circ\mathcal R_N=N(\mathrm{id}-\mathcal R_N)$, we get
$\|\nabla F_N\|_\infty\le\kappa_4+2N$, and, as $\mathfrak m\circ\mathcal R_N$
is monotone (write $x=\mathcal R_Nx+N^{-1}\mathfrak m(\mathcal R_Nx)$),
$v^\top\nabla F_N v\le\kappa_4|v|^2$.

Differentiating $\mathcal R_N+N^{-1}\mathfrak m\circ\mathcal R_N=\mathrm{id}$
$j\ge2$ times with the Fa\`a di Bruno formula and isolating the term
containing $\nabla^j\mathcal R_N$ expresses $\nabla^j\mathcal R_N$ as
$-(I+N^{-1}\nabla\mathfrak m(\mathcal R_N))^{-1}N^{-1}$ applied to a finite
sum of terms $\nabla^i\mathfrak m(\mathcal R_N)$, $2\le i\le j$, evaluated
on derivatives of $\mathcal R_N$ of order less than $j$. By induction on
$j$ and \eqref{eq:growth},
$|\nabla^j\mathcal R_N(x)|\le CN^{-1}(1+|x|)^{\mu}$ for
$2\le j\le\lfloor-\gamma\rfloor+5$. Similarly,
$|\mathcal R_Nx-x|=N^{-1}|\mathfrak m(\mathcal R_Nx)|$ and
$\nabla\mathcal R_N-I
=-(I+N^{-1}\nabla\mathfrak m(\mathcal R_N))^{-1}N^{-1}\nabla\mathfrak m(\mathcal R_N)$
are bounded by $CN^{-1}(1+|x|)^{\mu}$. Finally,
$\nabla^j(\mathfrak m\circ\mathcal R_N)$ is, by the Fa\`a di Bruno formula,
the sum of $\nabla^j\mathfrak m(\mathcal R_N)[\nabla\mathcal R_N,\ldots,\nabla\mathcal R_N]$
and of terms containing at least one factor $\nabla^i\mathcal R_N$ with
$i\ge2$. The bounds above give the uniform polynomial growth in
\eqref{eq:yosida-derivatives}, and, since $\mathcal R_N\to\mathrm{id}$ and
$\nabla\mathcal R_N\to I$ locally uniformly, only the first term survives in
the limit, which gives $\nabla^jF_N\to\nabla^jF$ locally uniformly.
\end{proof}

For $s\le\sigma\le t\le\bar t$, let $J_{\sigma,t}$ denote the linearised flow
along $\phi$, that is, $\partial_tJ_{\sigma,t}=\nabla F(\phi_t)J_{\sigma,t}$
and $J_{\sigma,\sigma}=I_d$. For every $\xi\in\R^d$,
\eqref{eq:one-sided-lipschitz} gives
$\frac{d}{dt}|J_{\sigma,t}\xi|^2
=2\langle\nabla F(\phi_t)J_{\sigma,t}\xi,J_{\sigma,t}\xi\rangle
\le2\kappa_4|J_{\sigma,t}\xi|^2$, and Gr\"onwall's inequality yields
\begin{align}\label{eq:J-Y-bound}
\begin{split}
 \|J_{\sigma,t}\|_{\mathrm{op}}\le e^{\kappa_4(t-\sigma)}.
\end{split}
\end{align}

\begin{lemma}\label{lem:DY}
Let $I=[s,\bar t]$ and $\phi=\phi^{s,x,\mathbf v}$ be as above. For every
$t\in I$ and every $p\in[1,\infty)$, one has
$\phi_t\in\mathbb D^{1,p}$. Moreover, there exist a jointly measurable map
$I^2\times\Omega\ni(\sigma,t,\omega)\mapsto D_\sigma\phi_t(\omega)$
and a set $\Omega_1\in\mathcal F$ with $\PP(\Omega_1)=1$ such that, for every
$\omega\in\Omega_1$ and all $(\sigma,t)\in I^2$,
\begin{align}\label{eq:DY-stmt}
\begin{split}
 D_\sigma\phi_t(\omega)
 =J_{\sigma,t}(\omega)\mathbf 1_{[s,t]}(\sigma)
 \quad\text{and}\quad
 D_\sigma\phi_t
 =\mathbf 1_{[s,t]}(\sigma)I_d
 +\int_s^t\nabla F(\phi_r)D_\sigma\phi_r\,dr.
\end{split}
\end{align}
For every fixed $t\in I$, the map
$(\sigma,\omega)\mapsto D_\sigma\phi_t(\omega)$ is a representative of the
abstract Malliavin derivative $D\phi_t\in L_p(\Omega;\mathcal H_{s,\bar t})$, and
\begin{align}\label{eq:DY-eq}
\begin{split}
 D\phi_t=\mathbf 1_{[s,t]}I_d
 +\int_s^t\nabla F(\phi_r)D\phi_r\,dr.
\end{split}
\end{align}
\end{lemma}

\begin{proof}
We regularise $F$, identify the derivative for the regularised equations,
and then use the closedness of $D$.

\smallskip
\emph{Regularisation.}
Let $F_N$, $N\ge1$, be the vector fields of \Cref{lem:yosida}. Let $\phi^N$ be the solution of \eqref{eq:flow-front}
with $F$ replaced by $F_N$.

Let $\mathfrak y_r:=\mathbf v_{s,r}+B_r^s$ and $\zeta_r^N:=\phi_r^N-\mathfrak y_r$. The
one-sided Lipschitz property between $\zeta_r^N+\mathfrak y_r$ and $\mathfrak y_r$ gives
\begin{align*}
 \frac{d}{dr}|\zeta_r^N|^2
 =2\langle \zeta_r^N,F_N(\zeta_r^N+\mathfrak y_r)\rangle
 \le C|\zeta_r^N|^2+C|F_N(\mathfrak y_r)|^2.
\end{align*}
Gr\"onwall's inequality, the uniform polynomial bound on $F_N$ from \eqref{eq:yosida-derivatives}, and the
Gaussian moments of the supremum $\sup_{r\in I}|B_r^s|$ yield
\begin{align}\label{eq:YN-moments}
\begin{split}
 \sup_{N\ge1}\EE\left[\sup_{r\in I}|\phi_r^N|^p\right]
 +\EE\left[\sup_{r\in I}|\phi_r|^p\right]<\infty,
 \qquad p<\infty,
\end{split}
\end{align}
the same argument applying to $\phi$.

\smallskip
\emph{Derivative of the regularised equation.}
Since $F_N$ is globally Lipschitz with bounded first derivative, the
Malliavin chain rule and differentiation under the time integral give
\begin{align*}
 D_\sigma \phi_t^N
 =\mathbf 1_{[s,t]}(\sigma)I_d
 +\int_s^t\nabla F_N(\phi_r^N)D_\sigma \phi_r^N\,dr.
\end{align*}
Adaptedness gives $D_\sigma \phi_r^N=0$ for $r<\sigma$. Therefore, uniqueness
for the linear equation yields
$D_\sigma \phi_t^N=J^N_{\sigma,t}\mathbf 1_{[s,t]}(\sigma)$, where
$\partial_tJ^N_{\sigma,t}=\nabla F_N(\phi_t^N)J^N_{\sigma,t}$ and
$J^N_{\sigma,\sigma}=I_d$. Exactly as in \eqref{eq:J-Y-bound},
\begin{align}\label{eq:JN-uniform-bound}
\begin{split}
 \|J^N_{\sigma,t}\|_{\mathrm{op}}
 \le e^{\kappa_4(t-\sigma)},
\end{split}
\end{align}
uniformly in $N$.

\smallskip
\emph{Passage to the limit.}
Subtracting the equations for $\phi^N$ and $\phi$, using the one-sided
Lipschitz bound for $F_N$, and applying Gr\"onwall's inequality gives
\begin{align}\label{eq:YN-state-convergence}
\begin{split}
 \sup_{t\in I}|\phi_t^N-\phi_t|
 \le e^{\kappa_4(\bar t-s)}\int_s^{\bar t}|F_N(\phi_r)-F(\phi_r)|\,dr.
\end{split}
\end{align}
The local convergence \eqref{eq:yosida-derivatives}, the polynomial bounds,
and \eqref{eq:YN-moments} imply that the right-hand side converges to zero
almost surely and in every $L_p$. Thus
\begin{align}\label{eq:YN-Lp-convergence}
\begin{split}
 \phi_t^N\longrightarrow \phi_t
 \quad\text{in }L_p(\Omega),
 \qquad p<\infty.
\end{split}
\end{align}
The coefficients in the linearised equations converge in $L_1(I)$ almost
surely, that is,
$\int_s^{\bar t}\|\nabla F_N(\phi_r^N)-\nabla F(\phi_r)\|\,dr\to0$.
Indeed, this follows from \eqref{eq:yosida-derivatives} on the random
compact set containing the paths $\phi^N$ and $\phi$, together with
\eqref{eq:YN-state-convergence}. Stability of linear ODEs and
\eqref{eq:JN-uniform-bound} therefore yield
\begin{align}\label{eq:JN-convergence}
\begin{split}
 \sup_{s\le\sigma\le t\le\bar t}
 \|J^N_{\sigma,t}-J_{\sigma,t}\|
 \longrightarrow0
 \qquad\text{almost surely}.
\end{split}
\end{align}

Fix $t\in I$ and set $h_t^N(\sigma):=J^N_{\sigma,t}\mathbf 1_{[s,t]}(\sigma)$
and $h_t(\sigma):=J_{\sigma,t}\mathbf 1_{[s,t]}(\sigma)$.
The jump at $t$ cancels in $h_t^N-h_t$, since
$J^N_{t,t}=J_{t,t}=I_d$. To apply \Cref{lem:Hinterp}, we estimate its
Lipschitz seminorm. The cocycle property gives, for
$s\le\sigma_1<\sigma_2\le t$,
\begin{align*}
 J^N_{\sigma_1,t}-J^N_{\sigma_2,t}
 =J^N_{\sigma_2,t}
 \int_{\sigma_1}^{\sigma_2}
 \nabla F_N(\phi_r^N)J^N_{\sigma_1,r}\,dr.
\end{align*}
Using \eqref{eq:JN-uniform-bound}, the uniform polynomial bound for
$\nabla F_N$, and \eqref{eq:YN-moments}, we obtain
$[J^N_{\cdot,t}]_{\operatorname{Lip}([s,t])}
\le C(1+\sup_{r\in I}|\phi_r^N|^{\varpi})$ for a deterministic exponent
$\varpi\ge0$, and the analogous estimate holds
for $J$. Consequently, \Cref{lem:Hinterp}, \eqref{eq:JN-convergence}, and
Vitali's convergence theorem, which applies because these Lipschitz
constants have moments bounded uniformly in $N$, imply
\begin{align}\label{eq:DYN-H-convergence}
\begin{split}
 h_t^N\longrightarrow h_t
 \quad\text{in }L_p(\Omega;\mathcal H),
 \qquad p<\infty.
\end{split}
\end{align}
Here $h_t\in L_p(\Omega;\mathcal H)$ because
$h_t=\mathbf 1_{[s,t]}I_d+(J_{\cdot,t}-I_d)\mathbf 1_{[s,t]}$, and the second
term has no jump at $t$ and is controlled by \Cref{lem:Hinterp}.

Combining \eqref{eq:YN-Lp-convergence} and
\eqref{eq:DYN-H-convergence}, the closedness of the Malliavin derivative gives
$\phi_t\in\mathbb D^{1,p}$ and $D\phi_t=h_t$ for every fixed $t$. Choose a set
$\Omega_1$ of probability one on which $\phi$ has continuous trajectories and
the random linear equation for $J$ is well defined for
all $s\le\sigma\le t\le\bar t$. Parameter-dependent ODE theory, or equivalently
the Picard expansion of $J_{\sigma,t}$, shows that
$(\sigma,t,\omega)\mapsto J_{\sigma,t}(\omega)$ is jointly measurable.
Define $D_\sigma\phi_t$ by the first identity in \eqref{eq:DY-stmt} on
$\Omega_1$ and set it equal to zero on the complement. For every fixed $t$,
the identification $D\phi_t=h_t$ shows that this kernel represents
the abstract Malliavin derivative. Since it is defined through the linearised
flow, both identities in \eqref{eq:DY-stmt} hold on the same set
$\Omega_1$, simultaneously for all $(\sigma,t)$. The Hilbert-valued identity
\eqref{eq:DY-eq} follows immediately.
\end{proof}

For $\mathbf u=(u_1,\ldots,u_\ell)$, we write $\bar u:=\max_{1\le i\le\ell}u_i$
and, for $A\subset\{1,\ldots,\ell\}$, $\mathbf u_A:=(u_i)_{i\in A}$. We denote
by $\mathcal P_\ell$ the set of partitions of $\{1,\ldots,\ell\}$.

\begin{proposition}[Higher-order Malliavin derivatives]
\label{prop:higher-DY}
Let $I=[s,\bar t]$ and $\phi=\phi^{s,x,\mathbf v}$ be as above. One has
$\phi_t\in\mathbb D^{\ell,p}$ for every
$\ell\in\{1,\ldots,\lfloor-\gamma\rfloor+5\}$, $p<\infty$, and $t\in I$.
Moreover, there exist jointly measurable maps
$I^\ell\times I\times\Omega\ni(\mathbf u,t,\omega)
\mapsto D^\ell_{\mathbf u}\phi_t(\omega)$, $1\le\ell\le\lfloor-\gamma\rfloor+5$,
and a set $\Omega_*\in\mathcal F$ with $\PP(\Omega_*)=1$ such that the
following holds for every $\omega\in\Omega_*$, every
$1\le\ell\le\lfloor-\gamma\rfloor+5$, every $\mathbf u\in I^\ell$, and every
$t\in I$. The kernel vanishes whenever $t<\bar u$, and, for $\ell\ge2$ and
$t\ge\bar u$,
\begin{align}\label{eq:higher-DY-partitions}
\begin{split}
 D^\ell_{\mathbf u}\phi_t
 =\int_{\bar u}^tJ_{r,t}
 \sum_{\substack{\pi\in\mathcal P_\ell\\|\pi|\ge2}}
 \nabla^{|\pi|}F(\phi_r)
 \left[
 D^{|A|}_{\mathbf u_A}\phi_r:A\in\pi
 \right]dr,
\end{split}
\end{align}
where the arguments indexed by the blocks $A\in\pi$ are inserted into the
multilinear map $\nabla^{|\pi|}F(\phi_r)$ in any fixed order. For every fixed
$t$, this kernel is a representative of
$D^\ell \phi_t\in L_p(\Omega;\mathcal H_{s,\bar t}^{\otimes\ell})$.
\end{proposition}

\begin{proof}
The case $\ell=1$ is \Cref{lem:DY}; let
$2\le\ell\le\lfloor-\gamma\rfloor+5$. As for $\ell=1$, we approximate $F$ by
the vector fields $F_N$ of \Cref{lem:yosida}. Since $F_N$ is globally
Lipschitz but its higher derivatives need not be bounded, we first add a
cutoff, then remove it for fixed $N$, and finally let $N\to\infty$.

\smallskip
\emph{A priori bound.}
Let $G\in C^\infty(\R^d;\R^d)$ be globally Lipschitz with
$|\nabla^jG(x)|\le C(1+|x|)^{\mu}$ for $j\le\ell$, let $\phi^G$ be the
solution of \eqref{eq:flow-front} with $F$ replaced by $G$, and let $J^G$ be
its linearised flow, with $\|J^G_{\sigma,t}\|\le\Theta$ for
$s\le\sigma\le t\le\bar t$. Assume that $\phi^G_t\in\mathbb D^{\ell,p}$ for
all $t$ and $p$, and that, for $2\le i\le\ell$,
\begin{align}\label{eq:higher-DY-partitions-H}
\begin{split}
 D^i\phi_t^G
 =\int_s^tJ^G_{r,t}
 \sum_{\substack{\pi\in\mathcal P_i\\|\pi|\ge2}}
 \nabla^{|\pi|}G(\phi_r^G)
 \bigl[D^{|A|}\phi_r^G:A\in\pi\bigr]dr
\end{split}
\end{align}
as a Bochner integral in $\mathcal H^{\otimes i}$, where $D^{|A|}\phi_r^G$
acts on the Malliavin variables indexed by $A$. Since $D\phi^G_r
=\mathbf 1_{[s,r]}I_d+(J^G_{\cdot,r}-I_d)\mathbf 1_{[s,r]}$, the cocycle
identity gives $[J^G_{\cdot,r}]_{\operatorname{Lip}([s,r])}
\le\Theta^2\sup_I|\nabla G(\phi^G)|$, and \Cref{lem:Hinterp} bounds
$\|D\phi_r^G\|_{\mathcal H}$. The $\mathcal H^{\otimes i}$-norm of a tensor
product of kernels in disjoint groups of variables is the product of their
norms, so induction on $i$ in \eqref{eq:higher-DY-partitions-H} yields
\begin{align}\label{eq:higher-DY-apriori}
\begin{split}
 \sup_{t\in I}\|D^i\phi_t^G\|_{\mathcal H^{\otimes i}}
 \le C(1+\Theta)^{3i}\bigl(1+\sup_{r\in I}|\phi_r^G|\bigr)^{\mu_i},
 \qquad 1\le i\le\ell,
\end{split}
\end{align}
with $C,\mu_i$ depending only on $\ell$, $H$, and the growth constants of
$G$. The argument only uses the recursion, so the same bound holds for any
family of kernels defined recursively by the right-hand side of
\eqref{eq:higher-DY-partitions-H}, starting from $J^G_{\cdot,t}\mathbf
1_{[s,t]}$.

\smallskip
\emph{Step 1: cutoff.}
Fix a smooth nondecreasing $\mathfrak h:[0,\infty)\to[0,2]$ with
$\mathfrak h(r)=r$ on $[0,1]$, $\mathfrak h\equiv2$ on $[3,\infty)$, and
$\mathfrak h'\le1$, and set $\mathfrak c_R(x):=R\,\mathfrak h(|x|/R)\,x/|x|$
for $R\ge1$. The map $\mathfrak c_R$ is smooth, equals the identity on the
ball $B_R$, satisfies $|\mathfrak c_R(x)|\le|x|\wedge2R$, is $1$-Lipschitz
(its Jacobian is symmetric with eigenvalues in $[0,1]$), and
$\|\nabla^i\mathfrak c_R\|_\infty\le C_iR^{1-i}$. Hence
$G_{N,R}:=F_N\circ\mathfrak c_R$ has bounded derivatives of every order,
$\|\nabla G_{N,R}\|_\infty\le\kappa_4+2N$, and, by \Cref{lem:yosida} and the
Fa\`a di Bruno formula, $|\nabla^jG_{N,R}(x)|\le C(1+|x|)^{\mu}$ for
$j\le\ell$, uniformly in $N$ and $R\ge1$. Let $\phi^{N,R}$ and $J^{N,R}$ be
the corresponding flow and linearised flow. For a drift with bounded
derivatives of all orders, the classical argument for smooth equations
(Picard iteration, the chain rule, and closedness of $D$; see
\cite[Section~2.2]{nualart2006malliavin}) applies verbatim, the noise
entering additively. It gives $\phi_t^{N,R}\in\mathbb D^{\ell,p}$ for every
$p<\infty$, and $D^\ell\phi_t^{N,R}$ is obtained by differentiating the
equation $\ell$ times. Removing the one-block partition by variation of
constants gives \eqref{eq:higher-DY-partitions-H} for $G=G_{N,R}$.

\smallskip
\emph{Step 2: removal of the cutoff for fixed $N$.}
Let $k^{N,1}_t:=D\phi_t^N$ (\Cref{lem:DY}) and define $k_t^{N,i}$, $2\le
i\le\ell$, recursively by the right-hand side of
\eqref{eq:higher-DY-partitions-H} with $(G,\phi^G,J^G)$ replaced by
$(F_N,\phi^N,J^N)$ and $D^{|A|}\phi_r^G$ by $k_r^{N,|A|}$. On the event
$\Omega_R:=\{\sup_I|\phi^N|\le R\}$, uniqueness gives $\phi^{N,R}=\phi^N$
and $J^{N,R}=J^N$ on $I$, so, by induction on $i$,
$D^i\phi_t^{N,R}=k_t^{N,i}$ on $\Omega_R$. By \eqref{eq:higher-DY-apriori}
with $\Theta=e^{\kappa_4+2N}$, and since $\sup_I|\phi^{N,R}|$ has moments of
every order uniformly in $R$ (the drift $G_{N,R}$ has linear growth,
uniformly in $R$), $\phi_t^{N,R}$ is bounded in $\mathbb D^{\ell,p}$
uniformly in $R$, and $k_t^{N,i}$ satisfies the same bound. As
$\PP(\Omega_R)\to1$, H\"older's inequality gives $\phi_t^{N,R}\to\phi_t^N$
and $D^i\phi_t^{N,R}\to k_t^{N,i}$ in $L_p$ for every $p<\infty$. The
closedness of $D^i$ gives $\phi_t^N\in\mathbb D^{\ell,p}$ and
$D^i\phi_t^N=k_t^{N,i}$, that is, \eqref{eq:higher-DY-partitions-H} holds for
$G=F_N$.

\smallskip
\emph{Step 3: the limit $N\to\infty$.}
Let $k_t^1:=D\phi_t$ and define $k_t^i$, $2\le i\le\ell$, recursively as in
Step 2 with $(F,\phi,J)$ in place of $(F_N,\phi^N,J^N)$. By
\eqref{eq:JN-uniform-bound} the a priori bound applies to $G=F_N$ with
$\Theta=e^{\kappa_4}$, independently of $N$; this is where the one-sided
bound in \Cref{lem:yosida} is needed. Together with \eqref{eq:YN-moments},
this shows that $D^i\phi_t^N$ is bounded in
$L_p(\Omega;\mathcal H^{\otimes i})$ uniformly in $N$ and $t$, and the same
argument applies to $k_t^i$. We show by induction on $i$ that
$D^i\phi_t^N\to k_t^i$ in $L_p(\Omega;\mathcal H^{\otimes i})$ for every $t$
and $p<\infty$. For $i=1$ this is \eqref{eq:DYN-H-convergence}. For $i\ge2$,
the difference $D^i\phi_t^N-k_t^i$ is a sum of integrals over $[s,t]$ whose
integrands contain one of the factors $J^N_{r,t}-J_{r,t}$,
$\nabla^{j}F_N(\phi_r^N)-\nabla^{j}F(\phi_r)$, or
$D^{|A|}\phi_r^N-k_r^{|A|}$, the other factors being bounded in every $L_q$.
The first two factors converge to zero almost surely, uniformly in $r$, by
\eqref{eq:JN-convergence}, \eqref{eq:YN-state-convergence}, and
\eqref{eq:yosida-derivatives}, and they are dominated by
$C(1+\sup_I|\phi^N|+\sup_I|\phi|)^{\mu}$, whose moments are bounded uniformly
in $N$. The third converges for every $r$ by the induction hypothesis, with
moments bounded uniformly in $r$. H\"older's inequality and Vitali's
convergence theorem (in $\omega$ and in $r$), based on these uniform moment
bounds, give the claim. Since
$\phi_t^N\to\phi_t$ in every $L_p$ by \eqref{eq:YN-Lp-convergence}, the
closedness of $D^\ell$ gives $\phi_t\in\mathbb D^{\ell,p}$ and
$D^\ell\phi_t=k_t^\ell$.

\smallskip
\emph{Pointwise kernels.}
Define the pointwise kernels recursively by \eqref{eq:higher-DY-partitions},
starting from \eqref{eq:DY-stmt}. Their joint measurability follows from
measurability of parameter-dependent Lebesgue integrals, they vanish for
$t<\bar u$ because $D_\sigma\phi_r=0$ for $r<\sigma$, and, for every fixed
$t$, they represent $k_t^\ell=D^\ell\phi_t$, since the $\mathcal
H^{\otimes\ell}$-valued integral in \eqref{eq:higher-DY-partitions-H} is
represented by the pointwise integral of the kernels. Intersecting the
finitely many full-probability sets for $1\le\ell\le\lfloor-\gamma\rfloor+5$
gives the common set $\Omega_*$.
\end{proof}

\subsection{Malliavin integration-by-parts estimates}\label{sec:Malliavin-IBP}

The preceding construction identifies the relevant Malliavin derivatives and
Hilbert-space norms. We now turn these objects into quantitative
integration-by-parts estimates at the critical regularity. For negative
spatial Besov--H\"older regularity, the number of Malliavin integrations by
parts is chosen as a larger integer and is used only to make the
high-frequency series summable. At non-negative regularity, the one-point
estimate follows directly from boundedness, while the same dyadic argument
yields the two-point and four-point estimates needed below.

Let $Z=(Z^1,\ldots,Z^d)$ be a random vector with Malliavin covariance matrix
$\mathcal M_Z^{i,i'}:=\langle DZ^i,DZ^{i'}\rangle_{\mathcal H}$,
$1\le i,i'\le d$. Matrix norms in this subsection are operator norms, and
norms of $\mathcal H$-valued matrices are Hilbert--Schmidt norms in the
Euclidean indices. We denote the divergence operator by $\boldsymbol\delta$.
For a real-valued random variable $G$ and $i\in\{1,\ldots,d\}$, define
\begin{align*}
 \mathscr H_i(Z,G)
 :=\boldsymbol\delta\left(
 G\sum_{i'=1}^d(\mathcal M_Z^{-1})_{i,i'}DZ^{i'}
 \right),
\end{align*}
whenever the right-hand side belongs to the domain of $\boldsymbol\delta$.
For a multi-index $\beta=(\beta_1,\ldots,\beta_n)$, define recursively
\begin{align}\label{eq:IBP-weight-multi}
\begin{split}
 \mathscr H_\beta(Z,G)
 :=\mathscr H_{\beta_n}
 \big(Z,\mathscr H_{(\beta_1,\ldots,\beta_{n-1})}(Z,G)\big),
 \qquad
 \mathscr H_\varnothing(Z,G):=G.
\end{split}
\end{align}
Then, for every smooth bounded $\varphi:\R^d\to\R$ with bounded derivatives,
whenever the weights in \eqref{eq:IBP-weight-multi} are well defined (see
\Cref{lem:IBP-weight-bound} below), Malliavin integration by parts gives
\begin{align}\label{eq:Malliavin-IBP-formula}
\begin{split}
 \EE\big[\partial_\beta\varphi(Z)G\big]
 =\EE\big[\varphi(Z)\mathscr H_\beta(Z,G)\big].
\end{split}
\end{align}
This is \cite[Proposition 2.1.4]{nualart2006malliavin}, which is stated for
$Z$ and $G$ in $\mathbb D^\infty$. Its proof only uses the identity
$\partial_i\varphi(Z)=\sum_{i'}(\mathcal M_Z^{-1})_{i,i'}
\langle D\varphi(Z),DZ^{i'}\rangle_{\mathcal H}$, which follows from the
chain rule for $Z\in\mathbb D^{1,p}$, and, at each of the $|\beta|$ steps,
the duality relation $\EE[\langle D\Phi,u\rangle_{\mathcal H}]
=\EE[\Phi\,\boldsymbol\delta(u)]$ for $\Phi=\partial_{\beta'}\varphi(Z)$, $\beta'$ a
sub-multi-index of $\beta$, and $u$ in the domain of $\boldsymbol\delta$. It
therefore applies under the finite regularity and integrability provided by
\Cref{lem:IBP-weight-bound}.

We next record the quantitative estimate used below. For each $n\in\N$, fix
an integer $\mathsf N_n=\mathsf N_n(n,d)$ large enough for the product
estimates generated by $n$ Malliavin differentiations, and set
\begin{align}\label{eq:Lambda-def}
\begin{split}
 \Lambda_{n,q}(Z)
 :=1+\sum_{a=1}^{\mathsf N_n}
 \|\mathcal M_Z^{-1}\|_{L_{q\mathsf N_n}}^a
 \|DZ\|_{\mathbb D^{n,q\mathsf N_n}(\mathcal H)}^{2a-1}.
\end{split}
\end{align}
The value of $\mathsf N_n$ plays no role, but it is fixed once $n$ and $d$
are fixed, and we choose it nondecreasing in $n$. Then $\Lambda_{n,q}(Z)$ is
nondecreasing in $n$ and in $q$. We set $\Lambda_{n,q}(Z)=+\infty$ if $Z$ does not have the
required Malliavin regularity, if $\mathcal M_Z$ is not invertible almost
surely, or if one of the norms in \eqref{eq:Lambda-def} is infinite.

\begin{lemma}[Estimate for the integration-by-parts weights]
\label{lem:IBP-weight-bound}
Let $n_0\in\N\cup\{0\}$, $p>1$, and let $\beta$ be a non-empty multi-index.
Set $n:=n_0+|\beta|$.
There exist $q>1$, $\bar q>1$, and a constant $C$, depending only on
$n_0,p,|\beta|$, and $d$, such that, whenever
$\Lambda_{n,q}(Z)<\infty$ and
$G\in\mathbb D^{n,\bar q}$, all the weights in
\eqref{eq:IBP-weight-multi} are well defined in the required divergence
domains and
\begin{align}\label{eq:IBP-weight-bound}
\begin{split}
 \|\mathscr H_\beta(Z,G)\|_{\mathbb D^{n_0,p}}
 \le C\|G\|_{\mathbb D^{n,\bar q}}
 \Lambda_{n,q}(Z)^{|\beta|}.
\end{split}
\end{align}
\end{lemma}

\begin{proof}
For $|\beta|=1$, the continuity of the divergence operator gives
\begin{align*}
 \|\mathscr H_i(Z,G)\|_{\mathbb D^{n_0,p}}
 \lesssim
 \left\|G\sum_{i'=1}^d
 (\mathcal M_Z^{-1})_{i,i'}DZ^{i'}
 \right\|_{\mathbb D^{n_0+1,p}(\mathcal H)},
\end{align*}
for each $i\in\{1,\ldots,d\}$.
Malliavin--Sobolev product estimates reduce the right-hand side to norms of
$G$, $DZ$, and Malliavin derivatives of $\mathcal M_Z^{-1}$. The identities
\begin{align*}
 D\mathcal M_Z^{i,i'}
 =\langle D^2Z^i,DZ^{i'}\rangle_{\mathcal H}
 +\langle DZ^i,D^2Z^{i'}\rangle_{\mathcal H},
 \qquad
 D\mathcal M_Z^{-1}
 =-\mathcal M_Z^{-1}(D\mathcal M_Z)\mathcal M_Z^{-1}
\end{align*}
show, after iteration, that each derivative of $\mathcal M_Z^{-1}$ is a finite
sum of products containing powers of $\mathcal M_Z^{-1}$ and derivatives of
$Z$. H\"older's inequality and the choice of $\mathsf N_n$ in
\eqref{eq:Lambda-def} show both that the vector field inside
$\boldsymbol\delta$ belongs to
$\mathbb D^{n_0+1,p}(\mathcal H)\subset\operatorname{Dom}\boldsymbol\delta$
and that \eqref{eq:IBP-weight-bound} holds for $|\beta|=1$. In other words,
for every $n_0\ge0$ and $p>1$ there exist $q,\bar q>1$ and $C$, depending
only on $n_0,p,d$, such that, for $i\in\{1,\ldots,d\}$,
\begin{align}\label{eq:IBP-one-step}
\begin{split}
 \|\mathscr H_i(Z,G)\|_{\mathbb D^{n_0,p}}
 \le C\|G\|_{\mathbb D^{n_0+1,\bar q}}\Lambda_{n_0+1,q}(Z).
\end{split}
\end{align}
For a general multi-index, write $\beta=(\beta',\beta_{|\beta|})$, so that
$\mathscr H_\beta(Z,G)=\mathscr H_{\beta_{|\beta|}}(Z,\mathscr H_{\beta'}(Z,G))$
by \eqref{eq:IBP-weight-multi}. Applying \eqref{eq:IBP-one-step} with
$G$ replaced by $\mathscr H_{\beta'}(Z,G)$ requires
$\mathscr H_{\beta'}(Z,G)\in\mathbb D^{n_0+1,\bar q}$, which is the induction
hypothesis with $(n_0,p)$ replaced by $(n_0+1,\bar q)$. It gives
\begin{align*}
 \|\mathscr H_\beta(Z,G)\|_{\mathbb D^{n_0,p}}
 \le C\|G\|_{\mathbb D^{n_0+|\beta|,\bar q'}}
 \Lambda_{n_0+|\beta|,q'}(Z)^{|\beta|-1}\Lambda_{n_0+1,q}(Z)
\end{align*}
for suitable $q',\bar q'>1$, and in particular shows that the last
divergence is well defined. Monotonicity of $\Lambda_{n,q}(Z)$ in $n$ and $q$
gives \eqref{eq:IBP-weight-bound} with $n=n_0+|\beta|$, after renaming the
integrability exponents ($\bar q:=\bar q'$ and $q:=q\vee q'$).
\end{proof}

We use an inhomogeneous Littlewood--Paley decomposition
$(\Delta_j)_{j\ge-1}$. For negative orders, the heat-semigroup definition
\eqref{eq:alpha-norm} is equivalent to the standard dyadic definition;
the same holds for the classical H\"older norm at orders in $(0,1)$.
At order zero we keep the convention $\mathcal C^0=L_\infty$, which is
not identified with $B^0_{\infty,\infty}$; boundedness of the dyadic
convolution operators on $L_\infty$ nevertheless gives the same block
bound. Thus, for every $a<1$ in the range considered in this paper,
\begin{align}\label{eq:endpoint-dyadic-holder}
\begin{split}
 \|\Delta_j f\|_\infty
 \lesssim 2^{-aj}\|f\|_{\mathcal C^a},
 \qquad j\ge-1.
\end{split}
\end{align}
We shall also use the following elementary dyadic lifting.

\begin{lemma}[Dyadic primitives]\label{lem:dyadic-primitives}
Let $n\in\N$. For every $j\ge0$ and every smooth function $h$ whose Fourier
transform is supported in a dyadic annulus of size $2^j$, there exist smooth
bounded functions $(h_\alpha)_{|\alpha|=n}$ such that
\begin{align}\label{eq:dyadic-primitives}
\begin{split}
 h=\sum_{|\alpha|=n}\partial^\alpha h_\alpha,
 \qquad
 \|h_\alpha\|_\infty\lesssim_n2^{-nj}\|h\|_\infty.
\end{split}
\end{align}
The same statement holds componentwise for vector-valued functions.
\end{lemma}

\begin{proof}
Let $\chi\in C_c^\infty(\R^d\setminus\{0\})$ equal one on the
reference annulus containing the Fourier support of  $h$.
For each multi-index $\alpha$ of length $n$, the multiplier
$\mu_\alpha(\xi):=\chi(\xi)\frac{n!}{\alpha!}(-i\xi)^\alpha|\xi|^{-2n}$ is
smooth and compactly supported, and the multinomial formula shows that
\begin{align*}
 \sum_{|\alpha|=n}(i\xi)^\alpha \mu_\alpha(\xi)
 =\frac{\chi(\xi)}{|\xi|^{2n}}
 \sum_{|\alpha|=n}\frac{n!}{\alpha!}\xi^{2\alpha}
 =\chi(\xi).
\end{align*}
Let $\check\mu_\alpha:=\mathcal F^{-1}\mu_\alpha$,
$\check\mu_{\alpha,j}(x):=2^{jd}\check\mu_\alpha(2^jx)$, and
$h_\alpha:=2^{-nj}\check\mu_{\alpha,j}*h$.
On the Fourier support of $h$, the multiplier of
$\sum_{|\alpha|=n}\partial^\alpha h_\alpha$ is
$\sum_{|\alpha|=n}2^{-nj}(i\xi)^\alpha \mu_\alpha(2^{-j}\xi)
=\chi(2^{-j}\xi)=1$.
This proves the representation. Since $\check\mu_\alpha$ is Schwartz,
$\|\check\mu_{\alpha,j}\|_{L_1}=\|\check\mu_\alpha\|_{L_1}<\infty$, and hence
$\|h_\alpha\|_\infty\le2^{-nj}\|\check\mu_\alpha\|_{L_1}\|h\|_\infty$.
Derivatives of the kernels are also integrable, so the $h_\alpha$ are
smooth. The constants depend only on $n$, the dimension, and the fixed
annular cutoffs. Applying the construction componentwise proves the
vector-valued assertion.
\end{proof}


\begin{lemma}[Integration by parts at negative regularity]
\label{lem:real-order-IBP}
Let $\sigma>0$ and let $\ell\in\N$ satisfy $\ell>\sigma+1$. There exist
$q>1$, $\bar q>1$, and $C\ge0$, depending only on
$\sigma,\ell$, and $d$, such that, for every
$f\in C_b^\infty(\R^d;\R^d)$ and all $\R^d$-valued random variables
$G,Z_1,Z_2$,
\begin{align}\label{eq:odd-one-point-IBP}
\begin{split}
 \big|\EE\langle G,f(Z_1)\rangle\big|
 &\le C\|f\|_{\mathcal C^{-\sigma}}
 \|G\|_{\mathbb D^{\ell,\bar q}}
 \Lambda_{\ell,q}(Z_1)^\sigma,
\end{split}
\end{align}
and, with $Z_\theta:=\theta Z_1+(1-\theta)Z_2$ for $\theta\in[0,1]$,
\begin{align}
\label{eq:odd-two-point-IBP}
 \big|\EE\langle G,f(Z_1)-f(Z_2)\rangle\big| &\le C\|f\|_{\mathcal C^{-\sigma}}
 \|G\|_{\mathbb D^{\ell,\bar q}}
 \|Z_1-Z_2\|_{\mathbb D^{\ell,\bar q}}
 \sup_{\theta\in[0,1]}
 \Lambda_{\ell,q}(Z_\theta)^{\sigma+1}.
\end{align}
\end{lemma}

\begin{proof}
If one of the quantities on the right-hand side is infinite, the corresponding
estimate is immediate. We therefore assume throughout the proof that all of
them are finite. To justify the dyadic expansions without any convergence
issue, we apply the estimates below first to the finite sums
$\sum_{j=-1}^{N}\Delta_jf$; since $f\in C_b^\infty$, these sums converge to
$f$ together with all derivatives required, the estimates are uniform in $N$,
and passage to the limit proves the stated bounds.

We first prove \eqref{eq:odd-one-point-IBP}. Set
$\mathfrak L_1:=\Lambda_{\ell,q}(Z_1)\ge1$ and  $\lambda_1:=\mathfrak L_1^{-1}\in(0,1]$. Write
$f=\sum_{j\ge-1}\Delta_jf$ and split the sum at $j_1$ chosen so that $2^{-j_1-1}<\lambda_1\le2^{-j_1}$. For the low
frequencies, H\"older's inequality and \eqref{eq:endpoint-dyadic-holder} give
\begin{align}\label{eq:endpoint-one-low}
\begin{split}
 \sum_{j=-1}^{j_1}
 \big|\EE\langle G,\Delta_jf(Z_1)\rangle\big|
 &\le \|G\|_{L_1}
 \sum_{j=-1}^{j_1}\|\Delta_jf\|_\infty
 \lesssim
 2^{\sigma j_1}
 \|f\|_{\mathcal C^{-\sigma}}\|G\|_{L_1}
 \lesssim
 \mathfrak L_1^\sigma
 \|f\|_{\mathcal C^{-\sigma}}
 \|G\|_{\mathbb D^{\ell,\bar q}},
\end{split}
\end{align}
where the geometric summation uses $\sigma>0$.

Let now $j>j_1$. Applying \Cref{lem:dyadic-primitives} with order
$\ell$ to each component of $\Delta_jf$ gives
\begin{align*}
 \Delta_jf_i
 =\sum_{|\alpha|=\ell}\partial^\alpha \psi_{j,i,\alpha},
 \qquad
 \|\psi_{j,i,\alpha}\|_\infty
 \lesssim 2^{-(\ell-\sigma)j}
 \|f\|_{\mathcal C^{-\sigma}},
\end{align*}
since \eqref{eq:dyadic-primitives} gives $2^{-\ell j}$, while
\eqref{eq:endpoint-dyadic-holder} contributes $2^{\sigma j}$. By
\eqref{eq:Malliavin-IBP-formula}, after choosing $q,\bar q$ large enough and
applying \Cref{lem:IBP-weight-bound} with $n_0=0$ and $|\beta|=\ell$,
\begin{align*}
 \big|\EE[G_i\partial^\alpha \psi_{j,i,\alpha}(Z_1)]\big|
 &\lesssim
 2^{-(\ell-\sigma)j}\|f\|_{\mathcal C^{-\sigma}}
 \|G\|_{\mathbb D^{\ell,\bar q}}
 \mathfrak L_1^\ell.
\end{align*}
Since $\ell-\sigma>0$, summing over $j>j_1$ 
yields
\begin{align}
\label{eq:endpoint-one-high}
 \sum_{j>j_1}
 \big|\EE\langle G,\Delta_jf(Z_1)\rangle\big|
 &\lesssim
 \lambda_1^{\ell-\sigma}\|f\|_{\mathcal C^{-\sigma}}
 \|G\|_{\mathbb D^{\ell,\bar q}}\mathfrak L_1^\ell =\|f\|_{\mathcal C^{-\sigma}}
 \|G\|_{\mathbb D^{\ell,\bar q}}\mathfrak L_1^\sigma.
\end{align}
Combining \eqref{eq:endpoint-one-low} and
\eqref{eq:endpoint-one-high} proves \eqref{eq:odd-one-point-IBP}.

For \eqref{eq:odd-two-point-IBP}, we set
$\mathfrak L_2:=\sup_{\theta\in[0,1]}\Lambda_{\ell,q}(Z_\theta)$,
$\lambda_2:=\mathfrak L_2^{-1}$, and split the frequency sum at
$j_2$ chosen so that $2^{-j_2-1}<\lambda_2\le2^{-j_2}$.
For every $j\ge-1$, the fundamental theorem of calculus gives
\begin{align}\label{eq:endpoint-two-ftc}
\begin{split}
 \Delta_jf_i(Z_1)-\Delta_jf_i(Z_2)
 =\sum_{i'=1}^d\int_0^1
 \partial_{i'}\Delta_jf_i(Z_\theta)(Z_1^{i'}-Z_2^{i'})\,d\theta.
\end{split}
\end{align}
For $j\le j_2$, Bernstein's inequality,
\eqref{eq:endpoint-dyadic-holder}, and the Malliavin--Sobolev product
estimate imply
\begin{align}
\label{eq:endpoint-two-low}
 \sum_{j=-1}^{j_2}
 \big|\EE\langle G,\Delta_jf(Z_1)-\Delta_jf(Z_2)\rangle\big| &\lesssim
 \mathfrak L_2^{\sigma+1}
 \|f\|_{\mathcal C^{-\sigma}}
 \|G\|_{\mathbb D^{\ell,\bar q}}
 \|Z_1-Z_2\|_{\mathbb D^{\ell,\bar q}}.
\end{align}

For $j>j_2$, \Cref{lem:dyadic-primitives} with order $\ell$ applied to
$\partial_{i'}\Delta_jf_i$ gives
\begin{align*}
 \partial_{i'}\Delta_jf_i
 =\sum_{|\alpha|=\ell}\partial^\alpha
 \psi_{j,i,i',\alpha},
 \qquad
 \|\psi_{j,i,i',\alpha}\|_\infty
 \lesssim2^{-(\ell-\sigma-1)j}
 \|f\|_{\mathcal C^{-\sigma}}.
\end{align*}
The exponent is positive because $\ell>\sigma+1$. Using
\eqref{eq:endpoint-two-ftc}, applying
\eqref{eq:Malliavin-IBP-formula} at $Z_\theta$ with multiplier
$G_i(Z_1^{i'}-Z_2^{i'})$, and then using \Cref{lem:IBP-weight-bound} and the
Malliavin--Sobolev product estimate gives
\begin{align*}
 \big|\EE\big[
 G_i(Z_1^{i'}-Z_2^{i'})\partial^\alpha
 \psi_{j,i,i',\alpha}(Z_\theta)
 \big]\big| &\lesssim
 2^{-(\ell-\sigma-1)j}
 \|f\|_{\mathcal C^{-\sigma}}
 \|G\|_{\mathbb D^{\ell,\bar q}}
 \|Z_1-Z_2\|_{\mathbb D^{\ell,\bar q}}
 \mathfrak L_2^\ell.
\end{align*}
Integrating over $\theta$, summing over $j>j_2$, and using the relation between $j_2$ and $\lambda_2$ gives
\begin{align*}
 \sum_{j>j_2}
 \big|\EE\langle G,\Delta_jf(Z_1)-\Delta_jf(Z_2)\rangle\big|
 &\lesssim
 \lambda_2^{\ell-\sigma-1}
 \|f\|_{\mathcal C^{-\sigma}}
 \|G\|_{\mathbb D^{\ell,\bar q}}
 \|Z_1-Z_2\|_{\mathbb D^{\ell,\bar q}}\mathfrak L_2^\ell\\
 &=
 \|f\|_{\mathcal C^{-\sigma}}
 \|G\|_{\mathbb D^{\ell,\bar q}}
 \|Z_1-Z_2\|_{\mathbb D^{\ell,\bar q}}\mathfrak L_2^{\sigma+1}.
\end{align*}
Together with \eqref{eq:endpoint-two-low}, this proves
\eqref{eq:odd-two-point-IBP}.
\end{proof}

For non-negative spatial regularity, the one-point estimate is immediate from
boundedness, whereas the two-point estimate still benefits from one spatial
difference.

\begin{lemma}[Two-point estimate at non-negative order]
\label{lem:nonnegative-two-point-IBP}
Let $a\in[0,1)$ and let $\ell\in\N$ satisfy $\ell>1-a$. There exist
$q>1$, $\bar q>1$, and $C\ge0$, depending only on $a,\ell$, and $d$, such
that, for every $f\in C_b^\infty(\R^d;\R^d)$ and all $\R^d$-valued random
variables $G,Z_1,Z_2$, with $Z_\theta:=\theta Z_1+(1-\theta)Z_2$,
\begin{align}
\label{eq:nonnegative-two-point-IBP}
 \big|\EE\langle G,f(Z_1)-f(Z_2)\rangle\big| &\le C\|f\|_{\mathcal C^a}
 \|G\|_{\mathbb D^{\ell,\bar q}}
 \|Z_1-Z_2\|_{\mathbb D^{\ell,\bar q}}
 \sup_{\theta\in[0,1]}
 \Lambda_{\ell,q}(Z_\theta)^{1-a}.
\end{align}
\end{lemma}

\begin{proof}
If the supremum on the right-hand side is infinite, there is nothing to
prove. Otherwise, denote it by $\mathfrak L$, set $\lambda=\mathfrak L^{-1}$,
and let $j_\lambda$ be chosen so that $2^{-j_\lambda-1}<\lambda\le2^{-j_\lambda}$. As before, we work with finite
Littlewood--Paley sums and pass to the limit at the end.
Each dyadic block satisfies the identity \eqref{eq:endpoint-two-ftc}.
For $j\le j_\lambda$, Bernstein's inequality and
\eqref{eq:endpoint-dyadic-holder} yield
$\|\nabla\Delta_jf\|_\infty\lesssim2^{(1-a)j}\|f\|_{\mathcal C^a}$, and
H\"older's inequality, followed by summation of the geometric series, gives
\begin{align*}
 \sum_{j=-1}^{j_\lambda}
 \big|\EE\langle G,\Delta_jf(Z_1)-\Delta_jf(Z_2)\rangle\big| &\lesssim
 \mathfrak L^{1-a}\|f\|_{\mathcal C^a}
 \|G\|_{\mathbb D^{\ell,\bar q}}
 \|Z_1-Z_2\|_{\mathbb D^{\ell,\bar q}}.
\end{align*}
For $j>j_\lambda$, \Cref{lem:dyadic-primitives} with order $\ell$ applied to
$\partial_{i'}\Delta_jf_i$ gives
\begin{align*}
 \partial_{i'}\Delta_jf_i
 =\sum_{|\alpha|=\ell}\partial^\alpha \psi_{j,i,i',\alpha},
 \qquad
 \|\psi_{j,i,i',\alpha}\|_\infty
 \lesssim2^{-(\ell+a-1)j}\|f\|_{\mathcal C^a}.
\end{align*}
The exponent is positive by assumption. Applying
\eqref{eq:Malliavin-IBP-formula} at $Z_\theta$, with multiplier
$G_i(Z_1^{i'}-Z_2^{i'})$, and then using \Cref{lem:IBP-weight-bound} and the
Malliavin--Sobolev product estimate gives
\begin{align*}
 \big|\EE\big[
 G_i(Z_1^{i'}-Z_2^{i'})\partial^\alpha
 \psi_{j,i,i',\alpha}(Z_\theta)\big]\big| &\lesssim
 2^{-(\ell+a-1)j}\|f\|_{\mathcal C^a}
 \|G\|_{\mathbb D^{\ell,\bar q}}
 \|Z_1-Z_2\|_{\mathbb D^{\ell,\bar q}}\mathfrak L^\ell.
\end{align*}
Integrating in $\theta$ and summing over $j>j_\lambda$ produces the factor
$\lambda^{\ell+a-1}\mathfrak L^\ell=\mathfrak L^{1-a}$.
Adding the low- and high-frequency estimates proves
\eqref{eq:nonnegative-two-point-IBP}.
\end{proof}

We now specialise these estimates to the regularity $\gamma$ of the singular
drift. Since $\lfloor-\gamma\rfloor+3>2-\gamma$, this number of
integrations by parts suffices for all the one-point, two-point, and
four-point estimates below.

\begin{lemma}\label{lem:malliavin-1}
There exist $q>1$, $\bar q>1$, and $C\ge0$, depending only on
$\gamma$ and $d$, such that, for every
$b\in C_b^\infty(\R^d;\R^d)$ and all $\R^d$-valued random variables $G,Z_1$,
\begin{align*}
 \big|\EE\langle G,b(Z_1)\rangle\big|
 \le C\|b\|_{\mathcal C^\gamma}
 \|G\|_{\mathbb D^{\lfloor-\gamma\rfloor+3,\bar q}}
 \begin{cases}
 \Lambda_{\lfloor-\gamma\rfloor+3,q}(Z_1)^{-\gamma},&\gamma<0,\\
 1,&\gamma\ge0.
 \end{cases}
\end{align*}
\end{lemma}

\begin{proof}
If $\gamma<0$, apply \eqref{eq:odd-one-point-IBP} with
$\sigma=-\gamma$ and $\ell=\lfloor-\gamma\rfloor+3$. If $\gamma\ge0$, use
$|\EE\langle G,b(Z_1)\rangle|\le\|b\|_\infty\|G\|_{L_1}
\le\|b\|_{\mathcal C^\gamma}\|G\|_{\mathbb D^{\lfloor-\gamma\rfloor+3,\bar q}}$.
\end{proof}

\begin{lemma}\label{lem:malliavin-2}
There exist $q>1$, $\bar q>1$, and $C\ge0$, depending only on
$\gamma$ and $d$, such that, for every
$b\in C_b^\infty(\R^d;\R^d)$ and all $\R^d$-valued random variables
$G,Z_1,Z_2$, with $Z_\theta:=\theta Z_1+(1-\theta)Z_2$,
\begin{align*}
 \big|\EE\langle G,b(Z_1)-b(Z_2)\rangle\big| &\le C\|b\|_{\mathcal C^\gamma}
 \|G\|_{\mathbb D^{\lfloor-\gamma\rfloor+3,\bar q}}
 \|Z_1-Z_2\|_{\mathbb D^{\lfloor-\gamma\rfloor+3,\bar q}}
 \sup_{\theta\in[0,1]}
 \Lambda_{\lfloor-\gamma\rfloor+3,q}(Z_\theta)^{1-\gamma}.
\end{align*}
\end{lemma}

\begin{proof}
If $\gamma<0$, apply \eqref{eq:odd-two-point-IBP} with
$\sigma=-\gamma$ and $\ell=\lfloor-\gamma\rfloor+3$. If $\gamma\ge0$, apply
\eqref{eq:nonnegative-two-point-IBP} with
$a=\gamma$ and $\ell=\lfloor-\gamma\rfloor+3$.
\end{proof}

\begin{lemma}\label{lem:malliavin-3}
There exist $q>1$, $\bar q>1$, and $C\ge0$, depending only on
$\gamma$ and $d$, such that, for every
$b\in C_b^\infty(\R^d;\R^d)$ and all
$\R^d$-valued random variables $G,Z_0,Z_1,Z_2$, with
\begin{align*}
 Z_{\theta,\vartheta}
 :=(\theta+\vartheta-1)Z_0
 +(1-\theta)Z_1+(1-\vartheta)Z_2,
 \qquad\theta,\vartheta\in[0,1],
\end{align*}
one has
\begin{align*}
 &\left|\EE\left\langle G,
 b(Z_0)-b(Z_1)-b(Z_2)+b(Z_1+Z_2-Z_0)
 \right\rangle\right|\\
 &\quad\le C\|b\|_{\mathcal C^\gamma}
 \|G\|_{\mathbb D^{\lfloor-\gamma\rfloor+3,\bar q}}
 \|Z_0-Z_1\|_{\mathbb D^{\lfloor-\gamma\rfloor+3,\bar q}}
 \|Z_0-Z_2\|_{\mathbb D^{\lfloor-\gamma\rfloor+3,\bar q}}
 \sup_{\theta,\vartheta\in[0,1]}
 \Lambda_{\lfloor-\gamma\rfloor+3,q}(Z_{\theta,\vartheta})^{2-\gamma}.
\end{align*}
\end{lemma}

\begin{proof}
If $\mathfrak L:=\sup_{\theta,\vartheta\in[0,1]}
\Lambda_{\lfloor-\gamma\rfloor+3,q}(Z_{\theta,\vartheta})$
is infinite, the estimate is immediate. Otherwise, set
$\lambda:=\mathfrak L^{-1}$ and define $j_\lambda$ as in \Cref{lem:nonnegative-two-point-IBP}. As before, we work with finite
Littlewood--Paley sums and pass to the limit at the end.
Since $Z_{\theta,\vartheta}=Z_0-(1-\theta)(Z_0-Z_1)-(1-\vartheta)(Z_0-Z_2)$,
for every smooth $h:\R^d\to\R^d$ the two-parameter fundamental theorem of
calculus gives
\begin{align}
\label{eq:two-parameter-FTC}
 h(Z_0)-h(Z_1)-h(Z_2)+h(Z_1+Z_2-Z_0) &=
 \int_0^1\int_0^1
 \nabla^2h(Z_{\theta,\vartheta})
 [Z_0-Z_1,Z_0-Z_2] \,d\theta\,d\vartheta.
\end{align}

Write $b=\sum_{j\ge-1}\Delta_jb$ and apply
\eqref{eq:two-parameter-FTC} to each block. For $j\le j_\lambda$,
Bernstein's inequality and \eqref{eq:endpoint-dyadic-holder} give
$\|\nabla^2\Delta_jb\|_\infty\lesssim2^{(2-\gamma)j}\|b\|_{\mathcal C^\gamma}$.
Using the Malliavin--Sobolev product estimate for
$G,Z_0-Z_1,Z_0-Z_2$ and summing the geometric series,
\begin{align}\label{eq:four-endpoint-low}
\begin{split}
 &\sum_{j=-1}^{j_\lambda}
 \left|\EE\left\langle G,
 \Delta_jb(Z_0)-\Delta_jb(Z_1)-\Delta_jb(Z_2)
 +\Delta_jb(Z_1+Z_2-Z_0)
 \right\rangle\right|\\
 &\quad\lesssim
 \mathfrak L^{2-\gamma}\|b\|_{\mathcal C^\gamma}
 \|G\|_{\mathbb D^{\lfloor-\gamma\rfloor+3,\bar q}}
 \|Z_0-Z_1\|_{\mathbb D^{\lfloor-\gamma\rfloor+3,\bar q}}
 \|Z_0-Z_2\|_{\mathbb D^{\lfloor-\gamma\rfloor+3,\bar q}}.
\end{split}
\end{align}

Let $j>j_\lambda$. For each component $i$ and indices $i_1,i_2$,
\Cref{lem:dyadic-primitives} with order $\lfloor-\gamma\rfloor+3$ applied to
$\partial_{i_1i_2}\Delta_jb_i$ gives
\begin{align*}
 \partial_{i_1i_2}\Delta_jb_i
 =\sum_{|\alpha|=\lfloor-\gamma\rfloor+3}\partial^\alpha
 \psi_{j,i,i_1,i_2,\alpha},
 \qquad
 \|\psi_{j,i,i_1,i_2,\alpha}\|_\infty
 \lesssim
 2^{-(\lfloor-\gamma\rfloor+\gamma+1)j}
 \|b\|_{\mathcal C^\gamma},
\end{align*}
and the exponent $\lfloor-\gamma\rfloor+\gamma+1$ is positive. Expanding
\eqref{eq:two-parameter-FTC} componentwise and applying
\eqref{eq:Malliavin-IBP-formula} at $Z_{\theta,\vartheta}$, with multiplier
$G_i(Z_0-Z_1)^{i_1}(Z_0-Z_2)^{i_2}$, gives, by
\Cref{lem:IBP-weight-bound} and the product estimate,
\begin{align*}
 &\left|
 \EE\left[
 G_i(Z_0-Z_1)^{i_1}(Z_0-Z_2)^{i_2}
 \partial^\alpha
 \psi_{j,i,i_1,i_2,\alpha}(Z_{\theta,\vartheta})
 \right]\right|\\
 &\quad\lesssim
 2^{-(\lfloor-\gamma\rfloor+\gamma+1)j}
 \|b\|_{\mathcal C^\gamma}
 \|G\|_{\mathbb D^{\lfloor-\gamma\rfloor+3,\bar q}}
 \|Z_0-Z_1\|_{\mathbb D^{\lfloor-\gamma\rfloor+3,\bar q}}
 \|Z_0-Z_2\|_{\mathbb D^{\lfloor-\gamma\rfloor+3,\bar q}}
 \mathfrak L^{\lfloor-\gamma\rfloor+3}.
\end{align*}
Integrating over $(\theta,\vartheta)$ and summing over $j>j_\lambda$ yields
\begin{align}\label{eq:four-endpoint-high}
\begin{split}
 &\sum_{j>j_\lambda}
 \left|\EE\left\langle G,
 \Delta_jb(Z_0)-\Delta_jb(Z_1)-\Delta_jb(Z_2)
 +\Delta_jb(Z_1+Z_2-Z_0)
 \right\rangle\right|\\
 &\quad\lesssim
 \lambda^{\lfloor-\gamma\rfloor+\gamma+1}
 \mathfrak L^{\lfloor-\gamma\rfloor+3}
 \|b\|_{\mathcal C^\gamma}
 \|G\|_{\mathbb D^{\lfloor-\gamma\rfloor+3,\bar q}}
 \|Z_0-Z_1\|_{\mathbb D^{\lfloor-\gamma\rfloor+3,\bar q}}
 \|Z_0-Z_2\|_{\mathbb D^{\lfloor-\gamma\rfloor+3,\bar q}},
\end{split}
\end{align}
and $\lambda^{\lfloor-\gamma\rfloor+\gamma+1}\mathfrak L^{\lfloor-\gamma\rfloor+3}
=\mathfrak L^{2-\gamma}$. Combining \eqref{eq:four-endpoint-low} and
\eqref{eq:four-endpoint-high} proves the estimate.
\end{proof}

\subsection{Conditioning and restarting fractional SDEs}

The integration-by-parts estimates above concern random variables at a fixed
time. To apply them inside conditional sewing arguments, we also need a
precise description of the equation after a conditioning time. This
subsection provides the required fixed-time restarting identities.

In this subsection, $b\in C_b^\infty(\R^d;\R^d)$ and $F$ satisfies
\Cref{ass:dissipative-drift}. For $s<T$, $x\in\R^d$, and
$\zeta\in C([s,T];\R^d)$, let $\mathcal S_{F+b,s}(x,\zeta)$ denote the
solution of
\begin{align*}
 \mathcal S_{F+b,s}(x,\zeta)_t
 =x+\int_s^t(F+b)\big(\mathcal S_{F+b,s}(x,\zeta)_r\big)\,dr
 +\zeta_{s,t},
 \qquad t\in[s,T].
\end{align*}
The coercivity of $F$, the boundedness of $b$, and local ODE stability imply
that $\mathcal S_{F+b,s}:\R^d\times C([s,T];\R^d)\to C([s,T];\R^d)$
is continuous when the path spaces are equipped with the uniform topology.
Indeed, solutions corresponding to data in a bounded set remain in a common
compact set by the coercivity estimate; on that compact set $F+b$ is
Lipschitz, and continuity follows from Gr\"onwall's inequality.

For a deterministic input $\mathbf v$, the solution of the smooth equation
started from $x$ at time $s$ is
\begin{align}\label{eq:X-Ito-map}
\begin{split}
 X^{s,x,\mathbf v;b}
 :=\mathcal S_{F+b,s}(x,\mathbf v+B^s),
\end{split}
\end{align}
and, as for the flow, we use the same notation when $x$ and $\mathbf v$ are
replaced by $\mathcal F_s$-measurable random variables. We omit $b$ from the
notation when the drift is clear from the context. For $b=0$, this is the
regularising flow: $X^{s,x,\mathbf v;0}=\phi^{s,x,\mathbf v}$. Since $B^s$
is independent of $\mathcal F_s$, we have the following fixed-time
conditioning identity.

\begin{lemma}[Conditioning at the initial time]
\label{lem:conditioning-initial}
Let $\xi$ and $\mathbf w$ be $\mathcal F_s$-measurable random variables with
values in $\R^d$ and $C([s,T];\R^d)$, respectively. For every fixed
$t\in[s,T]$ and every bounded Borel function $f:\R^d\to\R$,
\begin{align}\label{eq:conditioning-initial}
\begin{split}
 \EE^s\!\left[f\bigl(X_t^{s,\xi,\mathbf w;b}\bigr)\right]
 =\left.
 \EE\!\left[f\bigl(X_t^{s,x,\mathbf v;b}\bigr)\right]
 \right|_{x=\xi,\,\mathbf v=\mathbf w}.
\end{split}
\end{align}
The same identity holds for non-negative Borel $f$, with both sides allowed
to take the value $+\infty$.
\end{lemma}

\begin{proof}
By \eqref{eq:X-Ito-map},
$X_t^{s,\xi,\mathbf w;b}=\mathcal S_{F+b,s}(\xi,\mathbf w+B^s)_t$.
The pair $(\xi,\mathbf w)$ is $\mathcal F_s$-measurable, whereas $B^s$ is
independent of $\mathcal F_s$. The fixed-time conditioning identity
\cite[Lemma~A.3]{DHL25}, applied to the measurable map
$(x,\mathbf v,\mathbf y)\mapsto
f(\mathcal S_{F+b,s}(x,\mathbf v+\mathbf y)_t)$, gives
\eqref{eq:conditioning-initial}. The non-negative case follows by applying
the bounded identity to $f\wedge n$, $n\in\N$, and using conditional
monotone convergence.
\end{proof}

We now restart a solution originally started at $t_0$. Let
$X=X^{t_0,X_{t_0},V;b}$ and fix $s\in[t_0,T]$. Subtracting the equation at time
$s$ from the equation at a fixed time $t\ge s$ and using
\eqref{eq:restart-RL} gives
\begin{align*}
 X_t
 =X_s+\int_s^t(F+b)(X_r)\,dr
 +V_{s,t}+\bar B_t^{t_0,s}+B_t^s.
\end{align*}
Since $\bar B_s^{t_0,s}=0$, the right-hand side is the smooth equation
started at time $s$ from $X_s$ with input $V+\bar B^{t_0,s}$, and both $X_s$
and $V+\bar B^{t_0,s}$ are $\mathcal F_s$-measurable.
Pathwise uniqueness for the smooth equation therefore gives, for every fixed
$t\in[s,T]$, $X_t=X_t^{s,X_s,V+\bar B^{t_0,s};b}$ almost surely.
Combining this identity with \Cref{lem:conditioning-initial}, applied
jointly to finitely many equations, gives the form used throughout the
paper.

\begin{corollary}[Conditional restarting identity]
\label{cor:conditioning-restart}
Let $M\ge1$ and, for $i=1,\ldots,M$, let
$X^i=X^{t_0,X_{t_0}^i,V^i;b_i}$ be solutions of smooth equations driven by the
same Riemann--Liouville process, with drifts $F+b_i$,
$b_i\in C_b^\infty(\R^d;\R^d)$, $\mathcal F_{t_0}$-measurable initial
conditions $X_{t_0}^i$, and inputs $V^i$. Let $s\in[t_0,T]$. For every fixed
$t\in[s,T]$ and every Borel function $f:(\R^d)^M\to\R$ which is bounded,
non-negative, or such that $f(X_t^1,\ldots,X_t^M)$ is integrable,
\begin{align*}
 \EE^s\!\left[f(X_t^1,\ldots,X_t^M)\right] &=
 \left.
 \EE\!\left[
 f\bigl(X_t^{s,x_1,\mathbf v_1;b_1},\ldots,
 X_t^{s,x_M,\mathbf v_M;b_M}\bigr)
 \right]
 \right|_{x_i=X_s^i,\,\mathbf v_i=V^i+\bar B^{t_0,s}}.
\end{align*}
(In the integrable case, apply the non-negative case to $f^+$ and $f^-$;
both are then finite almost surely.) The same identity holds when some of
the $X^i$ are replaced by regularising flows $\phi^{u,\xi,\mathbf w}$ restarted at a time $u\le s$, with $\mathbf w$
and $\xi$ being $\mathcal F_u$-measurable. Indeed, these are the case $b_i=0$
of the smooth equation started at time $u$ instead of $t_0$, and the
corresponding substitution is $x_i=\phi_s^{u,\xi,\mathbf w}$ and
$\mathbf v_i=\mathbf w+\bar B^{u,s}$.
\end{corollary}

\subsection{Gaussian estimates for the history and the local scale}

When an interval is longer than the local scale \eqref{eq:local-scale} of an
input, we compensate by a high power of its length, paid for with moments of
the H\"older seminorm of the input. Since $\tau_{\mathbf v;I}^{\chi_\eps}
=(1\vee[\mathbf v]_{C^{H-\eps}(I)})^{-1}$, one has, for every $h\in(0,1]$,
$\nu\ge0$, and $\rho\ge1$,
\begin{align}
\label{eq:bad-event-basic}
 \ind_{\{h>\tau_{\mathbf v;I}\}}
 \bigl(1+[\mathbf v]_{C^{H-\eps}(I)}\bigr)^\nu &\le
 h^{\chi_\eps\rho}
 \bigl(1+[\mathbf v]_{C^{H-\eps}(I)}\bigr)^{\nu+\rho}.
\end{align}
The inputs to which \eqref{eq:bad-event-basic} is applied are the prescribed
input $V$ and the histories of the noise, whose H\"older seminorms have
Gaussian moments.

\begin{lemma}[H\"older moments of the history]\label{lem:history-holder-moments}
Let $H\in(0,1)$ and $\eps\in(0,H)$. There exists
$C=C(H,\eps,d)$ such that, for every $q\in[2,\infty)$,
\begin{align}\label{eq:full-history-holder-moments}
\begin{split}
 \sup_{a\le u\le u'}
 \left\|[\bar B^{a,u}]_{C^{H-\eps}([u',u'+1])}\right\|_{L_q}
 +\sup_{u\le u'}
 \left\|[\bar B^u]_{C^{H-\eps}([u',u'+1])}\right\|_{L_q}
 \le C\sqrt q.
\end{split}
\end{align}
\end{lemma}

\begin{proof}
When $H=1/2$, all the history processes in the statement vanish identically,
and the result is immediate.  We therefore assume that $H\ne1/2$ and first
consider $\bar B^{a,u}$ for fixed $a\le u\le u'$. Let
$u'\le t_1<t_2\le u'+1$ and put $h:=t_2-t_1$. Since the integrand defining $\bar B^{a,u}$ is
deterministic, the It\^o isometry, applied componentwise and then summed over
the $d$ coordinates, gives
\begin{align*}
 \EE\left|\bar B_{t_2}^{a,u}-\bar B_{t_1}^{a,u}\right|^2 &=c_H^2d\int_{a}^u
 \left|(t_2-r)^{H-\frac12}-(t_1-r)^{H-\frac12}\right|^2dr.
\end{align*}
Set $z_0:=t_1-u\ge0$ (which may be large) and substitute $r=u-y$. Enlarging the integration
interval from $[0,u-a]$ to $[0,\infty)$ and then scaling $y=hy'$ yields
\begin{align*}
 \EE\left|\bar B_{t_2}^{a,u}-\bar B_{t_1}^{a,u}\right|^2
 \le c_H^2d\,h^{2H}\mathcal I_H(z_0/h),
 \qquad
 \mathcal I_H(z):=\int_0^\infty
 \left|(z+1+y)^{H-\frac12}-(z+y)^{H-\frac12}\right|^2dy.
\end{align*}
The quantity $\mathcal I_H(z)$ is bounded uniformly in $z\ge0$.  Indeed, for
$y\in(0,1)$, the integrand is bounded by $(z+y)^{2H-1}\le y^{2H-1}$ if
$H<1/2$, and by $1$ if $H>1/2$, by concavity and subadditivity of
$y\mapsto y^{H-1/2}$; in both cases it is integrable on $(0,1)$. For
$y\ge1$, the mean value theorem bounds the integrand by
$C(z+y)^{2H-3}\le Cy^{2H-3}$, and $2H-3<-1$ because $H<1$. Therefore
\begin{align}\label{eq:history-L2-increments}
\begin{split}
 \sup_{a\le u\le t_1}
 \|\bar B_{t_2}^{a,u}-\bar B_{t_1}^{a,u}\|_{L_2}
 \le C_{H,d}|t_2-t_1|^H,
\end{split}
\end{align}
the constant being independent of $a$, $u$, and $t_1$: the parameters $a$
and $u$ disappear after the change of variables, and the distance $t_1-u$
only enters through $z_0$, in which $\mathcal I_H$ is uniformly bounded.

The increment is a centred Gaussian vector. Gaussian moment equivalence and
\eqref{eq:history-L2-increments} imply, for every $q\ge2$,
\begin{align}\label{eq:history-Lq-increments}
\begin{split}
 \sup_{a\le u\le t_1}
 \|\bar B_{t_2}^{a,u}-\bar B_{t_1}^{a,u}\|_{L_q}
 \le C_{H,d}\sqrt q\,|t_2-t_1|^H.
\end{split}
\end{align}
We now use the standard dyadic proof of Kolmogorov's continuity theorem on
the unit interval. For $j\ge0$ and $0\le i<2^j$, let
$t_{j,i}:=u'+i2^{-j}$ and
$\mathfrak M_j:=\max_{0\le i<2^j}
|\bar B_{t_{j,i+1}}^{a,u}-\bar B_{t_{j,i}}^{a,u}|$.
The deterministic chaining estimate gives
\begin{align}\label{eq:history-dyadic-chaining}
\begin{split}
 [\bar B^{a,u}]_{C^{H-\eps}([u',u'+1])}
 \le C_{H,\eps}\sum_{j\ge0}2^{j(H-\eps)}\mathfrak M_j.
\end{split}
\end{align}
Moreover, the elementary inequality
$\|\max_i \zeta_i\|_{L_q}\le(\sum_i\|\zeta_i\|_{L_q}^q)^{1/q}$ and
\eqref{eq:history-Lq-increments} give
$\|\mathfrak M_j\|_{L_q}\le C_{H,d}\sqrt q\,2^{j/q}2^{-jH}$.
Choose $q_0\ge2/\eps$. For $q\ge q_0$, Minkowski's inequality and
\eqref{eq:history-dyadic-chaining} yield
\begin{align*}
 \left\|[\bar B^{a,u}]_{C^{H-\eps}([u',u'+1])}\right\|_{L_q}
 &\le C\sqrt q\sum_{j\ge0}2^{-j(\eps-1/q)} \le C_{H,\eps,d}\sqrt q,
\end{align*}
since $\eps-1/q\ge\eps/2$. For $2\le q<q_0$, monotonicity of $L_q$-norms
reduces the estimate to the case $q=q_0$, after changing the constant by a
factor depending only on $q_0$, hence only on $\eps$. This proves the bound
for $\bar B^{a,u}$. For the full history $\bar B^u$ on $[u',u'+1]$, $u\le u'$, the same change
of variables starts with the integration range $[0,\infty)$ rather than
enlarging a finite interval to it, and the rest of the argument is unchanged.
\end{proof}

We combine \eqref{eq:bad-event-basic} with the Gaussian moments in the
following form. Recall that the seminorm of a path restarted at time $u$ on
an interval $I\ni u$ is computed after extending it constantly to the left
of $u$; since $I\cap[u,\infty)\subset[u,u+1]$ when $|I|\le1$,
\Cref{lem:history-holder-moments} applies on every such interval.

\begin{lemma}[Bad-event moments]\label{lem:bad-event-moments}
Let $\eps\in(0,H)$, $\mathfrak a\ge0$, $\nu\ge0$, and $q\ge1$. There exists a
constant $C$, depending only on $\mathfrak a,\nu,q,H,\eps$, and $k$, such that
the following holds. Let $I$ be an interval of length at most one, and let
$\mathbf w$ be a random path with
$\|[\mathbf w]_{C^{H-\eps}(I)}\|_{L_{q'}}\le\mathfrak a\sqrt{q'}$ for every
$q'\ge2$. Then, for every $\rho\ge1$ and every $h\in(0,1]$,
\begin{align}\label{eq:bad-event-moments}
 \big\|\ind_{\{h>\tau_{\mathbf w;I}\}}
 (1+[\mathbf w]_{C^{H-\eps}(I)})^\nu\big\|_{L_q}
 \le C^\rho\rho^{\rho/2}h^{\chi_\eps\rho}.
\end{align}
For any fixed $\rho'\ge0$, the bound also holds with
$h^{\chi_\eps(\rho+\rho')}$ in place of $h^{\chi_\eps\rho}$, after changing
$C$. By \Cref{lem:history-holder-moments}, the assumption on
$\mathbf w$ holds with $\mathfrak a=C(H,\eps,d)$ for every history
$\bar B^{a,u}$ and $\bar B^u$, and, under \Cref{ass:input}, with
$\mathfrak a=\Lambda_{V,\eps}+C(H,\eps,d)$ for $\mathbf w=V+\bar B^{t_0,u}$ on
$I_u$.
\end{lemma}

\begin{proof}
By \eqref{eq:bad-event-basic} and the moment assumption,
\begin{align*}
 \big\|\ind_{\{h>\tau_{\mathbf w;I}\}}(1+[\mathbf w]_{C^{H-\eps}(I)})^\nu
 \big\|_{L_q}
 &\le h^{\chi_\eps\rho}
 \|1+[\mathbf w]_{C^{H-\eps}(I)}\|_{L_{q(\nu+\rho)}}^{\nu+\rho}
 \le h^{\chi_\eps\rho}
 \bigl(C\sqrt{q(\nu+\rho)\vee2}\bigr)^{\nu+\rho},
\end{align*}
where moments of order below two are bounded by the second moment. Since $\nu$
and $q$ are fixed, $(1+\nu/\rho)^{\rho/2}\le e^{\nu/2}$ and the remaining
fixed powers of $\rho$ are absorbed into $C^\rho$. Fixed shifts of $\rho$ are
absorbed in the same way. No independence between $\mathbf w$ and other
factors is needed when this estimate is combined with H\"older's inequality.
\end{proof}

\subsection{Stochastic sewing}

We recall the stochastic sewing lemma of \cite{le2020stochastic} in the form
used in the a priori and stability arguments. For a two-parameter process
$A$ and $s<u<t$, we write $\delta A_{s,u,t}:=A_{s,t}-A_{s,u}-A_{u,t}$.

\begin{lemma}[Stochastic sewing lemma {\cite[Theorem~2.1]{le2020stochastic}}]
\label{lem:stochastic-sewing}
Let $S<S'$ and $m\in[2,\infty)$. Let $A:\Delta_{[S,S']}\to L_m$ be such that
$A_{s,t}$ is $\mathcal F_t$-measurable. Assume that there exist
$\mathsf C_1,\mathsf C_2\ge0$ and $\varsigma_1,\varsigma_2>0$ such that, for
every $S\le s<t\le S'$ and $u:=(s+t)/2$,
\begin{align}
 \|\EE^s\delta A_{s,u,t}\|_{L_m}
 &\le\mathsf C_1(t-s)^{1+\varsigma_1},\label{eq:stochastic-sewing-conditional-defect}\\
 \|\delta A_{s,u,t}\|_{L_m}
 &\le\mathsf C_2(t-s)^{1/2+\varsigma_2}.\label{eq:stochastic-sewing-defect}
\end{align}
Then there exists a unique, up to modification, adapted process
$\mathcal Y=(\mathcal Y_t)_{t\in[S,S']}$ with $\mathcal Y_S=0$ for which there
exist $\mathsf C'<\infty$ and $\varsigma'>0$ such that, for all
$S\le s<t\le S'$,
\begin{align*}
 \|\mathcal Y_t-\mathcal Y_s-A_{s,t}\|_{L_m}
 \le\mathsf C'(t-s)^{1/2+\varsigma'}
 \quad\text{and}\quad
 \|\EE^s(\mathcal Y_t-\mathcal Y_s-A_{s,t})\|_{L_m}
 \le\mathsf C'(t-s)^{1+\varsigma'}.
\end{align*}
Moreover, there exists a constant $C$, depending only on $m$, $\varsigma_1$,
and $\varsigma_2$, such that, for all $S\le s<t\le S'$,
\begin{align}
 \|\mathcal Y_t-\mathcal Y_s-A_{s,t}\|_{L_m}
 &\le C\mathsf C_1(t-s)^{1+\varsigma_1}
 +C\mathsf C_2(t-s)^{1/2+\varsigma_2},\label{eq:stochastic-sewing-error}\\
 \|\EE^s(\mathcal Y_t-\mathcal Y_s-A_{s,t})\|_{L_m}
 &\le C\mathsf C_1(t-s)^{1+\varsigma_1}.\label{eq:stochastic-sewing-conditional-error}
\end{align}
In particular,
$\|\mathcal Y_t-\mathcal Y_s\|_{L_m}\le\|A_{s,t}\|_{L_m}
+C\mathsf C_1(t-s)^{1+\varsigma_1}+C\mathsf C_2(t-s)^{1/2+\varsigma_2}$.
\end{lemma}

The lemma is due to \cite[Theorem~2.1]{le2020stochastic}; a proof of the
statement in the above form, in which the sewn process is characterised by
the two bounds on $\mathcal Y_t-\mathcal Y_s-A_{s,t}$, can be found in
\cite{dareiotis2024regularisation}.

In all our applications, the drift is smooth and the process $\mathcal Y$ is
known in advance: it is an explicit Lebesgue integral along the solutions.
The smoothness of the drift shows directly that $\mathcal Y_t-\mathcal Y_s$ is
close to $A_{s,t}$ in the sense of the characterisation above, with
constants $\mathsf C'$ depending on the $C^1$ norm of the drift. (When the
solutions along which $\mathcal Y$ is integrated have a singular drift, as
in \Cref{sec:stab}, this also uses their comparison estimate; see the proof
of \Cref{lem:stability-sewing}.) Hence
$\mathcal Y$ is the process given by \Cref{lem:stochastic-sewing}, and the
estimates \eqref{eq:stochastic-sewing-error}--\eqref{eq:stochastic-sewing-conditional-error}
hold with constants depending only on $\mathsf C_1$ and $\mathsf C_2$, which
do not involve the smoothness of the drift.

\section{The regularising flow}\label{sec:flow}

We now establish the estimates for the flow introduced in
\eqref{eq:flow-front}. These estimates provide the analytic input for the
Malliavin integration-by-parts arguments used in the a priori estimates. We
first state the complete result and then prove it in four stages: pathwise
estimates, Malliavin derivative estimates, spatial and four-point estimates
for those derivatives, and finally non-degeneracy of the Malliavin
covariance.

Throughout this section, $F$ satisfies \Cref{ass:dissipative-drift} and
$\eps\in(0,H)$ is fixed. For an interval $I$ of length at most one,
$(s,t)\in\Delta_I$, and paths
$\mathbf v^1,\ldots,\mathbf v^L\in C^{H-\eps}(I;\R^d)$, we define the input
control
\begin{align}\label{eq:input-control}
 R_{s,t}^{\mathbf v^1,\ldots,\mathbf v^L;I}
 :=1+\sum_{i=1}^L
 \ind_{\{t-s>\tau_{\mathbf v^i;I}\}}
 [\mathbf v^i]_{C^{H-\eps}(I)}.
\end{align}
Only inputs whose local scale is shorter than $t-s$ contribute. The inputs
may be random when this expression is used pathwise; in particular,
$R_{s,t}^{B^s,\mathbf v;I}$ includes the indicator associated with the local
scale of the fresh noise $B^s$. The interval $I$ is kept fixed throughout
each estimate.

In this section, $C_1,C_2$ denote deterministic constants depending only on
$d,H,\eps,k,\boldsymbol\kappa$, the number $L$ of inputs, and the
derivative orders involved; they are uniform in the initial points. Pure
Malliavin derivatives are estimated up to order $\lfloor-\gamma\rfloor+5$,
and a mixed estimate with $n$ spatial and $\ell$ Malliavin derivatives
requires $n+\ell\le\lfloor-\gamma\rfloor+5$. This is why the next
proposition is stated for $0\le j\le\lfloor-\gamma\rfloor+3$: the norm of
$D\Phi$ uses Malliavin derivatives up to order $j+1$, and the four-point
estimate uses two additional spatial derivatives.

\begin{proposition}[Estimates for the regularising flow]
\label{prop:flow-estimates}
Let $p\in[1,\infty)$,
$j\in\{0,\ldots,\lfloor-\gamma\rfloor+3\}$, $L,M\in\N$, and
$\lambda\in\mathcal A_M$. There exist constants $C_1,C_2<\infty$, depending
only on $d$, $H$, $\eps$, $k$, $j$, $p$, $L$, $M$, $\|\lambda\|_{\ell^1}$, and $\boldsymbol\kappa$,
such that the following assertions hold. Let $I=[s,\bar t]$ be an interval of
length at most one, let $t\in(s,\bar t]$, let
$\mathbf v^1,\ldots,\mathbf v^L\in C^{H-\eps}(I;\R^d)$ be deterministic
inputs, put $\mathbf v:=\sum_{i=1}^L\mathbf v^i$, and let
$\mathbf x=(x_1,\ldots,x_M)\in(\R^d)^M$.

First, using the natural Hilbert-valued extension of the Malliavin--Sobolev
norm,
\begin{align}\label{eq:upperbound-D}
\begin{split}
 \big\|D\Phi_t^{s,\mathbf x,\lambda,\mathbf v}
 \big\|_{\mathbb D^{j,p}(\mathcal H_{s,t})}
 \le C_1\big(R_{s,t}^{\mathbf v^1,\ldots,\mathbf v^L;I}\big)^{C_2}
 (t-s)^H.
\end{split}
\end{align}
Second, the Malliavin covariance of the affine combination is invertible and
satisfies
\begin{align}\label{eq:inv-matrix}
\begin{split}
 \big\|\mathcal M_{\Phi_t^{s,\mathbf x,\lambda,\mathbf v}}^{-1}
 \big\|_{L_p}
 \le C_1\big(R_{s,t}^{\mathbf v^1,\ldots,\mathbf v^L;I}\big)^{C_2}
 (t-s)^{-2H}.
\end{split}
\end{align}
Moreover, for every $x,y\in\R^d$,
\begin{align}\label{eq:upperbound-D-x}
\begin{split}
 \big\|\phi_t^{s,x,\mathbf v}-\phi_t^{s,y,\mathbf v}
 \big\|_{\mathbb D^{j,p}}
 \le C_1\big(R_{s,t}^{\mathbf v^1,\ldots,\mathbf v^L;I}\big)^{C_2}|x-y|,
\end{split}
\end{align}
and, for every $x_1,x_2,x_3,x_4\in\R^d$,
\begin{align}\label{eq:upperbound-4}
\begin{split}
 &\big\|\phi_t^{s,x_1,\mathbf v}-\phi_t^{s,x_2,\mathbf v}
 -\phi_t^{s,x_3,\mathbf v}+\phi_t^{s,x_4,\mathbf v}
 \big\|_{\mathbb D^{j,p}}\\
 &\quad\le
 C_1\big(R_{s,t}^{\mathbf v^1,\ldots,\mathbf v^L;I}\big)^{C_2}
 \Big(
 |x_1-x_2-x_3+x_4|
 +|x_1-x_2|\,|x_1-x_3|
 \Big).
\end{split}
\end{align}
\end{proposition}

The input dependence in \Cref{prop:flow-estimates} is determined by the
coming-down-from-infinity estimate, \Cref{lem:coming-down-infinity}, which we
prove for $L$ prescribed inputs in addition to the fresh noise $B^s$, keeping
the local scales of the inputs separate. All subsequent estimates use this bound only through
polynomial control of the flow and its derivatives. For readability, we
therefore write them for a single prescribed input $\mathbf v$; the proofs
are identical for several inputs, with
$R_{s,t}^{B^s,\mathbf v^1,\ldots,\mathbf v^L;I}$ in place of
$R_{s,t}^{B^s,\mathbf v;I}$.

These results are also applied to flows started at a time $u\in[s,\bar t)$
inside a fixed interval $I=[s,\bar t]$. They then hold, with the controls
computed on $I$: one applies them on $[u,\bar t]$, which has length at most
one, and uses that, with the left-extension convention of
\Cref{sec:framework}, $[\mathbf w]_{C^{H-\eps}([u,\bar t])}
\le[\mathbf w]_{C^{H-\eps}(I)}$ for every input $\mathbf w$, so that
$\tau_{\mathbf w;[u,\bar t]}\ge\tau_{\mathbf w;I}$ and the input controls
\eqref{eq:input-control} on $[u,\bar t]$ are bounded by those on $I$.

\begin{notation}\label{not:constantC}
In proofs only, $\mathbf C$ denotes a non-negative random quantity,
allowed to change from line to line, satisfying
\begin{align}\label{eq:bold-C-control}
 \mathbf C\le C_1\bigl(R_{s,t}^{B^s,\mathbf v;I}\bigr)^{C_2}
\end{align}
(with $R_{s,t}^{B^s,\mathbf v^1,\ldots,\mathbf v^L;I}$ in the case of several
inputs). The noise indicator is retained in pathwise bounds and removed when
taking moments: since the $C^{H-\eps}$-seminorm of $B^s$ has Gaussian moments
uniformly on intervals of length at most one, for every $p,q<\infty$,
\begin{align}\label{eq:flow-control-moments}
 \big\|\bigl(R_{s,t}^{B^s,\mathbf v^1,\ldots,\mathbf v^L;I}\bigr)^q
 \big\|_{L_p}
 \le C\bigl(R_{s,t}^{\mathbf v^1,\ldots,\mathbf v^L;I}\bigr)^q.
\end{align}
Statements below display the control factor explicitly; $\mathbf C$ is only
a proof abbreviation.
\end{notation}

\subsection{Pathwise estimates for the regularising flow}

We begin with the coming-down-from-infinity estimate. Besides controlling the
flow itself, its proof explains the choice of the local time scale
\eqref{eq:local-scale}. Recall that $\tau_{\mathbf v;I}
=(1\vee[\mathbf v]_{C^{H-\eps}(I)})^{-1/\chi_\eps}$ and note that
$1/\chi_\eps=(k-1)/(1+(k-1)(H-\eps))$.

\begin{lemma}[Coming down from infinity]\label{lem:coming-down-infinity}
Let $\omega\in\Omega$ and $L\in\N$. There is a deterministic constant $C$, depending only on
$H,\eps,k,\boldsymbol\kappa$, and $L$, such that, for every interval
$I=[s,\bar t]$ of length at most one, deterministic
inputs $\mathbf v^1,\ldots,\mathbf v^L\in C^{H-\eps}(I;\R^d)$ with sum
$\mathbf v$, for every $x\in\R^d$, every $t\in(s,\bar t]$, and
every $h\in(0,t-s]$,
\begin{align}\label{eq:coming-down-h}
 \bigl|\phi_t^{s,x,\mathbf v}(\omega) \bigr|
 \le C\Bigg(
 1+h^{-\frac1{k-1}}
 +h^{H-\eps}\Big(
 [B^s(\omega)]_{C^{H-\eps}(I)}
 +\sum_{i=1}^L[\mathbf v^i]_{C^{H-\eps}(I)}
 \Big)\Bigg).
\end{align}
As a consequence,
\begin{align}\label{eq:coming-down}
 \bigl|\phi_t^{s,x,\mathbf v}(\omega)\bigr|
 \le C\Bigg(
 1+(t-s)^{-\frac1{k-1}}
 +\ind_{\{t-s>\tau_{B^s;I}\}}
 [B^s(\omega)]_{C^{H-\eps}(I)}^{\frac1{(k-1)\chi_\eps}}
 +\sum_{i=1}^L
 \ind_{\{t-s>\tau_{\mathbf v^i;I}\}}
 [\mathbf v^i]_{C^{H-\eps}(I)}^{\frac1{(k-1)\chi_\eps}}
 \Bigg),
\end{align}
and in particular
\begin{align}\label{eq:coming-down-R}
 \bigl|\phi_t^{s,x,\mathbf v}(\omega)\bigr|
 \le C\Big(
 (t-s)^{-\frac1{k-1}}
 +(R_{s,t}^{B^s(\omega),\mathbf v^1,\ldots,\mathbf v^L;I})
 ^{\frac1{(k-1)\chi_\eps}}
 \Big).
\end{align}
\end{lemma}

\begin{proof} We may assume $[B^s(\omega)]_{C^{H-\eps}(I)}<\infty$ since otherwise the estimate is trivial. We omit the $\omega$-dependence in the notation.
We write
$S_I:=[B^s]_{C^{H-\eps}(I)}+\sum_{i=1}^L[\mathbf v^i]_{C^{H-\eps}(I)}$.
Fix $t\in(s,\bar t]$ and $h\in(0,t-s]$. We restart the
equation at $r:=t-h$ and, for $u\in[r,t]$, set
$\zeta_u:=\mathbf v_{r,u}+B_{r,u}^s$ and
$y_u:=\phi_u^{s,x,\mathbf v}-\zeta_u$.
Then $y$ is absolutely continuous on $[r,t]$, with
$\frac{d}{du}y_u=F(\phi_u^{s,x,\mathbf v})$ for almost every $u\in[r,t]$.
Writing $\phi_u:=\phi_u^{s,x,\mathbf v}=y_u+\zeta_u$, the coercivity and
growth assumptions on $F$ yield
\begin{align*}
 \frac12\frac{d}{du}|y_u|^2
 =
 \langle \phi_u,F(\phi_u)\rangle
 -\langle \zeta_u,F(\phi_u)\rangle
 \le
 -\kappa_1|\phi_u|^{k+1}
 +C+C|\zeta_u|(1+|\phi_u|^k).
\end{align*}
By Young's inequality,
$C|\zeta_u|(1+|\phi_u|^k)\le\frac{\kappa_1}{2}|\phi_u|^{k+1}
+C(1+|\zeta_u|^{k+1})$. Since
$|\phi_u|^{k+1}\ge c|y_u|^{k+1}-C|\zeta_u|^{k+1}$ and
$|\zeta_u|\le S_I(u-r)^{H-\eps}\le S_Ih^{H-\eps}$, we obtain, for some
$a>0$,
\begin{align*}
 \frac{d}{du}|y_u|^2
 \le
 -a|y_u|^{k+1}
 +C\bigl(1+S_I^{k+1}h^{(H-\eps)(k+1)}\bigr).
\end{align*}

We use the following elementary scalar comparison estimate.
If a nonnegative absolutely continuous function $\mathfrak q$ satisfies
$\mathfrak q'(u)\le-a\mathfrak q(u)^{\frac{k+1}{2}}+A$
almost everywhere on $[r,t]$, where $t-r=h$ and $A\ge0$,
then $\mathfrak q(t)\lesssim h^{-\frac2{k-1}}+A^{\frac2{k+1}}$.
Indeed, set $Q:=(2A/a)^{2/(k+1)}$. If $\mathfrak q(t)\le Q$,
the estimate is immediate. Otherwise, $\mathfrak q(u)>Q$ for every
$u\in[r,t]$: once $\mathfrak q$ reaches a value at most $Q$, it cannot
subsequently exceed $Q$, because $(\mathfrak q-Q)_+' \le0$ almost
everywhere. Hence $\mathfrak q'\le-\frac a2\mathfrak q^{\frac{k+1}{2}}$
almost everywhere on $[r,t]$, so that
$\frac{d}{du}\mathfrak q(u)^{-\frac{k-1}{2}}\ge\frac{a(k-1)}4$, and
integrating gives
$\mathfrak q(t)^{-\frac{k-1}{2}}\ge\frac{a(k-1)}4h$, which proves the
comparison estimate.

Applying it to $\mathfrak q(u)=|y_u|^2$ yields
$|y_t|\lesssim 1+h^{-\frac1{k-1}}+h^{H-\eps}S_I$.
Since $\phi_t=y_t+\zeta_t$ and $|\zeta_t|\le S_Ih^{H-\eps}$, this proves
\eqref{eq:coming-down-h}. The argument is pathwise and applies to every
admissible $t$ and $h$ on the same full-probability event.

We now deduce \eqref{eq:coming-down} by choosing $h$
according to the individual local scales. Write
$\mathfrak s_0:=[B^s]_{C^{H-\eps}(I)}$,
$\mathfrak s_i:=[\mathbf v^i]_{C^{H-\eps}(I)}$ for $i=1,\ldots,L$, and
$\tau_i:=(1\vee \mathfrak s_i)^{-1/\chi_\eps}$ for $i=0,\ldots,L$, and choose,
pathwise, $h:=(t-s)\wedge\tau_0\wedge\cdots\wedge\tau_L$.
For each $i$, if $t-s\le\tau_i$, then
$\mathfrak s_i h^{H-\eps}\le\mathfrak s_i(t-s)^{H-\eps}\le(t-s)^{-\frac1{k-1}}$,
by definition of $\tau_i$. If $t-s>\tau_i$, then $\mathfrak s_i>1$, since
$t-s\le1$, and
$\mathfrak s_i h^{H-\eps}\le\mathfrak s_i\tau_i^{H-\eps}
=\mathfrak s_i^{\frac1{(k-1)\chi_\eps}}$. Also,
\begin{align*}
 h^{-\frac1{k-1}}
 \le
 (t-s)^{-\frac1{k-1}}
 +\sum_{i=0}^L
 \ind_{\{t-s>\tau_i\}}
 \mathfrak s_i^{\frac1{(k-1)\chi_\eps}}.
\end{align*}
Substituting these bounds into \eqref{eq:coming-down-h} gives
\eqref{eq:coming-down}. Since $1/((k-1)\chi_\eps)\in(0,1)$, the constant and
the input terms in \eqref{eq:coming-down} are bounded by a constant multiple
of $(R_{s,t}^{B^s,\mathbf v^1,\ldots,\mathbf v^L;I})^{1/((k-1)\chi_\eps)}$.
This proves \eqref{eq:coming-down-R}.
\end{proof}

For the remainder of this section we fix one deterministic input
$\mathbf v\in C^{H-\eps}(I;\R^d)$, as explained after
\Cref{prop:flow-estimates}. The next consequence of
\Cref{lem:coming-down-infinity} is the form used for derivatives of $F$.

\begin{corollary}[Integrated polynomial control]
\label{cor:integrated-flow-control}
In the setting of \Cref{lem:coming-down-infinity}, for every
$\nu\in[0,k-1)$ there is a deterministic $C_\nu<\infty$ such that
\begin{align}\label{eq:integrated-flow-control}
 \sup_{x\in\R^d}\int_s^t(1+|\phi_r^{s,x,\mathbf v}|)^\nu\,dr
 \le C_\nu\Bigg(
 \frac{(t-s)^{1-\nu/(k-1)}}{1-\nu/(k-1)}
 +(t-s)(R_{s,t}^{B^s,\mathbf v;I})^{\frac{\nu}{(k-1)\chi_\eps}}
 \Bigg).
\end{align}
\end{corollary}

\begin{proof}
For $\nu=0$ the left-hand side is $t-s$. For $\nu>0$, apply
\eqref{eq:coming-down-R} at time $r\in(s,t]$ and note that
$R_{s,r}^{B^s,\mathbf v;I}\le R_{s,t}^{B^s,\mathbf v;I}$, since the
seminorms are taken on the same fixed interval $I$. Thus
$(1+|\phi_r^{s,x,\mathbf v}|)^\nu
\le C_\nu((r-s)^{-\nu/(k-1)}
+(R_{s,t}^{B^s,\mathbf v;I})^{\nu/((k-1)\chi_\eps)})$.
Integration gives the displayed estimate; the singular power is integrable
precisely for $\nu<k-1$.
\end{proof}

\begin{lemma}[Dependence on the initial condition]
\label{lem:reg-phi}
For every interval $I=[s,\bar t]$, every $t\in I$,
every deterministic input $\mathbf v\in C(I;\R^d)$, and every
$x,y\in\R^d$,
\begin{align}\label{eq:pathwise-x-Lipschitz}
\begin{split}
 |\phi_t^{s,x,\mathbf v}-\phi_t^{s,y,\mathbf v}|
 \le e^{\kappa_4(t-s)}|x-y|.
\end{split}
\end{align}
\end{lemma}

\begin{proof}
The additive inputs cancel in
$\mathfrak d_t:=\phi_t^{s,x,\mathbf v}-\phi_t^{s,y,\mathbf v}$, so
$\mathfrak d$ is absolutely continuous and, by
\eqref{eq:F-one-sided-global},
$\frac{d}{dt}|\mathfrak d_t|^2
=2\langle\mathfrak d_t,F(\phi_t^{s,x,\mathbf v})-F(\phi_t^{s,y,\mathbf v})\rangle
\le2\kappa_4|\mathfrak d_t|^2$.
Gr\"onwall's inequality gives \eqref{eq:pathwise-x-Lipschitz}.
\end{proof}

\subsection{Malliavin estimates for a single flow}

The pathwise estimates obtained above control the coefficients in the
Malliavin variational equations. We now derive pointwise and Hilbert-space
bounds for the Malliavin derivatives of a single flow.

Throughout this subsection, $I=[s,\bar t]$ is an interval of length at most
one, $t\in(s,\bar t]$, and $\mathbf v\in C^{H-\eps}(I;\R^d)$. We write
$D^\ell_{\mathbf u}\phi_t^{s,x,\mathbf v}$ for the jointly measurable version
selected in \Cref{lem:DY,prop:higher-DY}. These versions can be chosen
jointly measurable in $(x,\mathbf u,t,\omega)$; this follows from the
parameter-dependent variational equations and the smooth dependence of the
deterministic flow on $x$. All pointwise identities involving Malliavin
derivatives are understood on one full-probability set, independent of $x$,
$t$, and the Malliavin time variables. All the estimates below are uniform
in the initial condition $x$. We denote by $J_{u,r}$, $s\le u\le r\le t$, the
linearised flow along $\phi^{s,x,\mathbf v}$ introduced before
\Cref{lem:DY}; by \eqref{eq:J-Y-bound} and $t-s\le1$,
$|J_{u,r}|\le e^{\kappa_4(r-u)}\le e^{\kappa_4}$.

The Malliavin derivatives are controlled through the following integrated
form of the coming-down-from-infinity estimate. For $u,u'>s$, put
$\varrho_s(u,u'):=|u-u'|+|\log((u-s)/(u'-s))|$.

\begin{lemma}[Local integrated control]\label{lem:local-endpoint-control}
For every $s<a\le b\le t$,
\begin{align}\label{eq:local-endpoint-control}
\begin{split}
 \int_a^b\bigl(1+|\phi_r^{s,x,\mathbf v}|^{k-1}\bigr)\,dr
 \le C_1(R_{s,t}^{B^s,\mathbf v;I})^{C_2}\,\varrho_s(a,b).
\end{split}
\end{align}
\end{lemma}

\begin{proof}
By \eqref{eq:coming-down}, after enlarging $\mathbf C$ if necessary,
\begin{align*}
 1+|\phi_r^{s,x,\mathbf v}|^{k-1}
 \le \mathbf C\bigl(1+(r-s)^{-1}\bigr),
 \qquad r\in(s,t].
\end{align*}
Integrating this estimate over $[a,b]$ gives
\begin{align*}
 \int_a^b\bigl(1+|\phi_r^{s,x,\mathbf v}|^{k-1}\bigr)\,dr
 &\le \mathbf C\left(
 |b-a|+\left|\log\frac{b-s}{a-s}\right|
 \right),
\end{align*}
which is \eqref{eq:local-endpoint-control}.
\end{proof}

\paragraph{First Malliavin derivative.}

\begin{lemma}\label{lem:deriv1}
For every $u\in(s,t]$, one has $D_u\phi_t^{s,x,\mathbf v}=J_{u,t}$, and
hence
\begin{align}\label{eq:first-D-bound}
\begin{split}
 |D_u\phi_t^{s,x,\mathbf v}|\le e^{\kappa_4}.
\end{split}
\end{align}
Moreover, for every $u,u'\in(s,t]$,
\begin{align}\label{eq:first-D-time-increment}
\begin{split}
 |D_u\phi_t^{s,x,\mathbf v}-D_{u'}\phi_t^{s,x,\mathbf v}|
 \le C_1(R_{s,t}^{B^s,\mathbf v;I})^{C_2}\,\varrho_s(u,u').
\end{split}
\end{align}
The same estimates hold for an affine combination
$\Phi_t^{s,\mathbf x,\lambda,\mathbf v}$, with an additional factor
$\|\lambda\|_{\ell^1}$ on the right-hand side.
\end{lemma}

\begin{proof}
The identity $D_u\phi_t^{s,x,\mathbf v}=J_{u,t}$ is \eqref{eq:DY-stmt},
and \eqref{eq:first-D-bound} follows from \eqref{eq:J-Y-bound}.

It remains to prove \eqref{eq:first-D-time-increment}.  By symmetry, assume
that $u'\le u$.  Subtracting the two linear equations satisfied by
$D_u\phi^{s,x,\mathbf v}$ and $D_{u'}\phi^{s,x,\mathbf v}$ gives
\begin{align*}
 D_u\phi_t-D_{u'}\phi_t
 &=\int_u^t\nabla F(\phi_r)
   (D_u\phi_r-D_{u'}\phi_r)\,dr
   -\int_{u'}^u\nabla F(\phi_r)D_{u'}\phi_r\,dr,
\end{align*}
where the common superscripts have been omitted.  The variation-of-constants
formula, \eqref{eq:J-Y-bound}, and
\eqref{eq:first-D-bound} therefore yield
\begin{align*}
 |D_u\phi_t-D_{u'}\phi_t|
 &\lesssim\int_{u'}^u|\nabla F(\phi_r)|\,dr \lesssim\int_{u'}^u\bigl(1+|\phi_r|^{k-1}\bigr)\,dr.
\end{align*}
Applying \Cref{lem:local-endpoint-control} proves
\eqref{eq:first-D-time-increment}.  The assertion for affine combinations
follows from linearity of the Malliavin derivative and the triangle
inequality.
\end{proof}

\paragraph{Higher-order derivatives.}

The first-derivative estimate provides the base case for an induction over the
Malliavin order. For $\ell\ge2$ and
$\mathbf u=(u_1,\ldots,u_\ell)\in(s,t]^\ell$, define
\begin{align}\label{eq:higher-D-source}
\begin{split}
 \mathcal G_{\mathbf u}^{\ell}(r)
 :=\sum_{\substack{\pi\in\mathcal P_\ell\\|\pi|\ge2}}
 \nabla^{|\pi|}F(\phi_r^{s,x,\mathbf v})
 \left[
 D_{\mathbf u_A}^{|A|}\phi_r^{s,x,\mathbf v}:A\in\pi
 \right].
\end{split}
\end{align}
As usual, a Malliavin derivative is understood to vanish when one of its time
arguments is larger than the observation time.  The higher-order derivative
formula from \Cref{prop:higher-DY} may then be written as
\begin{align}\label{eq:higher-D-variation}
\begin{split}
 D_{\mathbf u}^{\ell}\phi_t^{s,x,\mathbf v}
 =\int_{\bar u}^tJ_{r,t}\mathcal G_{\mathbf u}^{\ell}(r)\,dr.
\end{split}
\end{align}

\begin{lemma}[Pointwise higher-order bounds]\label{lem:highderiv}
For every $\ell\in\{2,\ldots,\lfloor-\gamma\rfloor+5\}$,
\begin{align*}
 \sup_{\mathbf u\in(s,t]^\ell}
 |D_{\mathbf u}^{\ell}\phi_t^{s,x,\mathbf v}|
 \le C_1(R_{s,t}^{B^s,\mathbf v;I})^{C_2}.
\end{align*}
For an affine combination, the right-hand side is multiplied by
$\|\lambda\|_{\ell^1}$.
\end{lemma}

\begin{proof}
We argue by induction on $\ell$.  The first-order estimate is
\eqref{eq:first-D-bound}.  Assume that the assertion holds up to order
$\ell-1$.  Every partition occurring in
\eqref{eq:higher-D-source} has at least two blocks, and hence each derivative
appearing there has order at most $\ell-1$.  By the induction hypothesis and
\eqref{eq:growth},
\begin{align*}
 |\mathcal G_{\mathbf u}^{\ell}(r)|
 \le \mathbf C\bigl(1+|\phi_r^{s,x,\mathbf v}|^{((k-2)\vee0)}\bigr).
\end{align*}
Indeed, $|\pi|\ge2$ implies
$(k-|\pi|)\vee0\le((k-2)\vee0)$.  Combining this estimate with
\eqref{eq:higher-D-variation} and
\eqref{eq:J-Y-bound}, we obtain
\begin{align*}
 |D_{\mathbf u}^{\ell}\phi_t^{s,x,\mathbf v}|
 &\le \mathbf C\int_{\bar u}^t
 \bigl(1+|\phi_r^{s,x,\mathbf v}|^{((k-2)\vee0)}\bigr)\,dr
 \le \mathbf C,
\end{align*}
where the last bound follows from
\Cref{cor:integrated-flow-control}, since $((k-2)\vee0)<k-1$.  This closes the
induction.  The affine estimate again follows by linearity and the triangle
inequality.
\end{proof}

\paragraph{Malliavin Hilbert-space bounds.}

The tensor norms are estimated directly from the Hilbert-valued variational
equations, so that no increment estimate in the Malliavin time variables is
needed beyond the first-order one, \eqref{eq:first-D-time-increment}. The
latter is used in the base case $\ell=1$, through the following
deterministic lemma, and in the non-degeneracy argument below.

\begin{lemma}\label{lem:H-log-modulus}
Let $t-s\le1$, let $E$ be a finite-dimensional Hilbert space, and let
$f:(s,t]\to E$ satisfy $\sup_{u\in(s,t]}|f(u)|\le A$ and
$|f(u)-f(u')|\le A\varrho_s(u,u')$ for all $u,u'\in(s,t]$.
Extend $f$ by zero outside $[s,t]$.  Then
\begin{align}\label{eq:H-log-modulus}
\begin{split}
 \|f\|_{\mathcal H_{s,t}}\le C_HA(t-s)^H.
\end{split}
\end{align}
\end{lemma}

\begin{proof}
If $H=1/2$, then \eqref{eq:brownian-canonical-space} gives
$\|f\|_{\mathcal H_{s,t}}=\|f\|_{L_2([s,t])}\le A(t-s)^{1/2}$.
If $H>1/2$, \eqref{eq:smooth-H-embedding}, applied on $[s,t]$, gives
$\|f\|_{\mathcal H_{s,t}}\le C_H(t-s)^{H-\frac12}\|f\|_{L_2([s,t])}
\le C_HA(t-s)^H$.  We may therefore assume in the
remainder of the proof that $H<1/2$.
Using the representation of $\mathfrak K^*$ as a fractional derivative
(see \eqref{eq:H-fractional-norm} and the proof of \Cref{lem:Hinterp}) for a
function supported on $[s,t]$, for $u\in(s,t)$ we have
\begin{align*}
 |(\mathfrak K^*f)(u)|
 &\lesssim |f(u)|(t-u)^{H-\frac12}
 +\int_u^t|f(r)-f(u)|(r-u)^{H-\frac32}\,dr.
\end{align*}
The increment assumption gives
\begin{align}\label{eq:log-kernel-integral}
\begin{split}
 \int_u^t\varrho_s(r,u)(r-u)^{H-\frac32}\,dr
 \lesssim
 (t-u)^{H-\frac12}+(u-s)^{H-\frac12}.
\end{split}
\end{align}
Indeed, the contribution of $r-u$ is immediate.  For the logarithmic term,
set $a=u-s$ and $z=r-u$, split the integral at $z=a$, use
$\log(1+z/a)\le z/a$ for $z\le a$, and use
$\log(1+z/a)\lesssim(z/a)^\iota$ for $z>a$, with
$0<\iota<\frac12-H$.  This proves \eqref{eq:log-kernel-integral}.
Consequently,
$|(\mathfrak K^*f)(u)|\lesssim A((t-u)^{H-\frac12}+(u-s)^{H-\frac12})$.
Both singularities are square-integrable because $H>0$.  Integrating over
$u\in(s,t)$ and using the isometry defining $\mathcal H_{s,t}$ gives
\eqref{eq:H-log-modulus}.
\end{proof}

\begin{lemma}[Tensor Malliavin bounds]\label{lem:single-flow-Hnorm}
For every $\ell\in\{1,\ldots,\lfloor-\gamma\rfloor+5\}$,
\begin{align}\label{eq:single-flow-Hnorm-pathwise}
\begin{split}
 \|D^\ell\phi_t^{s,x,\mathbf v}\|_{\mathcal H_{s,t}^{\otimes\ell}}
 \le C_1(R_{s,t}^{B^s,\mathbf v;I})^{C_2}(t-s)^{\ell H}.
\end{split}
\end{align}
Consequently, for every $p\in[1,\infty)$,
\begin{align}\label{eq:single-flow-Hnorm-Lp}
\begin{split}
 \big\|D^\ell\phi_t^{s,x,\mathbf v}
 \big\|_{L_p(\mathcal H_{s,t}^{\otimes\ell})}
 \le C_1
 \bigl(R_{s,t}^{\mathbf v;I}\bigr)^{C_2}
 (t-s)^{\ell H}.
\end{split}
\end{align}
For an affine combination, both right-hand sides are multiplied by
$\|\lambda\|_{\ell^1}$.
\end{lemma}

\begin{proof}
If $H>1/2$, the pointwise estimates
\eqref{eq:first-D-bound} and \Cref{lem:highderiv}, followed by
\eqref{eq:smooth-H-tensor-embedding}, give directly
\begin{align*}
 \|D^\ell\phi_t\|_{\mathcal H_{s,t}^{\otimes\ell}}
 &\le C_H^\ell(t-s)^{\ell(H-\frac12)}
 \|D^\ell\phi_t\|_{L_2((s,t)^\ell)} \le\mathbf C(t-s)^{\ell H}.
\end{align*}
Suppose now that $H\le1/2$.  For $\ell=1$, combine
\Cref{lem:deriv1} with \Cref{lem:H-log-modulus}.  Assume that
\eqref{eq:single-flow-Hnorm-pathwise} holds for all orders smaller than
$\ell$, where $\ell\ge2$.  The tensor-valued version of
\eqref{eq:higher-D-variation}, obtained from the Malliavin chain rule and the
variation-of-constants formula, gives
\begin{align*}
 \|D^\ell\phi_t\|_{\mathcal H^{\otimes\ell}}
 &\le e^{\kappa_4}
 \sum_{\substack{\pi\in\mathcal P_\ell\\|\pi|\ge2}}
 \int_s^t
 |\nabla^{|\pi|}F(\phi_r)|
 \prod_{A\in\pi}
 \|D^{|A|}\phi_r\|_{\mathcal H^{\otimes|A|}}\,dr.
\end{align*}
The induction hypothesis and the identity
$\sum_{A\in\pi}|A|=\ell$ imply
\begin{align*}
 \|D^\ell\phi_t\|_{\mathcal H^{\otimes\ell}}
 &\le \mathbf C\int_s^t
 \bigl(1+|\phi_r|^{((k-2)\vee0)}\bigr)(r-s)^{\ell H}\,dr.
\end{align*}
By \eqref{eq:coming-down},
$1+|\phi_r|^{(k-2)\vee0}\le\mathbf C(1+(r-s)^{-((k-2)\vee0)/(k-1)})$, and,
since $t-s\le1$,
$\int_s^t(1+(r-s)^{-((k-2)\vee0)/(k-1)})(r-s)^{\ell H}\,dr
\lesssim (t-s)^{\ell H}$.
This proves \eqref{eq:single-flow-Hnorm-pathwise} in the remaining cases.  The moment estimate
\eqref{eq:single-flow-Hnorm-Lp} follows from
\eqref{eq:bold-C-control} and the Gaussian moments of
$[B^s]_{C^{H-\eps}(I)}$.  The affine estimate follows by linearity.
\end{proof}

\subsection{Spatial regularity of the Malliavin hierarchy}

The preceding subsection controls the Malliavin hierarchy for a fixed initial
condition. The sewing arguments also compare flows started from different
points, so we next establish the corresponding spatial and four-point bounds.
They are obtained by differentiating the Hilbert-valued Malliavin
variational equations with respect to the initial condition and arguing
inductively in the total number of derivatives. The Malliavin variables are
distributed among disjoint blocks of a partition, so the Hilbert tensor norm
factorises and each term is bounded by the product of lower-order tensor
norms.

For $n\in\{1,2\}$ and $\ell\ge0$, we write
$\nabla_x^nD_{\mathbf u}^{\ell}\phi_t^{s,x,\mathbf v}$ for the $n$-th
Fr\'echet derivative with respect to $x$, viewed as an $n$-linear map.  When
$\ell=0$, the Malliavin variables are omitted.
The pointwise kernels of \Cref{lem:DY,prop:higher-DY} are given by
\eqref{eq:DY-stmt} and the recursion \eqref{eq:higher-DY-partitions} for
every $\omega$ such that $B^s(\omega)$ is continuous; this full-probability
set does not depend on $x$. For such $\omega$, $\phi^{s,x,\mathbf v}$ and $J$
are solutions of ordinary differential equations with smooth coefficients
and depend smoothly on $x$, with locally uniform bounds on their
$x$-derivatives on $I$. Hence the kernels are smooth in $x$, and their
$x$-derivatives are
obtained by differentiating \eqref{eq:DY-stmt} and
\eqref{eq:higher-DY-partitions} under the integral sign; this is the mixed
chain rule used below.

\begin{lemma}[Mixed spatial-Malliavin derivatives]
\label{lem:mixed-spatial-malliavin}
Let $n\in\{1,2\}$ and $\ell\in\N\cup\{0\}$ satisfy
$n+\ell\le\lfloor-\gamma\rfloor+5$. Then, for all $x\in\R^d$,
\begin{align}
 \sup_{\mathbf u\in(s,t]^\ell}
 \big|\nabla_x^nD_{\mathbf u}^{\ell}
 \phi_t^{s,x,\mathbf v}\big|
 &\le C_1(R_{s,t}^{B^s,\mathbf v;I})^{C_2},
 \label{eq:mixed-spatial-pointwise}
\end{align}
and, if $\ell\ge1$,
\begin{align}
 \big\|\nabla_x^nD^{\ell}
 \phi_t^{s,x,\mathbf v}\big\|_{\mathcal H_{s,t}^{\otimes\ell}}
 &\le C_1(R_{s,t}^{B^s,\mathbf v;I})^{C_2}(t-s)^{\ell H}.
 \label{eq:mixed-spatial-H}
\end{align}
\end{lemma}

\begin{proof}
We first consider $\ell=0$.  The Jacobian $\nabla_x\phi_t^{s,x,\mathbf v}$
is the linearised flow $J_{s,t}$, so $|\nabla_x\phi_t^{s,x,\mathbf v}|\le
e^{\kappa_4(t-s)}$ by \eqref{eq:J-Y-bound}. The second spatial derivative
$J^{(2)}_{s,t}:=\nabla_x^2\phi_t^{s,x,\mathbf v}$ satisfies
\begin{align*}
 \partial_tJ^{(2)}_{s,t}
 =\nabla F(\phi_t^{s,x,\mathbf v})J^{(2)}_{s,t}
   +\nabla^2F(\phi_t^{s,x,\mathbf v})
       [J_{s,t},J_{s,t}],
 \qquad
 J^{(2)}_{s,s}=0.
\end{align*}
Variation of constants, the bound on $J$, the growth assumption on
$\nabla^2F$, and \Cref{cor:integrated-flow-control} yield
\begin{align*}
 |J^{(2)}_{s,t}|
 &\le \mathbf C\int_s^t
 \bigl(1+|\phi_r^{s,x,\mathbf v}|^{(k-2)\vee0}\bigr)\,dr
 \le \mathbf C,
\end{align*}
since $(k-2)\vee0<k-1$.

We now argue by induction on the total order $n+\ell$.  After isolating the
term in which all derivatives fall on the inner function, the mixed chain
rule gives, for $\ell\ge1$,
\begin{align}
 \nabla_x^nD_{\mathbf u}^{\ell}
 \phi_t^{s,x,\mathbf v}
 =\int_{\bar u}^t
 J_{r,t}\,\mathcal G_{\mathbf u}^{n,\ell}(r)\,dr,
 \label{eq:mixed-spatial-voc}
\end{align}
where the remainder $\mathcal G_{\mathbf u}^{n,\ell}$ is a finite sum of
terms of the form
\begin{align}
 \nabla^\mu F(\phi_r^{s,x,\mathbf v})
 \big[
   \nabla_x^{n_1}D_{\mathbf u_{A_1}}^{|A_1|}\phi_r,
   \ldots,
   \nabla_x^{n_\mu}D_{\mathbf u_{A_\mu}}^{|A_\mu|}\phi_r
 \big],
 \label{eq:mixed-spatial-remainder}
\end{align}
with $\mu\ge2$, pairwise disjoint and possibly empty sets
$A_1,\ldots,A_\mu$ with union $\{1,\ldots,\ell\}$, and integers
$n_1+\cdots+n_\mu=n$ with $n_i+|A_i|\ge1$. The combinatorial
coefficients depend only on $n$ and $\ell$.  Since $\mu\ge2$, every factor
in \eqref{eq:mixed-spatial-remainder} has strictly smaller total order than
$n+\ell$.

By the induction hypothesis, all factors in
\eqref{eq:mixed-spatial-remainder} are pointwise bounded by $\mathbf C$.
Moreover, since $\mu\ge2$, \eqref{eq:growth} and \eqref{eq:coming-down} give
$|\nabla^\mu F(\phi_r)|\lesssim 1+|\phi_r|^{(k-2)\vee0}
\le \mathbf C(1+(r-s)^{-((k-2)\vee0)/(k-1)})$.
Substituting this estimate into \eqref{eq:mixed-spatial-voc} and using
$|J_{r,t}|\le e^{\kappa_4(t-r)}$ proves
\eqref{eq:mixed-spatial-pointwise}.

If $H>1/2$, \eqref{eq:smooth-H-tensor-embedding} and the pointwise estimate
just proved yield
\begin{align*}
 \big\|\nabla_x^nD^\ell\phi_t
 \big\|_{\mathcal H_{s,t}^{\otimes\ell}}
 &\le C_H^\ell(t-s)^{\ell(H-\frac12)}
 \big\|\nabla_x^nD^\ell\phi_t
 \big\|_{L_2((s,t)^\ell)} \le\mathbf C(t-s)^{\ell H}.
\end{align*}
We may therefore assume that $H\le1/2$.
In that case, take the $\mathcal H_{s,t}^{\otimes\ell}$-norm in
\eqref{eq:mixed-spatial-voc}.  In each term of
\eqref{eq:mixed-spatial-remainder}, the Malliavin variables are distributed
among the disjoint sets $A_1,\ldots,A_\mu$.  The tensor-product property of
the Hilbert norm and the induction hypothesis give
\begin{align*}
 \left\|
 \prod_{i=1}^\mu
 \nabla_x^{n_i}D_{\mathbf u_{A_i}}^{|A_i|}\phi_r
 \right\|_{\mathcal H_{s,r}^{\otimes\ell}} &\le
 \mathbf C
 \prod_{i:A_i\ne\varnothing}(r-s)^{|A_i|H}
 =\mathbf C(r-s)^{\ell H}.
\end{align*}
Consequently, using $t-s\le1$,
\begin{align*}
 \big\|\nabla_x^nD^{\ell}\phi_t
 \big\|_{\mathcal H_{s,t}^{\otimes\ell}}
 &\le \mathbf C\int_s^t
 \bigl(1+(r-s)^{-\frac{(k-2)\vee0}{k-1}}\bigr)
 (r-s)^{\ell H}\,dr \le \mathbf C(t-s)^{\ell H}.
\end{align*}
This closes the induction.
\end{proof}

\begin{corollary}[Spatial finite differences]
\label{cor:spatial-finite-differences}
For every $\ell\in\{0,\ldots,\lfloor-\gamma\rfloor+4\}$ and $x,y\in\R^d$,
\begin{align}
 \big\|D^{\ell}
 (\phi_t^{s,x,\mathbf v}-\phi_t^{s,y,\mathbf v})
 \big\|_{\mathcal H_{s,t}^{\otimes\ell}}
 \le C_1(R_{s,t}^{B^s,\mathbf v;I})^{C_2}|x-y|(t-s)^{\ell H},
 \label{eq:spatial-difference-H}
\end{align}
where, for $\ell=0$, the Hilbert norm and the factor
$(t-s)^{\ell H}$ are interpreted as the Euclidean norm and $1$,
respectively.
\end{corollary}

\begin{proof}
By the fundamental theorem of calculus in the initial condition,
$D^{\ell}\phi_t^{s,x,\mathbf v}-D^{\ell}\phi_t^{s,y,\mathbf v}
=\int_0^1\nabla_xD^{\ell}\phi_t^{s,y+\theta(x-y),\mathbf v}[x-y]\,d\theta$.
This identity holds pointwise in the Malliavin variables, since the kernels
are smooth in $x$. For $\ell\ge1$, the norm of $\mathcal H_{s,t}^{\otimes\ell}$
is the $L_2$-norm of the image under $(\mathfrak K^*)^{\otimes\ell}$, which
commutes with the $\theta$-integral, so Minkowski's integral inequality and
\eqref{eq:mixed-spatial-H} give \eqref{eq:spatial-difference-H}; the case $\ell=0$ follows from
\eqref{eq:mixed-spatial-pointwise}.
\end{proof}

For $x_1,x_2,x_3,x_4\in\R^d$, we write
$\square_{x_1,x_2,x_3,x_4}\phi_t^{s,\mathbf v}
:=\phi_t^{s,x_1,\mathbf v}-\phi_t^{s,x_2,\mathbf v}
-\phi_t^{s,x_3,\mathbf v}+\phi_t^{s,x_4,\mathbf v}$
for the four-point increment of the flow.

\begin{corollary}[Four-point Malliavin estimates]
\label{cor:four-point-malliavin}
Let $x_1,x_2,x_3,x_4\in\R^d$ and
$\mathfrak e(\mathbf x):=|x_1-x_2-x_3+x_4|+|x_1-x_2|\,|x_1-x_3|$.
Then, for every $\ell\in\{0,\ldots,\lfloor-\gamma\rfloor+3\}$,
\begin{align}
 \big\|D^{\ell}
 \square_{x_1,x_2,x_3,x_4}\phi_t^{s,\mathbf v}
 \big\|_{\mathcal H_{s,t}^{\otimes\ell}}
 \le C_1(R_{s,t}^{B^s,\mathbf v;I})^{C_2}\mathfrak e(\mathbf x)(t-s)^{\ell H},
 \label{eq:four-point-H}
\end{align}
with the same convention as above when $\ell=0$.
\end{corollary}

\begin{proof}
Let $\mathfrak f:\R^d\to E$ be a $C^2$ map with values in a Banach space $E$,
and set $z:=x_2+x_3-x_1$.  The fundamental theorem of calculus gives
\begin{align*}
 &\mathfrak f(x_1)-\mathfrak f(x_2)-\mathfrak f(x_3)+\mathfrak f(x_4)
 \\
 &\quad=
 \int_0^1\int_0^1
 \nabla^2\mathfrak f
 \bigl(x_1-\theta(x_1-x_2)-\vartheta(x_1-x_3)\bigr)
 [x_1-x_2,x_1-x_3] \,d\theta\,d\vartheta
 +\mathfrak f(x_4)-\mathfrak f(z).
\end{align*}
Consequently,
\begin{align}
 \|\mathfrak f(x_1)-\mathfrak f(x_2)-\mathfrak f(x_3)+\mathfrak f(x_4)\|_E
 \le
 \|\nabla \mathfrak f\|_{\infty}|x_1-x_2-x_3+x_4|
 +\|\nabla^2\mathfrak f\|_{\infty}|x_1-x_2|\,|x_1-x_3|.
 \label{eq:banach-four-point}
\end{align}
Apply the preceding identity pointwise in the Malliavin variables to the
kernels of $D^{\ell}\phi_t^{s,x,\mathbf v}$, which are smooth in $x$. As
in the proof of \Cref{cor:spatial-finite-differences}, Minkowski's integral
inequality then gives \eqref{eq:banach-four-point} with
$E=\mathcal H_{s,t}^{\otimes\ell}$ and
$\mathfrak f(x)=D^{\ell}\phi_t^{s,x,\mathbf v}$.  Estimate
\eqref{eq:mixed-spatial-H} proves \eqref{eq:four-point-H}; for $\ell=0$ use
\eqref{eq:mixed-spatial-pointwise}.
\end{proof}

\subsection{Non-degeneracy of affine combinations}

With the Malliavin and spatial estimates available, it remains to prove the
non-degeneracy required by the integration-by-parts formulas. The affine
constraint $\sum_{i=1}^M\lambda_i=1$ is used only in this
subsection: it ensures that the Malliavin derivative of the affine
combination at the terminal time is the identity matrix. We first show that
it remains close to the identity on a sufficiently short terminal
interval.

\begin{lemma}[Terminal control of the Malliavin derivative]
\label{lem:terminal-D-control}
Work in the setting of \Cref{prop:flow-estimates} with a single input.
There exist constants $C_1,C_2<\infty$, depending only on the parameters
listed in that proposition, such that, with
$\mathfrak A:=C_1(1+\|\lambda\|_{\ell^1})(R_{s,t}^{B^s,\mathbf v;I})^{C_2}$,
one has, for all $u,r\in[s+(t-s)/2,t]$,
\begin{align}
 \bigl\|D_u\Phi_t^{s,\mathbf x,\lambda,\mathbf v}-I_d\bigr\|
 &\le \mathfrak A\frac{t-u}{t-s},
 \label{eq:terminal-D-identity}\\
 \bigl\|D_r\Phi_t^{s,\mathbf x,\lambda,\mathbf v}
       -D_u\Phi_t^{s,\mathbf x,\lambda,\mathbf v}\bigr\|
 &\le \mathfrak A\frac{|r-u|}{t-s}.
 \label{eq:terminal-D-Lipschitz}
\end{align}
The estimates are uniform in the initial points $x_1,\ldots,x_M$.
\end{lemma}

\begin{proof}
Since $D_t\phi_t^{s,x_i,\mathbf v}=I_d$ and
$\sum_{i=1}^M\lambda_i=1$, we have
$D_t\Phi_t^{s,\mathbf x,\lambda,\mathbf v}=I_d$. By \eqref{eq:first-D-time-increment} and linearity,
\begin{align*}
 \bigl\|D_u\Phi_t-I_d\bigr\|
 &\le \|\lambda\|_{\ell^1}\mathbf C\,\varrho_s(u,t),\\
 \bigl\|D_r\Phi_t-D_u\Phi_t\bigr\|
 &\le \|\lambda\|_{\ell^1}\mathbf C\,\varrho_s(r,u).
\end{align*}
If $u,u'\in[s+(t-s)/2,t]$, then, using $t-s\le1$ and
$(u\wedge u')-s\ge(t-s)/2$,
\begin{align*}
 \varrho_s(u,u')
 &=|u-u'|+
 \left|\log\frac{u-s}{u'-s}\right| \le \frac{|u-u'|}{t-s}
 +\frac{|u-u'|}{(u\wedge u')-s}
 \le 3\frac{|u-u'|}{t-s}.
\end{align*}
Combining this estimate with \eqref{eq:bold-C-control} proves
\eqref{eq:terminal-D-identity} and \eqref{eq:terminal-D-Lipschitz}, after
increasing the deterministic constants in the definition of $\mathfrak A$.
\end{proof}

\begin{lemma}[Quantitative non-degeneracy]
\label{lem:quantitative-nondegeneracy}
In the setting of \Cref{prop:flow-estimates} with a single input, there exist
deterministic constants $c,C_2>0$ such that, almost surely,
\begin{align}
 y^\top\mathcal M_{\Phi_t^{s,\mathbf x,\lambda,\mathbf v}}y
 \ge
 c\bigl(1+\|\lambda\|_{\ell^1}\bigr)^{-C_2}
 \bigl(R_{s,t}^{B^s,\mathbf v;I}\bigr)^{-C_2}
 (t-s)^{2H}|y|^2
 \label{eq:pathwise-covariance-lower}
\end{align}
for every $y\in\R^d$.  In particular, the Malliavin covariance is invertible
almost surely.
\end{lemma}

\begin{proof}
Let $\mathfrak A$ be as in \Cref{lem:terminal-D-control}; by increasing its
deterministic prefactor, we may assume that $\mathfrak A\ge2$.  Fix a
constant $\theta_*\in(0,1/4)$, to be chosen below, and define
$h_*:=\theta_*(t-s)/\mathfrak A$ and $a:=t-h_*$, so that
$a\ge s+(t-s)/2$.  For $y\in\R^d$ and $u\in[s,t]$, set
$\mathfrak g_y(u):=(D_u\Phi_t^{s,\mathbf x,\lambda,\mathbf v})^\top y$ and
$\mathfrak r_y(u):=\mathfrak g_y(u)-y$.
By \eqref{eq:terminal-D-identity}--\eqref{eq:terminal-D-Lipschitz}, for
$u,r\in[a,t]$,
\begin{align}
 |\mathfrak r_y(u)|
 &\le \mathfrak A\frac{t-u}{t-s}|y|,
 \label{eq:terminal-error-size}\\
 |\mathfrak r_y(r)-\mathfrak r_y(u)|
 &\le \mathfrak A\frac{|r-u|}{t-s}|y|.
 \label{eq:terminal-error-Lipschitz}
\end{align}

If $H=1/2$, then $\mathfrak K^*=\mathrm{Id}$ and, for $u\in[a,t]$,
\begin{align*}
 |\mathfrak g_y(u)|
 \ge |y|-|\mathfrak r_y(u)|
 \ge(1-\theta_*)|y|.
\end{align*}
Taking, for instance, $\theta_*\le1/2$, we obtain
\begin{align*}
 y^\top\mathcal M_{\Phi_t}y
 =\int_s^t|\mathfrak g_y(u)|^2\,du
 \ge\frac14|y|^2h_*
 =\frac14|y|^2h_*^{2H}.
\end{align*}
Substituting the definition of $h_*$ and the bound on $\mathfrak A$
proves \eqref{eq:pathwise-covariance-lower} in the Brownian case.

Suppose next that $H>1/2$.  By
\eqref{eq:smooth-K-star}, for $u\in[a,t]$,
\begin{align*}
 (\mathfrak K^*\mathfrak g_y)(u)
 &=c_H\Bigl(H-\frac12\Bigr)\int_u^t\mathfrak g_y(r)(r-u)^{H-\frac32}\,dr =c_H(t-u)^{H-\frac12}y+(\mathfrak K^*\mathfrak r_y)(u).
\end{align*}
For $r\in[a,t]$, \eqref{eq:terminal-error-size} and the definition of
$h_*$ imply
\begin{align*}
 |\mathfrak r_y(r)|
 \le \mathfrak A\frac{h_*}{t-s}|y|
 =\theta_*|y|.
\end{align*}
Consequently,
\begin{align*}
 |(\mathfrak K^*\mathfrak r_y)(u)|
 &\le c_H\Bigl(H-\frac12\Bigr)\theta_*|y|
 \int_u^t(r-u)^{H-\frac32}\,dr =c_H\theta_*|y|(t-u)^{H-\frac12}.
\end{align*}
Taking $\theta_*\le1/2$ therefore yields
\begin{align*}
 |(\mathfrak K^*\mathfrak g_y)(u)|
 \ge\frac{c_H}{2}|y|(t-u)^{H-\frac12},
 \qquad u\in[a,t].
\end{align*}
It follows that
\begin{align*}
 y^\top\mathcal M_{\Phi_t}y
 &=\int_s^t|(\mathfrak K^*\mathfrak g_y)(u)|^2\,du \ge c_H|y|^2\int_a^t(t-u)^{2H-1}\,du
 =c_H'|y|^2h_*^{2H}.
\end{align*}
Substituting the definition of $h_*$ and the bound on $\mathfrak A$
proves \eqref{eq:pathwise-covariance-lower} for $H>1/2$.

We henceforth assume that $H<1/2$.

For $u\in[a,t]$, the right-sided representation of $\mathfrak K^*$ gives
\begin{align*}
 (\mathfrak K^*\mathfrak g_y)(u)
 &=c_H(t-u)^{H-\frac12}y+(\mathfrak K^*\mathfrak r_y)(u),
\end{align*}
where
\begin{align*}
 |(\mathfrak K^*\mathfrak r_y)(u)|
 &\lesssim_H
 |\mathfrak r_y(u)|(t-u)^{H-\frac12}
 +\int_u^t
 |\mathfrak r_y(r)-\mathfrak r_y(u)|(r-u)^{H-\frac32}\,dr.
\end{align*}
Using \eqref{eq:terminal-error-size} and
\eqref{eq:terminal-error-Lipschitz}, we obtain
\begin{align*}
 |(\mathfrak K^*\mathfrak r_y)(u)|
 &\lesssim_H
 \mathfrak A\frac{h_*}{t-s}|y|
 (t-u)^{H-\frac12} \le C_H\theta_* |y|(t-u)^{H-\frac12}.
\end{align*}
Choose $\theta_*$ so small that $C_H\theta_*\le c_H/2$.  It follows that
\begin{align*}
 |(\mathfrak K^*\mathfrak g_y)(u)|
 \ge \frac{c_H}{2}|y|(t-u)^{H-\frac12},
 \qquad u\in[a,t].
\end{align*}
Using the isometry between $\mathcal H_{s,t}$ and the range of $\mathfrak K^*$, we
therefore have
\begin{align*}
 y^\top\mathcal M_{\Phi_t}y
 &=\|\mathfrak g_y\|_{\mathcal H_{s,t}}^2
 =\int_s^t|(\mathfrak K^*\mathfrak g_y)(u)|^2\,du \ge c_H|y|^2\int_a^t(t-u)^{2H-1}\,du
 =c_H'|y|^2h_*^{2H}.
\end{align*}
Substituting the definition of $h_*$ and the bound on $\mathfrak A$ gives
\eqref{eq:pathwise-covariance-lower}.
\end{proof}

\begin{corollary}[Inverse covariance estimate]
\label{cor:inverse-covariance}
In the single-input setting of \Cref{prop:flow-estimates},
\begin{align}
 \bigl\|\mathcal M_{\Phi_t^{s,\mathbf x,\lambda,\mathbf v}}^{-1}\bigr\|_{L_p}
 \le C_1
 \bigl(R_{s,t}^{\mathbf v;I}\bigr)^{C_2}
 (t-s)^{-2H}.
 \label{eq:inverse-covariance-corollary}
\end{align}
\end{corollary}

\begin{proof}
The lower bound \eqref{eq:pathwise-covariance-lower} implies the pathwise
estimate
\begin{align*}
 \|\mathcal M_{\Phi_t}^{-1}\|_{\mathrm{op}}
 \le C_1(1+\|\lambda\|_{\ell^1})^{C_2}(R_{s,t}^{B^s,\mathbf v;I})^{C_2}
 (t-s)^{-2H}.
\end{align*}
Taking the $L_p$-norm and applying
\eqref{eq:flow-control-moments} proves \eqref{eq:inverse-covariance-corollary}.
\end{proof}

\subsection{Final proof of the regularising-flow estimates}

\begin{proof}[Proof of \Cref{prop:flow-estimates}]
We first collect the estimates for one input. Since
$0\le j\le\lfloor-\gamma\rfloor+3$, we have $j+2\le\lfloor-\gamma\rfloor+5$,
so every derivative used below is covered by the preceding lemmas.

\medskip\noindent\emph{Proof of \eqref{eq:upperbound-D}.}
By definition of the Hilbert-valued Malliavin--Sobolev norm,
\begin{align*}
 \big\|D\Phi_t^{s,\mathbf x,\lambda,\mathbf v}
 \big\|_{\mathbb D^{j,p}(\mathcal H_{s,t})}
 \lesssim
 \sum_{\ell=1}^{j+1}
 \big\|D^\ell\Phi_t^{s,\mathbf x,\lambda,\mathbf v}
 \big\|_{L_p(\mathcal H_{s,t}^{\otimes\ell})}.
\end{align*}
Applying \Cref{lem:single-flow-Hnorm} to each flow in the affine combination
gives
\begin{align*}
 \big\|D^\ell\Phi_t^{s,\mathbf x,\lambda,\mathbf v}
 \big\|_{L_p(\mathcal H_{s,t}^{\otimes\ell})}
 \le C_1\|\lambda\|_{\ell^1}
 \bigl(R_{s,t}^{\mathbf v;I}\bigr)^{C_2}
 (t-s)^{\ell H}.
\end{align*}
Since $t-s\le1$ and $\ell\ge1$, $(t-s)^{\ell H}\le(t-s)^H$.
Summing over $\ell=1,\ldots,j+1$ proves \eqref{eq:upperbound-D}.

\medskip\noindent\emph{Proof of \eqref{eq:inv-matrix}.}
This is \Cref{cor:inverse-covariance}.

\medskip\noindent\emph{Proof of \eqref{eq:upperbound-D-x}.}
For $\ell=0,\ldots,j$, \Cref{cor:spatial-finite-differences} gives
\begin{align*}
 \big\|D^{\ell}
 (\phi_t^{s,x,\mathbf v}-\phi_t^{s,y,\mathbf v})
 \big\|_{\mathcal H_{s,t}^{\otimes\ell}}
 \le \mathbf C|x-y|(t-s)^{\ell H}.
\end{align*}
Taking $L_p$-norms, using \eqref{eq:bold-C-control}, and summing over
$\ell=0,\ldots,j$ yields
\begin{align*}
 \big\|\phi_t^{s,x,\mathbf v}-\phi_t^{s,y,\mathbf v}
 \big\|_{\mathbb D^{j,p}}
 \le C_1
 \bigl(R_{s,t}^{\mathbf v;I}\bigr)^{C_2}|x-y|,
\end{align*}
because $t-s\le1$.  This is \eqref{eq:upperbound-D-x}.

\medskip\noindent\emph{Proof of \eqref{eq:upperbound-4}.}
For $\ell=0,\ldots,j$, \Cref{cor:four-point-malliavin} implies
\begin{align*}
 &\big\|D^{\ell}
 \square_{x_1,x_2,x_3,x_4}\phi_t^{s,\mathbf v}
 \big\|_{\mathcal H_{s,t}^{\otimes\ell}}
 \\
 &\qquad\le
 \mathbf C
 \Bigl(
 |x_1-x_2-x_3+x_4|
 +|x_1-x_2|\,|x_1-x_3|
 \Bigr)(t-s)^{\ell H}.
\end{align*}
Taking $L_p$-norms and summing over $\ell=0,\ldots,j$ proves
\eqref{eq:upperbound-4}.

\medskip\noindent\emph{Several inputs.}
Let now $\mathbf v=\sum_{i=1}^L\mathbf v^i$. By
\Cref{lem:coming-down-infinity}, the pathwise bounds hold with
$R_{s,t}^{B^s,\mathbf v;I}$ replaced by
$R_{s,t}^{B^s,\mathbf v^1,\ldots,\mathbf v^L;I}$, and
\Cref{cor:integrated-flow-control} uses only this bound and its monotonicity
in $t$ on the fixed interval $I$. The variational equations for the
Malliavin and spatial derivatives do not involve the inputs, so all
subsequent estimates carry the same factor, with constants depending also on
$L$. Finally, \eqref{eq:flow-control-moments} removes the fresh noise from
this factor when taking $L_p$ norms.
\end{proof}

\section{A priori estimates}\label{sec:a-priori}

We now transfer the regularising-flow estimates to solutions of the equation.
The first subsection treats one smooth equation and yields bounds uniform in
the smooth approximation of the singular drift. The second subsection compares
two solutions with singular drifts and is used for existence (through the
Cauchy property of the smooth approximations), uniqueness, and stability.

The argument is organised as follows. \Cref{prop:regsol-inf} compares a
solution of a smooth equation with the regularising flow below the local
scale, by stochastic sewing closed through a bootstrap over short windows;
\Cref{prop:regsol} transfers this comparison to the restarted, conditional
setting, which yields the class $\mathcal V_\eps(\infty)$ uniformly in the
smooth approximation. \Cref{lem:smooth-solution-moments} and
\Cref{prop:regsol-usual} add uniform moments and temporal regularity. In
\Cref{sec:stab}, \Cref{lem:stability-dissipativity} controls the deviation
of two solutions from their restarted flows in terms of the drift
difference, \Cref{lem:stability-sewing} controls the drift difference in
terms of this deviation, and \Cref{lem:stability-energy-sewing} provides the
energy inequality; \Cref{prop:stab} combines the three. In
\Cref{sec:proof-main}, \Cref{prop:stab} gives the Cauchy property of the
smooth approximations (existence), uniqueness, and stability.

Besides $\alpha_\gamma=1+(\gamma\wedge0)H$, defined in
\Cref{sec:framework}, the
sewing arguments involve the exponents $\beta_\gamma:=1-(1-\gamma)H$ and
$\zeta_\gamma:=1-(2-\gamma)H$, which arise from two-point and four-point
integration by parts, respectively. The subcritical condition
\eqref{eq:subcritical} implies
\begin{align}\label{eq:exponent-relations}
 \alpha_\gamma>\frac12,
 \qquad
 \beta_\gamma>\frac12,
 \qquad
 \zeta_\gamma>0,
 \qquad
 \alpha_\gamma+\beta_\gamma>1,
 \qquad
 \alpha_\gamma+\zeta_\gamma>1.
\end{align}
Indeed, the first inequality in \eqref{eq:subcritical} gives
$\beta_\gamma>1/2$, while the second gives $\zeta_\gamma>0$.  If
$\gamma\ge0$, then $\alpha_\gamma=1$, so all the remaining assertions are
immediate.  If $\gamma<0$, the first inequality is equivalent to
$\gamma H>H-1/2$, so $\alpha_\gamma>H+1/2>1/2$,
$\alpha_\gamma+\beta_\gamma=2-H+2\gamma H>1+H>1$, and
$\alpha_\gamma+\zeta_\gamma=2-2H+2\gamma H>1$.

In the sewing arguments below, the defect of a germ over $[u,v]$ is
evaluated at the midpoint $c:=(u+v)/2$, as in \Cref{lem:stochastic-sewing}.
The roles of the exponents are summarised as follows.
\begin{center}
\renewcommand{\arraystretch}{1.25}
\begin{tabular}{lll}
 \hline
 Exponent & Origin & Main use\\
 \hline
 $\alpha_\gamma=1+(\gamma\wedge0)H$ & one-point estimate &
 comparison with the flow, temporal regularity\\
 $\beta_\gamma=1-(1-\gamma)H$ & two-point estimate &
 regularity of the drift difference\\
 $\zeta_\gamma=1-(2-\gamma)H$ & four-point estimate &
 conditional defects in \Cref{sec:stab}\\
 $\chi_\eps=H-\eps+\frac1{k-1}$ & local scale \eqref{eq:local-scale} &
 bad-event bounds $C^\rho\rho^{\rho/2}h^{\chi_\eps\rho}$\\
 \hline
\end{tabular}
\end{center}
Throughout, constants of the form $C^\rho\rho^{\rho/2}$ have $C$ independent
of $\rho\ge1$; this uniformity is used in the final choice of $\rho$ in
\Cref{prop:stab}.

\subsection{A priori estimates for smooth drifts}

We begin with smooth singular drifts and derive estimates which are uniform at
the critical regularity. Fix $\Xi>0$ and let
$\mathscr B_{\gamma,\Xi}:=\{b\in C_b^\infty(\R^d;\R^d):
\|b\|_{\mathcal C^\gamma}\le\Xi\}$.
For $b\in\mathscr B_{\gamma,\Xi}$, we write $X^{s,x,\mathbf v}=X^{s,x,\mathbf v;b}$
for the solution of the smooth equation with drift $F+b$, see
\eqref{eq:X-Ito-map}, and $\phi^{s,x,\mathbf v}$ for the regularising flow.

We first record the form in which the flow estimates enter the
integration-by-parts lemmas. Let $I=[s,\bar t]$ have length at most one, let
$t\in(s,\bar t]$, fix $\ell\in\{1,\ldots,\lfloor-\gamma\rfloor+3\}$ and
$q>1$, and let $\mathbf w^1,\ldots,\mathbf w^L$ be deterministic inputs on
$I$ with sum $\mathbf w$. Then every affine combination of regularising flows
satisfies
\begin{align}\label{eq:Lambda-flow-scale}
\begin{split}
 \Lambda_{\ell,q}
 \bigl(\Phi_t^{s,\mathbf x,\lambda,\mathbf w}\bigr)
 \le C\bigl(R_{s,t}^{\mathbf w^1,\ldots,\mathbf w^L;I}\bigr)^C(t-s)^{-H},
 \qquad \lambda\in\mathcal A_M,
\end{split}
\end{align}
with $C$ depending on $\|\lambda\|_{\ell^1}$ and the parameters of
\Cref{prop:flow-estimates}. Indeed, by \eqref{eq:upperbound-D} and
\eqref{eq:inv-matrix}, the term of index $a\ge1$ in \eqref{eq:Lambda-def} is
bounded by $(R_{s,t}^{\mathbf w^1,\ldots,\mathbf w^L;I})^C
(t-s)^{-2aH}(t-s)^{(2a-1)H}
=(R_{s,t}^{\mathbf w^1,\ldots,\mathbf w^L;I})^C(t-s)^{-H}$,
and, since $t-s\le1$, the constant term satisfies the same bound. In the
applications, the inputs $\mathbf w^i$ are frozen $\mathcal F_s$-measurable
paths, and we keep in $R$ each component that carries its own local-scale
condition.

\begin{proposition}[Two-input comparison with the regularising flow]
\label{prop:regsol-inf}
Let $\eps\in(0,H)$, $m\in[2,\infty)$, and $\Xi>0$, and let $F$ satisfy
\Cref{ass:dissipative-drift}. There exists a constant
$C=C(d,H,\gamma,k,\boldsymbol\kappa,\eps,m,\Xi)<\infty$
such that the following holds. Let $b\in\mathscr B_{\gamma,\Xi}$, let
$I=[s,\bar t]$ be an interval of length at most one, let $x\in\R^d$, and let
$\mathbf v,\mathbf z\in C^{H-\eps}(I;\R^d)$. Then, for every $t\in(s,\bar t]$
with $t-s\le\tau_{\mathbf v;I}$,
\begin{align}\label{eq:two-input-flow-comparison}
\begin{split}
 \|X_t^{s,x,\mathbf v+\mathbf z}
 -\phi_t^{s,x,\mathbf v+\mathbf z}\|_{L_m}
 \le
 C\bigl(1+[\mathbf z]_{C^{H-\eps}(I)}\bigr)^C
 (t-s)^{\alpha_\gamma}.
\end{split}
\end{align}
Consequently, for every $\mathbf v\in C^{H-\eps}(I;\R^d)$ and
$t\in(s,\bar t]$,
\begin{align}\label{eq:local-flow-comparison}
\begin{split}
 \ind_{\{t-s\le\tau_{\mathbf v;I}\}}
 \|X_t^{s,x,\mathbf v}-\phi_t^{s,x,\mathbf v}\|_{L_m}
 \le C(t-s)^{\alpha_\gamma}
\end{split}
\end{align}
and
\begin{align}\label{eq:weighted-flow-comparison}
\begin{split}
 \|X_t^{s,x,\mathbf v}-\phi_t^{s,x,\mathbf v}\|_{L_m}
 \le C\bigl(1+[\mathbf v]_{C^{H-\eps}(I)}\bigr)^C
 (t-s)^{\alpha_\gamma}.
\end{split}
\end{align}
\end{proposition}

\begin{proof}
We first prove \eqref{eq:two-input-flow-comparison}; the two other
estimates are special cases, see the end of the proof. To keep the number of
indices manageable, we write the proof for $\bar t=s+1$, that is,
$I=[s,s+1]$; the general case is identical, with $s+1$ replaced by $\bar t$
throughout. The proof is done in three steps.

\medskip\noindent
\emph{Step 1: dissipativity.}
Write $X_r:=X_r^{s,x,\mathbf v+\mathbf z}$,
$\phi_r:=\phi_r^{s,x,\mathbf v+\mathbf z}$, and $\mathcal Z_r:=X_r-\phi_r$.
Since the two equations have the same additive signal,
$\frac{d}{dr}\mathcal Z_r=F(X_r)-F(\phi_r)+b(X_r)$. Using the one-sided Lipschitz estimate \eqref{eq:F-one-sided-global},
\begin{align*}
 \frac{d}{dr}|\mathcal Z_r|^2
 &=2\langle \mathcal Z_r,F(X_r)-F(\phi_r)\rangle
 +2\langle \mathcal Z_r,b(X_r)\rangle \le2\kappa_4|\mathcal Z_r|^2+2\langle \mathcal Z_r,b(X_r)\rangle.
\end{align*}
Hence
\begin{align}\label{eq:regZ-0}
\begin{split}
 |\mathcal Z_t|^2
 \le2\int_s^t e^{2\kappa_4(t-r)}
 \langle \mathcal Z_r,b(X_r)\rangle\,dr.
\end{split}
\end{align}

\medskip\noindent
\emph{Step 2: stochastic sewing.}
Assume that $t-s\le\tau_{\mathbf v;I}$. Recall that histories restarted
after $s$ are extended constantly to the left when their seminorm on $I$ is
taken; see \Cref{sec:framework}. For $u\in[s,t]$, let
$\mathbf w^u:=\mathbf v+\mathbf z+\bar B^{s,u}$ be the input frozen at time
$u$, and, for $s\le u<v\le t$, consider the germ
\begin{align*}
 A_{u,v}:=\int_u^v\EE^u\!\left[
 e^{2\kappa_4(t-r)}
 \left\langle
 \phi_r^{u,X_u,\mathbf w^u}-
 \phi_r^{u,\phi_u,\mathbf w^u},
 b\bigl(\phi_r^{u,X_u,\mathbf w^u}\bigr)
 \right\rangle
 \right]dr.
\end{align*}
For $x_1,x_2\in\R^d$ and deterministic paths
$\mathbf w^1,\mathbf w^2,\mathbf w^3$, write
$\mathbf w=\mathbf w^1+\mathbf w^2+\mathbf w^3$ and define
\begin{align*}
 \mathcal I_{u,r}(x_1,x_2,\mathbf w)
 :=\EE\left[
 \left\langle
 \phi_r^{u,x_1,\mathbf w}-\phi_r^{u,x_2,\mathbf w},
 b(\phi_r^{u,x_1,\mathbf w})
 \right\rangle\right].
\end{align*}
\Cref{cor:conditioning-restart}, applied to the flows restarted at $u$,
gives
\begin{align*}
 A_{u,v}
 =\EE^u\int_u^v e^{2\kappa_4(t-r)}
 \mathcal I_{u,r}(X_u,\phi_u,\mathbf w^u)\,dr.
\end{align*}
We apply \Cref{lem:malliavin-1} with $G=\phi_r^{u,x_1,\mathbf w}-\phi_r^{u,x_2,\mathbf w}$
and $Z_1=\phi_r^{u,x_1,\mathbf w}$, and \Cref{prop:flow-estimates} with the
three inputs kept separate. By
\eqref{eq:upperbound-D-x} and \eqref{eq:Lambda-flow-scale},
\begin{align*}
 \|\phi_r^{u,x_1,\mathbf w}-\phi_r^{u,x_2,\mathbf w}
 \|_{\mathbb D^{\lfloor-\gamma\rfloor+3,\bar q}}
 &\lesssim
 \bigl(R_{u,r}^{\mathbf w^1,\mathbf w^2,\mathbf w^3;I}\bigr)^C|x_1-x_2|,\\
 \Lambda_{\lfloor-\gamma\rfloor+3,q}(\phi_r^{u,x_1,\mathbf w})^{-(\gamma\wedge0)}
 &\lesssim
 \bigl(R_{u,r}^{\mathbf w^1,\mathbf w^2,\mathbf w^3;I}\bigr)^C
 (r-u)^{(\gamma\wedge0)H}.
\end{align*}

We thus obtain
\begin{align}\label{eq:regularity-one-point-after-flow}
\begin{split}
 |\mathcal I_{u,r}(x_1,x_2,\mathbf w)|
 \lesssim
 \Xi\bigl(R_{u,r}^{\mathbf w^1,\mathbf w^2,\mathbf w^3;I}\bigr)^C
 |x_1-x_2|(r-u)^{(\gamma\wedge0)H}.
\end{split}
\end{align}
After substituting
$(\mathbf w^1,\mathbf w^2,\mathbf w^3)
=(\mathbf v,\mathbf z,\bar B^{s,u})$, let
\begin{align*}
 G_{u,r}^{(0)}
 :=\{r-u\le\tau_{\mathbf z;I}\}
 \cap\{r-u\le\tau_{\bar B^{s,u};I}\}.
\end{align*}
The input $\mathbf v$ is below its scale because
$r-u\le t-s\le\tau_{\mathbf v;I}$. Thus the preceding bound can be
written on the entire probability space as
\begin{align}\label{eq:regularity-germ-indicators}
 |\mathcal I_{u,r}(X_u,\phi_u,\mathbf w^u)|
 \lesssim |\mathcal Z_u|(r-u)^{\alpha_\gamma-1}
 \Big[\ind_{G_{u,r}^{(0)}}
 +\ind_{(G_{u,r}^{(0)})^c}
 (R_{u,r}^{\mathbf v,\mathbf z,\bar B^{s,u};I})^C\Big].
\end{align}
For the second term use
\begin{align*}
 \ind_{(G_{u,r}^{(0)})^c}
 \le\ind_{\{r-u>\tau_{\mathbf z;I}\}}
 +\ind_{\{r-u>\tau_{\bar B^{s,u};I}\}}.
\end{align*}
Apply \eqref{eq:bad-event-basic} to $\mathbf z$ and
\Cref{lem:bad-event-moments} to $\bar B^{s,u}$. The other input factors
have fixed polynomial powers, so H\"older's inequality gives
\begin{align*}
 \big\|\ind_{(G_{u,r}^{(0)})^c}
 (R_{u,r}^{\mathbf v,\mathbf z,\bar B^{s,u};I})^C\big\|_{L_{2m}}\le C^\rho\rho^{\rho/2}
 (1+[\mathbf z]_{C^{H-\eps}(I)})^{C+\rho}
 (r-u)^{\chi_\eps\rho},\qquad\rho\ge1.
\end{align*}
Retaining $|\mathcal Z_u|$ in both terms of
\eqref{eq:regularity-germ-indicators}, taking norms by Cauchy--Schwarz,
and integrating in $r$ yields
\begin{align}\label{eq:regularity-germ-bound-detailed}
\begin{split}
 \|A_{u,v}\|_{L_m}
 \lesssim\|\mathcal Z_u\|_{L_{2m}}\Big(& (v-u)^{\alpha_\gamma}
 +C^\rho\rho^{\rho/2}
 (1+[\mathbf z]_{C^{H-\eps}(I)})^{C+\rho}
 (v-u)^{\alpha_\gamma+\chi_\eps\rho}\Big).
\end{split}
\end{align}

We now estimate the conditional sewing defect at the midpoint
$c=(u+v)/2$. The tower property and the flow property, together with
\eqref{eq:history-cocycle}, give
$\EE^u\delta A_{u,c,v}=\Theta_1+\Theta_2$, where
\begin{align*}
 \Theta_1
 &:=\EE^u\int_c^v e^{2\kappa_4(t-r)}
 \left\langle
 \phi_r^{u,X_u,\mathbf w^u}-\phi_r^{u,\phi_u,\mathbf w^u},
 b(\phi_r^{u,X_u,\mathbf w^u})-
 b(\phi_r^{c,X_c,\mathbf w^c})
 \right\rangle dr,\\
 \Theta_2
 &:=\EE^u\int_c^v e^{2\kappa_4(t-r)}
 \left\langle
 \phi_r^{u,X_u,\mathbf w^u}-\phi_r^{c,X_c,\mathbf w^c},
 b(\phi_r^{c,X_c,\mathbf w^c})
 \right\rangle dr.
\end{align*}
Introduce the $\mathcal F_c$-measurable starting points
$\xi_1:=\phi_c^{u,X_u,\mathbf w^u}$,
$\xi_2:=\phi_c^{u,\phi_u,\mathbf w^u}=\phi_c$, and $\xi_3:=X_c$.
For deterministic $x_1,x_2,x_3$ and
$\mathbf w=\mathbf w^1+\mathbf w^2+\mathbf w^3$, set
\begin{align*}
 \mathcal J_{1,r}(x_1,x_2,x_3,\mathbf w)
 &:=\EE\left[
 \left\langle
 \phi_r^{c,x_1,\mathbf w}-\phi_r^{c,x_2,\mathbf w},
 b(\phi_r^{c,x_1,\mathbf w})-b(\phi_r^{c,x_3,\mathbf w})
 \right\rangle\right],\\
 \mathcal J_{2,r}(x_1,x_3,\mathbf w)
 &:=\EE\left[
 \left\langle
 \phi_r^{c,x_1,\mathbf w}-\phi_r^{c,x_3,\mathbf w},
 b(\phi_r^{c,x_3,\mathbf w})
 \right\rangle\right].
\end{align*}
Conditioning at $c$ gives
\begin{align*}
 \Theta_1&=\EE^u\int_c^v e^{2\kappa_4(t-r)}
 \mathcal J_{1,r}(\xi_1,\xi_2,\xi_3,\mathbf w^c)\,dr
 \quad\text{and}\quad
 \Theta_2&=\EE^u\int_c^v e^{2\kappa_4(t-r)}
 \mathcal J_{2,r}(\xi_1,\xi_3,\mathbf w^c)\,dr.
\end{align*}
To close the estimate, we use a bootstrap quantity which involves two
structural exponents, chosen in the following order. First, $\rho_0\ge1$
is fixed so large that $\alpha_\gamma+\chi_\eps\rho_0>1$, which makes all
bad-event terms below admissible in the sewing lemma. Second, let $C_*$ be the
largest exponent of $1+[\mathbf z]_{C^{H-\eps}(I)}$ produced by the flow
estimates in the bounds below; $C_*$ depends only on the parameters of
\Cref{prop:flow-estimates} and not on the bootstrap quantity. We fix an
integer $L_*\ge C_*+\rho_0$. With these choices, for $p\ge2$ and an interval
$\mathbb J\subset[s,t]$, define
\begin{align*}
 \mathcal N_p(\mathbb J)
 :=\sup_{\substack{(a,r)\in\Delta_{\mathbb J},\ y\in\R^d\\
 \mathbf v',\mathbf z'\in C^{H-\eps}(I;\R^d)}}
 \frac{
 \ind_{\{r-a\le\tau_{\mathbf v';I}\}}
 \|X_r^{a,y,\mathbf v'+\mathbf z'}
 -\phi_r^{a,y,\mathbf v'+\mathbf z'}\|_{L_p}}
 {\bigl(1+[\mathbf z']_{C^{H-\eps}(I)}\bigr)^{L_*}
 (r-a)^{\alpha_\gamma}}.
\end{align*}
The two equations use the restrictions of $\mathbf v'$ and $\mathbf z'$
after the restarting time $a$. Thus $\mathcal N_p(\mathbb J)$ is a
bootstrap seminorm over all equations restarted inside $\mathbb J$, with the
local-scale condition imposed only on the first input $\mathbf v'$, which
we call the good input; the second input $\mathbf z'$ is the residual
input.
This quantity is finite before it is bootstrapped. Indeed, $b$ is smooth
and bounded in the present proposition, and the one-sided estimate used in
\eqref{eq:regZ-0}, applied to
$X^{a,y,\mathbf v'+\mathbf z'}-\phi^{a,y,\mathbf v'+\mathbf z'}$, gives
\begin{align*}
 \sup_{r'\in[a,r]}
 |X_{r'}^{a,y,\mathbf v'+\mathbf z'}
 -\phi_{r'}^{a,y,\mathbf v'+\mathbf z'}|
 \le
 e^{\kappa_4(r-a)}\|b\|_\infty(r-a).
\end{align*}
Since $\alpha_\gamma\le1$ and $r-a\le1$, it follows that
$\mathcal N_p(\mathbb J)<\infty$ for every finite $p$. We run the bootstrap at the
moment order $2m$; the resulting estimate in $L_{2m}$ implies the stated
$L_m$ estimate.

Because $c-u\le v-u\le t-s\le\tau_{\mathbf v;I}$, the first input
$\mathbf v$ satisfies the local-scale condition in the definition of
$\mathcal N_{2m}(\mathbb J)$.  The second input
$\mathbf z+\bar B^{s,u}$ is kept as the residual input; no local-scale
condition is imposed on its sum.
We first estimate $\Theta_1$. Apply \Cref{lem:malliavin-2} with
$G=\phi_r^{c,x_1,\mathbf w}-\phi_r^{c,x_2,\mathbf w}$,
$Z_1=\phi_r^{c,x_1,\mathbf w}$, and $Z_2=\phi_r^{c,x_3,\mathbf w}$.
For $\theta\in[0,1]$, the interpolation
$Z_\theta=\theta Z_1+(1-\theta)Z_2$ is an affine combination of two flows.
Therefore, \eqref{eq:upperbound-D-x} and
\eqref{eq:Lambda-flow-scale} give
\begin{align*}
 \|G\|_{\mathbb D^{\lfloor-\gamma\rfloor+3,\bar q}}
 &\lesssim
 \bigl(R_{c,r}^{\mathbf w^1,\mathbf w^2,\mathbf w^3;I}\bigr)^C|x_1-x_2|,\\
 \|Z_1-Z_2\|_{\mathbb D^{\lfloor-\gamma\rfloor+3,\bar q}}
 &\lesssim
 \bigl(R_{c,r}^{\mathbf w^1,\mathbf w^2,\mathbf w^3;I}\bigr)^C|x_1-x_3|,\\
 \sup_{\theta\in[0,1]}
 \Lambda_{\lfloor-\gamma\rfloor+3,q}(Z_\theta)^{1-\gamma}
 &\lesssim
 \bigl(R_{c,r}^{\mathbf w^1,\mathbf w^2,\mathbf w^3;I}\bigr)^C
 (r-c)^{-(1-\gamma)H}.
\end{align*}
Hence
\begin{align}\label{eq:J1-IBP-bound}
\begin{split}
 |\mathcal J_{1,r}(x_1,x_2,x_3,\mathbf w)|
 \lesssim
 \Xi\bigl(R_{c,r}^{\mathbf w^1,\mathbf w^2,\mathbf w^3;I}\bigr)^C
 |x_1-x_2|\,|x_1-x_3|(r-c)^{-(1-\gamma)H}.
\end{split}
\end{align}
We combine the good and complementary regions before taking norms.
Let $G(r):=G_{c,r}^{(0)}$ be the event defined before
\eqref{eq:regularity-germ-indicators}, with $u$ replaced by $c$; it is
$\mathcal F_c$-measurable, and $\ind_{G(r)}+\ind_{G(r)^c}=1$. In all the flow factors after substitution, the inputs are
kept separate as $(\mathbf v,\mathbf z,\bar B^{s,c})$.

On $G(r)^c$, \Cref{lem:malliavin-1}, applied to each
of the two evaluations of $b$ in $\mathcal J_{1,r}$, gives
\begin{align}\label{eq:J1-lower-order-bound}
 |\mathcal J_{1,r}(x_1,x_2,x_3,\mathbf w)|
 \lesssim\Xi(R_{c,r}^{\mathbf w^1,\mathbf w^2,\mathbf w^3;I})^C
 |x_1-x_2|(r-c)^{\alpha_\gamma-1}.
\end{align}
The same one-point estimate for $\mathcal J_{2,r}$ gives the factor
$|x_1-x_3|(r-c)^{\alpha_\gamma-1}$. By \Cref{lem:reg-phi}, and since
$\phi_c^{u,\phi_u,\mathbf w^u}=\phi_c$ by the flow property,
\begin{align*}
 |\xi_1-\xi_2|\le e^{\kappa_4}|\mathcal Z_u|,
 \qquad |\xi_1-\xi_3|\le|\xi_1-\xi_2|+|\phi_c-X_c|
 \le e^{\kappa_4}|\mathcal Z_u|+|\mathcal Z_c|.
\end{align*}
On $G(r)$ the separate-input control is one. Combining the two-point
and one-point bounds therefore yields the following two inequalities on
the entire probability space:
\begin{align}\label{eq:comparison-defect-indicators}
\begin{split}
 |\mathcal J_{1,r}(\xi_1,\xi_2,\xi_3,\mathbf w^c)|
 &\lesssim
 \ind_{G(r)}|\mathcal Z_u|\,|\xi_1-\xi_3|
 (r-c)^{\beta_\gamma-1}
 +\ind_{G(r)^c}(R_{c,r}^{\mathbf v,\mathbf z,\bar B^{s,c};I})^C
 |\mathcal Z_u|(r-c)^{\alpha_\gamma-1},\\
 |\mathcal J_{2,r}(\xi_1,\xi_3,\mathbf w^c)|
 &\lesssim
 \ind_{G(r)}|\xi_1-\xi_3|
 (r-c)^{\alpha_\gamma-1}
 +\ind_{G(r)^c}(R_{c,r}^{\mathbf v,\mathbf z,\bar B^{s,c};I})^C
 (|\mathcal Z_u|+|\mathcal Z_c|)(r-c)^{\alpha_\gamma-1}.
\end{split}
\end{align}
The constants include the fixed bound $\Xi$.

For the first terms, condition at $u$ and apply the bootstrap seminorm to
the good input $\mathbf v$ and residual input
$\mathbf z+\bar B^{s,u}$. Since $c-u\le\tau_{\mathbf v;I}$,
\begin{align}\label{eq:two-input-bootstrap-remainder}
 \big\|
 (\EE^u|\xi_1-\xi_3|^{2m})^{1/(2m)}\big\|_{L_{2m}}
 \le C\mathcal N_{2m}(\mathbb J)
 (1+[\mathbf z]_{C^{H-\eps}(I)})^{L_*}
 (c-u)^{\alpha_\gamma}.
\end{align}
Indeed, conditionally on $\mathcal F_u$, $\xi_3-\xi_1=X_c-\phi_c^{u,X_u,\mathbf w^u}$
is the comparison remainder of the equation restarted at $u$ from $X_u$,
with good input $\mathbf v$ and residual input $\mathbf z+\bar B^{s,u}$; the
history then contributes only a fixed Gaussian moment of order determined by
$m$ and $L_*$, by \Cref{lem:history-holder-moments}. Conditional Jensen and ordinary Cauchy--Schwarz
then control the first terms of
\eqref{eq:comparison-defect-indicators}.

The residual inputs enter only through the complement of $G(r)$, for
which
\begin{align*}
 \ind_{G(r)^c}
 \le(r-c)^{\chi_\eps\rho}
 \Big[(1+[\mathbf z]_{C^{H-\eps}(I)})^\rho
 +(1+[\bar B^{s,c}]_{C^{H-\eps}(I)})^\rho\Big]
\end{align*}
by \eqref{eq:bad-event-basic}.
We use \Cref{lem:bad-event-moments}, absorbing fixed extra powers of the
histories into $C^\rho\rho^{\rho/2}$.
Cauchy--Schwarz pairs these input factors with $\mathcal Z_u$ or $\mathcal Z_c$ in
$L_{2m}$. Integration of the single bound
\eqref{eq:comparison-defect-indicators} now gives
\begin{align}\label{eq:J1-final}
\begin{split}
 \|\Theta_1\|_{L_m}
 &\lesssim
 \sup_{a\in[s,t]}\|\mathcal Z_a\|_{L_{2m}}\,
 \mathcal N_{2m}(\mathbb J)
 (1+[\mathbf z]_{C^{H-\eps}(I)})^{L_*}
 (c-u)^{\alpha_\gamma}(v-c)^{\beta_\gamma}\\
 &\quad+C^\rho\rho^{\rho/2}
 (1+[\mathbf z]_{C^{H-\eps}(I)})^{C+\rho}
 \sup_{a\in[s,t]}\|\mathcal Z_a\|_{L_{2m}}
 (v-u)^{\alpha_\gamma+\chi_\eps\rho},
\end{split}
\end{align}
and
\begin{align}\label{eq:J2-final}
\begin{split}
 \|\Theta_2\|_{L_m}
 &\lesssim
 \mathcal N_{2m}(\mathbb J)
 (1+[\mathbf z]_{C^{H-\eps}(I)})^{L_*}
 (c-u)^{\alpha_\gamma}(v-c)^{\alpha_\gamma}\\
 &\quad+C^\rho\rho^{\rho/2}
 (1+[\mathbf z]_{C^{H-\eps}(I)})^{C+\rho}
 \sup_{a\in[s,t]}\|\mathcal Z_a\|_{L_{2m}}
 (v-u)^{\alpha_\gamma+\chi_\eps\rho}.
\end{split}
\end{align}
The complementary indicator supplies the power
$(v-c)^{\alpha_\gamma+\chi_\eps\rho}\le(v-u)^{\alpha_\gamma+\chi_\eps\rho}$.

Combining \eqref{eq:J1-final} and \eqref{eq:J2-final},
\begin{align*}
 \|\EE^u\delta A_{u,c,v}\|_{L_m}
 &\lesssim
 \sup_{a\in[s,t]}\|\mathcal Z_a\|_{L_{2m}}\,
 \mathcal N_{2m}(\mathbb J)
 \bigl(1+[\mathbf z]_{C^{H-\eps}(I)}\bigr)^{L_*}
 (c-u)^{\alpha_\gamma}(v-c)^{\beta_\gamma}\\
 &\quad+
 \mathcal N_{2m}(\mathbb J)
 \bigl(1+[\mathbf z]_{C^{H-\eps}(I)}\bigr)^{L_*}
 (c-u)^{\alpha_\gamma}(v-c)^{\alpha_\gamma}\\
 &\quad+C^\rho\rho^{\rho/2}
 \bigl(1+[\mathbf z]_{C^{H-\eps}(I)}\bigr)^{C+\rho}
 \sup_{a\in[s,t]}\|\mathcal Z_a\|_{L_{2m}}
 (v-u)^{\alpha_\gamma+\chi_\eps\rho}.
\end{align*}
The non-conditional defect is controlled by
\eqref{eq:regularity-germ-bound-detailed} and the triangle inequality. By
\eqref{eq:exponent-relations} and the choice of $\rho$, the sewing estimates
have the required exponents. Let
$\mathcal Y_r:=\int_s^r e^{2\kappa_4(t-r')}\langle \mathcal Z_{r'},b(X_{r'})\rangle\,dr'$.
Since $b$ is smooth, $\mathcal Y$ is the process associated with the germ $A$
by \Cref{lem:stochastic-sewing}; see the discussion following that lemma.
Hence \Cref{lem:stochastic-sewing} applies to $\mathcal Y$. Since $\mathcal Z_s=0$, the two
flows in the germ coincide for $u=s$, so that $A_{s,t}=0$; the
non-conditional defect still contributes the first term below. We obtain
\begin{align}\label{eq:sewn-interaction-bound}
\begin{split}
 \|\mathcal Y_t\|_{L_m}
 &\lesssim
 \sup_{a\in[s,t]}\|\mathcal Z_a\|_{L_{2m}}(t-s)^{\alpha_\gamma}\\
 &\quad+
 \sup_{a\in[s,t]}\|\mathcal Z_a\|_{L_{2m}}\,
 \mathcal N_{2m}(\mathbb J)
 \bigl(1+[\mathbf z]_{C^{H-\eps}(I)}\bigr)^{L_*}
 (t-s)^{\alpha_\gamma+\beta_\gamma}\\
 &\quad+\mathcal N_{2m}(\mathbb J)
 \bigl(1+[\mathbf z]_{C^{H-\eps}(I)}\bigr)^{L_*}
 (t-s)^{2\alpha_\gamma}\\
 &\quad+C^\rho\rho^{\rho/2}
 \bigl(1+[\mathbf z]_{C^{H-\eps}(I)}\bigr)^{C+\rho}
 \sup_{a\in[s,t]}\|\mathcal Z_a\|_{L_{2m}}
 (t-s)^{\alpha_\gamma+\chi_\eps\rho}.
\end{split}
\end{align}

\medskip\noindent
\emph{Step 3: conclusion.}
Let $\rho_0$ and $L_*$ be as fixed above.  Let
$\mathbb J=[S,S']\subset[s,t]$ and write $h:=S'-S$. Repeat Steps 1 and 2 with
arbitrary $(a,r)\in\Delta_{\mathbb J}$, initial condition $y$, good input
$\mathbf v'$, and residual input $\mathbf z'$ occurring in the definition of
$\mathcal N_{2m}(\mathbb J)$, and write
$W:=1+[\mathbf z']_{C^{H-\eps}(I)}$.
On the local-scale regime for $\mathbf v'$, the bootstrap bound gives
\begin{align*}
 \sup_{r'\in[a,r]}
 \|X_{r'}^{a,y,\mathbf v'+\mathbf z'}
 -\phi_{r'}^{a,y,\mathbf v'+\mathbf z'}\|_{L_{2m}}
 \le \mathcal N_{2m}(\mathbb J)W^{L_*}
 (r-a)^{\alpha_\gamma}.
\end{align*}
Apply \eqref{eq:sewn-interaction-bound} with $\rho=\rho_0$ and
$(\mathbf v,\mathbf z)=(\mathbf v',\mathbf z')$.  Combining it with
\eqref{eq:regZ-0}, dividing by $W^{2L_*}(r-a)^{2\alpha_\gamma}$, and taking the supremum over
the parameters in the definition of $\mathcal N_{2m}(\mathbb J)$ gives
\begin{align*}
 \mathcal N_{2m}(\mathbb J)^2
 &\lesssim
 h^{\beta_\gamma}\mathcal N_{2m}(\mathbb J)^2
 +\mathcal N_{2m}(\mathbb J)
 +C^{\rho_0}\rho_0^{\rho_0/2}
 h^{\chi_\eps\rho_0}
 \mathcal N_{2m}(\mathbb J).
\end{align*}
Here the choice $L_*\ge C_*+\rho_0$ ensures that every residual-input power
is cancelled by the normalising factor.  Choose $h_0>0$ small enough, depending only on the parameters of the
proposition, so that the first
term is absorbed whenever $h\le h_0$. Young's inequality applied to the two
linear terms then yields
\begin{align}\label{eq:two-input-short-window}
 \mathcal N_{2m}(\mathbb J)\le C
 \qquad\text{whenever }|\mathbb J|\le h_0.
\end{align}

We next propagate the two-input short-window estimate. For
$s\le a<r\le t$, set
\begin{align*}
 \mathfrak R_{a,r}
 :=X_r^{s,x,\mathbf v+\mathbf z}
 -\phi_r^{a,X_a^{s,x,\mathbf v+\mathbf z},
 \mathbf v+\mathbf z+\bar B^{s,a}}.
\end{align*}
In particular, $\mathfrak R_{s,r}=X_r^{s,x,\mathbf v+\mathbf z}
-\phi_r^{s,x,\mathbf v+\mathbf z}$. If $s\le S<S''<S'\le t$, the restarting identities for the smooth equation and the
regularising flow, together with \eqref{eq:history-cocycle}, give
\begin{align*}
 \mathfrak R_{S,S'}
 &=\mathfrak R_{S'',S'} +
 \phi_{S'}^{S'',X_{S''}^{s,x,\mathbf v+\mathbf z},
 \mathbf v+\mathbf z+\bar B^{s,S''}} -
 \phi_{S'}^{S'',\phi_{S''}^{S,X_S^{s,x,\mathbf v+\mathbf z},
 \mathbf v+\mathbf z+\bar B^{s,S}},
 \mathbf v+\mathbf z+\bar B^{s,S''}}.
\end{align*}
The Lipschitz dependence of the regularising flow on its initial condition
therefore yields the subadditivity estimate
\begin{align}\label{eq:comparison-subadditivity}
\begin{split}
 \|\mathfrak R_{S,S'}\|_{L_{2m}}
 \le\|\mathfrak R_{S'',S'}\|_{L_{2m}}
 +e^{\kappa_4(S'-S'')}\|\mathfrak R_{S,S''}\|_{L_{2m}}.
\end{split}
\end{align}
The first term can be estimated conditionally at $S''$, since
\begin{align*}
 \|\mathfrak R_{S'',S'}\|_{L_{2m}}
 =\bigl\|(\EE^{S''}|\mathfrak R_{S'',S'}|^{2m})^{1/(2m)}\bigr\|_{L_{2m}};
\end{align*}
thus
\eqref{eq:two-input-short-window} may be applied inside the conditional
expectation, with the state and the accumulated history frozen at time
$S''$.

Choose a uniform deterministic partition $s=S_0<\cdots<S_N=t$ with
$N:=\lceil (t-s)/h_0\rceil\vee1$ and $h_j:=S_{j+1}-S_j\le h_0$. The number
of intervals is bounded by a constant because $t-s\le1$, and
$h_j\le t-s\le\tau_{\mathbf v;I}$. Set $\mathbf z^j:=\mathbf z+\bar B^{s,S_j}$.
Conditional on $\mathcal F_{S_j}$, the state
$X_{S_j}^{s,x,\mathbf v+\mathbf z}$ and the residual input
$\mathbf z^j$ are fixed, whereas $B^{S_j}$ is independent of
$\mathcal F_{S_j}$. Applying \eqref{eq:two-input-short-window} conditionally,
with good input $\mathbf v$ and residual input $\mathbf z^j$, gives
\begin{align*}
 \bigl(\EE^{S_j}|\mathfrak R_{S_j,S_{j+1}}|^{2m}\bigr)^{1/(2m)} &\le
 C\bigl(1+[\mathbf z+\bar B^{s,S_j}]_{C^{H-\eps}(I)}
 \bigr)^{L_*}h_j^{\alpha_\gamma}.
\end{align*}
Taking the outer $L_{2m}$-norm, using the triangle inequality for the
H\"older seminorm and \Cref{lem:history-holder-moments}, yields
\begin{align*}
 \|\mathfrak R_{S_j,S_{j+1}}\|_{L_{2m}}
 \le
 C\bigl(1+[\mathbf z]_{C^{H-\eps}(I)}\bigr)^{L_*}
 h_j^{\alpha_\gamma}.
\end{align*}
Thus the fixed Gaussian moment of the history is absorbed into the
constant.  Iterating
\eqref{eq:comparison-subadditivity} over the partition now yields
\begin{align*}
 \|X_t^{s,x,\mathbf v+\mathbf z}
 -\phi_t^{s,x,\mathbf v+\mathbf z}\|_{L_{2m}} &\le
 C\bigl(1+[\mathbf z]_{C^{H-\eps}(I)}\bigr)^{L_*}
 \sum_{j=0}^{N-1}h_j^{\alpha_\gamma}.
\end{align*}
Since $N$ is bounded and $\sum_jh_j=t-s$, concavity gives
\begin{align*}
 \sum_{j=0}^{N-1}h_j^{\alpha_\gamma}
 \le N^{1-\alpha_\gamma}(t-s)^{\alpha_\gamma}
 \le C(t-s)^{\alpha_\gamma}.
\end{align*}
This proves \eqref{eq:two-input-flow-comparison} in $L_{2m}$, and hence in
$L_m$.

Finally, taking $\mathbf z=0$ in \eqref{eq:two-input-flow-comparison} proves
\eqref{eq:local-flow-comparison}. Applying
\eqref{eq:two-input-flow-comparison} with good input $0$ and residual input
$\mathbf v$, that is, with $(\mathbf v,\mathbf z)$ replaced by
$(0,\mathbf v)$, proves \eqref{eq:weighted-flow-comparison}, because
$\tau_{0;I}=1\ge t-s$.
\end{proof}

The preceding proposition is formulated for deterministic initial data and
input. The next statement transfers it to the conditional restarting
configuration used throughout the paper.

\begin{proposition}[Restarted comparison estimate]\label{prop:regsol}
Let $t_0<T$, $\eps\in(0,H)$, $m\in[2,\infty)$, and $\Xi>0$, and let $F$ satisfy
\Cref{ass:dissipative-drift}. There exist $q\ge m$ and constants
$C_1,C_2<\infty$, depending only on
$d,H,\gamma,k,\boldsymbol\kappa,\eps,m$, and $\Xi$, such that the following
holds. Let $b\in\mathscr B_{\gamma,\Xi}$, let $X_{t_0}$
be $\mathcal F_{t_0}$-measurable, and let $V$ be an
$\mathcal F_{t_0}$-measurable continuous random path. Let $X$ be the
corresponding smooth solution on $[t_0,T]$. Then, for every
$t_0\le s<t\le T$ with $t-s\le1$, setting
$I_s=[s,(s+1)\wedge T]$, one has
\begin{align}\label{eq:conditional-restarted-comparison}
\begin{split}
 &\left\|
 \ind_{\{t-s\le\tau_{V+\bar B^{t_0,s};I_s}\}}
 \bigl(\EE^s|X_t-
 \phi_t^{s,X_s,V+\bar B^{t_0,s}}|^m\bigr)^{1/m}
 \right\|_{L_\infty}
 \le C_1(t-s)^{\alpha_\gamma}.
\end{split}
\end{align}
If, in addition,
$\|[V]_{C^{H-\eps}(I_s)}\|_{L_q}<\infty$, then
\begin{align}\label{eq:unconditional-restarted-comparison}
\begin{split}
 \|X_t-
 \phi_t^{s,X_s,V+\bar B^{t_0,s}}\|_{L_m}
 \le C_1\left(1+
 \|[V]_{C^{H-\eps}(I_s)}\|_{L_q}\right)^{C_2}
 (t-s)^{\alpha_\gamma}.
\end{split}
\end{align}
\end{proposition}

\begin{proof}
We apply \Cref{prop:regsol-inf} on the interval $I_s=[s,\bar t]$ with
$\bar t=(s+1)\wedge T$; note that $t\in(s,\bar t]$ and that the frozen
input $V+\bar B^{t_0,s}$ is defined on $I_s$. For fixed $t$, the
conditioning identity of \Cref{cor:conditioning-restart} gives
\begin{align*}
 \EE^s\left[
 \left|X_t-
 \phi_t^{s,X_s,V+\bar B^{t_0,s}}\right|^m
 \right] &=
 \left.
 \EE\left[
 \left|X_t^{s,x,\mathbf v}-\phi_t^{s,x,\mathbf v}\right|^m
 \right]
 \right|_{x=X_s,\,\mathbf v=V+\bar B^{t_0,s}}.
\end{align*}
Applying \eqref{eq:local-flow-comparison} with $I=I_s$ pointwise in the
$\mathcal F_s$-measurable parameters proves
\eqref{eq:conditional-restarted-comparison}. Applying instead
\eqref{eq:weighted-flow-comparison}, taking the outer $L_m$-norm, and using
the triangle inequality
$[V+\bar B^{t_0,s}]_{C^{H-\eps}(I_s)}\le[V]_{C^{H-\eps}(I_s)}
+[\bar B^{t_0,s}]_{C^{H-\eps}(I_s)}$
together with \Cref{lem:history-holder-moments} proves
\eqref{eq:unconditional-restarted-comparison}.
\end{proof}

The restarted comparison controls the contribution of the singular drift. In
combination with the pathwise bounds for the regularising flow, it also yields
uniform moments for the full smooth solution.

\begin{lemma}[Uniform moments for smooth solutions]\label{lem:smooth-solution-moments}
Let $t_0<T$, $\eps\in(0,H)$, $p\in[2,\infty)$, and $\Xi>0$, and let $F$ satisfy
\Cref{ass:dissipative-drift}. There exist $q\ge p$ and constants
$C_1,C_2<\infty$, depending only on
$d,H,\gamma,k,\boldsymbol\kappa,\eps,p,\Xi$, and $T-t_0$, such that the
following holds. For every $b\in\mathscr B_{\gamma,\Xi}$,
every $\mathcal F_{t_0}$-measurable $X_{t_0}\in L_q(\Omega;\R^d)$, and every
$\mathcal F_{t_0}$-measurable continuous random path $V$ with
$\|[V]_{C^{H-\eps}([t_0,T])}\|_{L_q}<\infty$,
the corresponding smooth solution $X$ satisfies
\begin{align}\label{eq:smooth-solution-moments}
\begin{split}
 \sup_{r\in[t_0,T]}\|X_r\|_{L_p}
 \le C_1\left(
 1+\|X_{t_0}\|_{L_q}
 +\|[V]_{C^{H-\eps}([t_0,T])}\|_{L_q}
 \right)^{C_2}.
\end{split}
\end{align}
\end{lemma}

\begin{proof}
Choose a deterministic partition
$t_0=s_0<s_1<\cdots<s_N=T$ with $s_{j+1}-s_j\le1$. For
$r\in(s_j,s_{j+1}]$, restart at $s_j$ and write $X_r$ as the sum of
$\phi_r^{s_j,X_{s_j},V+\bar B^{t_0,s_j}}$ and the comparison remainder
$X_r-\phi_r^{s_j,X_{s_j},V+\bar B^{t_0,s_j}}$. The remainder is bounded by
\eqref{eq:unconditional-restarted-comparison}. For the flow, the additive
signal is $V_{s_j,r}+\bar B_r^{t_0,s_j}+B_r^{s_j}=V_{s_j,r}+B^{t_0}_{s_j,r}$
by \eqref{eq:restart-RL}, so the history is part of the increment of
$B^{t_0}$. Hence the pathwise energy estimate used in the proof of
\Cref{lem:coming-down-infinity}, without discarding the initial value, gives
\begin{align*}
 \sup_{r\in[s_j,s_{j+1}]}
 |\phi_r^{s_j,X_{s_j},V+\bar B^{t_0,s_j}}|
 \le C\bigl(1+|X_{s_j}|
 +[V]_{C^{H-\eps}([t_0,T])}
 +[B^{t_0}]_{C^{H-\eps}([t_0,T])}\bigr)^C.
\end{align*}
Taking moments, using Gaussian H\"older moments for $B^{t_0}$, and applying
H\"older's inequality gives the required bound on the first interval.
Induction over the finite partition gives the same bound on every interval,
with constants depending on $N$, hence only on $T-t_0$. This proves
\eqref{eq:smooth-solution-moments}.
\end{proof}

We can now combine the moment estimate with the regularisation property of the
singular drift to control ordinary temporal increments.

\begin{proposition}[Temporal regularity]\label{prop:regsol-usual}
Let $t_0<T$, $\eps\in(0,H)$, $m\in[2,\infty)$, and $\Xi>0$, and let $F$ satisfy
\Cref{ass:dissipative-drift}. There exist $q\ge m$ and constants
$C_1,C_2<\infty$, depending only on
$d,H,\gamma,k,\boldsymbol\kappa,\eps,m,\Xi$, and $T-t_0$, such that the
following holds. For every $b\in\mathscr B_{\gamma,\Xi}$,
every $\mathcal F_{t_0}$-measurable $X_{t_0}\in L_q(\Omega;\R^d)$, and every
$\mathcal F_{t_0}$-measurable continuous random path $V$ with
$\|[V]_{C^{H-\eps}([t_0,T])}\|_{L_q}<\infty$,
the corresponding smooth solution $X$ satisfies
\begin{align}\label{eq:smooth-temporal-regularity}
\begin{split}
 [X-V-B^{t_0}]_{C^{\alpha_\gamma}([t_0,T];L_m)}
 \le C_1\left(
 1+\|X_{t_0}\|_{L_q}
 +\|[V]_{C^{H-\eps}([t_0,T])}\|_{L_q}
 \right)^{C_2}.
\end{split}
\end{align}
\end{proposition}

\begin{proof}
Fix an interval $\mathbb J=[S,S']\subset[t_0,T]$ with $S'-S\le1$, and set
$\Upsilon:=X-V-B^{t_0}$ and
$\mathscr Y_m(\mathbb J):=[\Upsilon]_{C^{\alpha_\gamma}(\mathbb J;L_m)}$.
This seminorm is finite before it is absorbed. Indeed, $b$ is smooth and
bounded, while \Cref{lem:smooth-solution-moments} and the polynomial growth
of $F$ imply $\sup_{r\in[t_0,T]}\|F(X_r)+b(X_r)\|_{L_m}<\infty$.
Thus $\Upsilon$ is Lipschitz as an $L_m$-valued path, and
$\mathscr Y_m(\mathbb J)<\infty$ because $\alpha_\gamma \leq 1$.
For $S\le s<t\le S'$, the equation gives
\begin{align}\label{eq:temporal-regularity-equation}
\begin{split}
 \Upsilon_t-\Upsilon_s
 =\int_s^tF(X_r)\,dr+\int_s^tb(X_r)\,dr.
\end{split}
\end{align}
By \Cref{lem:smooth-solution-moments} and the polynomial growth of $F$,
there exists $q\ge m$ such that
\begin{align}\label{eq:F-term-temporal-regularity}
\begin{split}
 \left\|\int_s^tF(X_r)\,dr\right\|_{L_m}
 \le C\left(1+\|X_{t_0}\|_{L_q}
 +\|[V]_{C^{H-\eps}([t_0,T])}\|_{L_q}\right)^C(t-s).
\end{split}
\end{align}

We apply \cite[Proposition 3.3]{DHL25} to the smooth drift $b$, with
$\psi:=X-B^{t_0}=\Upsilon+V$, $q=m$, and H\"older exponent $\alpha_\gamma$ for
$\psi$. Its hypotheses are $\gamma>-1/(2H)$, which follows from
\eqref{eq:subcritical}, and $H(\gamma-1)+\alpha_\gamma>0$, which is
$\alpha_\gamma+\beta_\gamma>1$, see \eqref{eq:exponent-relations}. Since $V$
is $\mathcal F_{t_0}$-measurable, the centred conditional seminorm of $\psi$
introduced below coincides with that of $\Upsilon$. With
$\delta_\gamma:=\alpha_\gamma+\beta_\gamma-1>0$, this gives
\begin{align*}
 \left\|\int_s^tb(X_r)\,dr\right\|_{L_m}
 \le C(t-s)^{\alpha_\gamma}
 +C\left(1+[\Upsilon]_{\alpha_\gamma,m;\mathbb J}^{\circ}\right)(t-s)^{1+\delta_\gamma},
\end{align*}
where the centred conditional seminorm is
\begin{align*}
 [\Upsilon]_{\alpha_\gamma,m;\mathbb J}^{\circ}
 :=\sup_{(u,v)\in\Delta_{\mathbb J}}
 \frac{\|\Upsilon_v-\EE^u\Upsilon_v\|_{L_m}}{(v-u)^{\alpha_\gamma}}.
\end{align*}
Since $\Upsilon_u$ is $\mathcal F_u$-measurable,
$\Upsilon_v-\EE^u\Upsilon_v=\Upsilon_{u,v}-\EE^u\Upsilon_{u,v}$.
Conditional expectation is a contraction on $L_m$, and therefore
\begin{align*}
 \|\Upsilon_v-\EE^u\Upsilon_v\|_{L_m}
 &\le\|\Upsilon_v-\Upsilon_u\|_{L_m}
 +\|\EE^u(\Upsilon_v-\Upsilon_u)\|_{L_m} \le2\mathscr Y_m(\mathbb J)(v-u)^{\alpha_\gamma}.
\end{align*}
Thus $[\Upsilon]_{\alpha_\gamma,m;\mathbb J}^{\circ}\le2\mathscr Y_m(\mathbb J)$.
Consequently,
\begin{align}\label{eq:occupation-via-target-seminorm}
\begin{split}
 \left\|\int_s^tb(X_r)\,dr\right\|_{L_m}
 \le C(t-s)^{\alpha_\gamma}
 +C\bigl(1+\mathscr Y_m(\mathbb J)\bigr)(t-s)^{1+\delta_\gamma}.
\end{split}
\end{align}
Combining \eqref{eq:temporal-regularity-equation},
\eqref{eq:F-term-temporal-regularity}, and
\eqref{eq:occupation-via-target-seminorm}, dividing by
$(t-s)^{\alpha_\gamma}$, and taking the supremum over $\Delta_{\mathbb J}$, we get
\begin{align*}
 \mathscr Y_m(\mathbb J)
 \le C\left(1+\|X_{t_0}\|_{L_q}
 +\|[V]_{C^{H-\eps}([t_0,T])}\|_{L_q}\right)^C
 +C(S'-S)^{\beta_\gamma}\mathscr Y_m(\mathbb J).
\end{align*}
Choose $S'-S$ small enough to absorb the last term. Covering $[t_0,T]$ by
finitely many such intervals and splitting increments at the partition
points proves \eqref{eq:smooth-temporal-regularity}.
\end{proof}

\subsection{Stability estimates}\label{sec:stab}

We now compare two solutions with possibly different singular drifts. The
cutoff $p_k=8((k-2)\vee1)$ in the comparison class is the largest
conditional moment of the comparison remainder used below: the orders that
appear are $8$, $8((k-2)\vee0)$, and $4(k-1)$, arising in
\Cref{lem:stability-dissipativity} from the polynomial factors in
$\nabla^2F$ and $\nabla F$, and $4(k-1)\le8((k-2)\vee1)$; the other lemmas
of this subsection only use conditional moments of order at most two. The
cutoff does not depend on the outer moment parameter $m$.

\paragraph{Standing setting.}
Throughout this subsection, we fix $t_0<T$, $m\in[2,\infty)$, $\delta_0>0$,
$C_0<\infty$, $\Xi>0$, and $\eps\in(0,H)$; the input $V$ satisfies
\Cref{ass:input}, and $b,g\in\mathcal C^\gamma(\R^d;\R^d)$ satisfy
$\|b\|_{\mathcal C^\gamma}+\|g\|_{\mathcal C^\gamma}\le\Xi$. We let $X$ and
$Y$ be solutions driven by the same input and the same Riemann--Liouville
process, with singular drifts $b$ and $g$, respectively, and drift processes
$K^X$ and $K^Y$, and we write $\mathbb K:=K^X-K^Y$. We assume that
\begin{align}\label{eq:dissipativity-fixed-bounds}
\begin{split}
 &\sup_{r\in[t_0,T]}
 \bigl(\|X_r\|_{L_{2m+\delta_0}}+
       \|Y_r\|_{L_{2m+\delta_0}}\bigr)\le C_0,\\
 &\sup_{\substack{(s,t)\in\Delta_{[t_0,T]}\\t-s\le1}}
 \frac{\left\|
 \ind_{\{t-s\le\tau_{V+\bar B^{t_0,s};I_s}\}}
 \bigl(\EE^s|Z_t-
 \phi_t^{s,Z_s,V+\bar B^{t_0,s}}|^p\bigr)^{1/p}
 \right\|_{L_\infty}}
 {(t-s)^{\alpha_\gamma}}
 \le C_0
\end{split}
\end{align}
for $Z=X,Y$ and every $p\in[2,p_k]$; the second bound is the comparison
estimate \eqref{eq:classVinf} with a common constant. Unless stated
otherwise, constants depend only on
$d$, $H$, $\gamma$, $k$, $\boldsymbol\kappa$, $\eps$, $T-t_0$, $\Lambda_{V,\eps}$, $m$, $\delta_0$, $\Xi$,
and $C_0$. All regularising flows restarted at a time $s$ carry the frozen
input $V+\bar B^{t_0,s}$, which we omit from the notation when no confusion
can arise.

For $t_0\le s\le r\le T$, we introduce the double comparison remainder
\begin{align*}
 \mathsf D_{s,r}
 :=X_r-\phi_r^{s,X_s}-Y_r+\phi_r^{s,Y_s},
 \qquad
 \phi_r^{s,Z_s}:=\phi_r^{s,Z_s,V+\bar B^{t_0,s}},
\end{align*}
and, for an interval $\mathbb J\subset[t_0,T]$, the two seminorms
\begin{align}\label{eq:QK-def}
 \mathcal Q_m(\mathbb J)
 :=\sup_{(s,t)\in\Delta_{\mathbb J}}
 \frac{\|\EE^s|\mathsf D_{s,t}|^2\|_{L_m}}{t-s},
 \qquad
 \mathcal K_m(\mathbb J)
 :=\sup_{(s,t)\in\Delta_{\mathbb J}}
 \frac{\|\mathbb K_{s,t}\|_{L_{2m}}}{(t-s)^{\beta_\gamma}}.
\end{align}
The first measures how far the two solutions deviate, conditionally, from
the difference of their restarted regularising flows; the second is the
H\"older seminorm of the difference of the drift processes. The stability
estimate is obtained by combining three bounds: \Cref{lem:stability-dissipativity}
controls $\mathcal Q_m$ by $\mathcal K_m$ and the state difference,
\Cref{lem:stability-sewing} controls $\mathcal K_m$ by $\mathcal Q_m$, the
state difference, and $\|b-g\|_{\mathcal C^{\gamma-1}}$, and
\Cref{lem:stability-energy-sewing} provides the energy inequality for
$X-Y$. \Cref{prop:stab} combines the three bounds on a mesh of size $h$
through a discrete Gr\"onwall argument. The difficulty is that each local
estimate gains only a factor $h^{\beta_\gamma}$ with $\beta_\gamma<1$, so
that the iteration over the $h^{-1}$ intervals of the mesh produces a factor
$\exp(Ch^{-(1-\beta_\gamma)})$. The contributions of the events on which an
input exceeds its local scale are of size $C^\rho\rho^{\rho/2}h^{\chi_\eps\rho}$,
hence of order $\exp(-ah^{-2\chi_\eps})$ for $\rho$ of order
$h^{-2\chi_\eps}$. Choosing $h$ so that these contributions match the target
error, the Gr\"onwall factor is smaller than any negative power of that
error because $2\chi_\eps>1-\beta_\gamma$; this is where the range
\eqref{eq:stability-epsilon-range}, and hence the compatibility condition
\eqref{eq:H-k-iteration-condition}, enters.

Of condition (ii) in \Cref{def:solution}, only the convergence along the
heat approximations $P_{1/N}b$ and $P_{1/N}g$ is used in this subsection.
In particular, its results apply when $b,g\in C_b^\infty(\R^d;\R^d)$ and $X$
and $Y$ are classical solutions, since $P_{1/N}b\to b$ and $P_{1/N}g\to g$
uniformly in that case; this is how they are used in the existence proof.

We shall use the following growth bound for the regularising flow. Let
$u\in[t_0,T]$, $x\in\R^d$, and let $\mathbf w$ be a continuous input on
$I_u$. With $\mathfrak y_r:=\mathbf w_{u,r}+B^u_r$, the path
$\phi_r^{u,x,\mathbf w}-\mathfrak y_r$ starts at $x$ and has derivative
$F(\phi_r^{u,x,\mathbf w})$, so that \eqref{eq:F-one-sided-global} and
Young's inequality give
$\frac{d}{dr}|\phi_r^{u,x,\mathbf w}-\mathfrak y_r|^2
\le C|\phi_r^{u,x,\mathbf w}-\mathfrak y_r|^2+C|F(\mathfrak y_r)|^2$.
Gr\"onwall's inequality, \eqref{eq:growth}, and $|I_u|\le1$ yield, for
$r\in I_u$,
\begin{align}\label{eq:flow-pathwise-growth}
 |\phi_r^{u,x,\mathbf w}|
 \le C(1+|x|)
 \bigl(1+[\mathbf w]_{C^{H-\eps}(I_u)}+[B^u]_{C^{H-\eps}(I_u)}\bigr)^k.
\end{align}
For $\mathbf w=V+\bar B^{t_0,u}$, the seminorms on the right have moments of
every order by \Cref{ass:input} and \Cref{lem:history-holder-moments}, and
$[B^u]_{C^{H-\eps}(I_u)}$ is independent of $\mathcal F_u$ with Gaussian
moments. Hence, for every $\mathcal F_u$-measurable $\xi$ and every $q<\infty$,
\begin{align}\label{eq:flow-conditional-moments}
 \bigl(\EE^u|\phi_r^{u,\xi,\mathbf w}|^q\bigr)^{1/q}
 \le C(1+|\xi|)
 \bigl(1+[\mathbf w]_{C^{H-\eps}(I_u)}\bigr)^k,
\end{align}
with $C$ depending on $q$, and H\"older's inequality gives
$\|\phi_r^{u,\xi,\mathbf w}\|_{L_{q}}\le C(1+\|\xi\|_{L_{q_1}})$ whenever
$q<q_1$.

We shall repeatedly use the following elementary mixed-norm inequality. If
$2\le q\le2m$ and $\zeta\in L_{2m}$, then
\begin{align}\label{eq:conditional-mixed-moment}
 \left\|\bigl(\EE^s|\zeta|^q\bigr)^{2/q}\right\|_{L_m}
 \le\|\zeta\|_{L_{2m}}^2,
\end{align}
since, after raising both sides to the power $m$, this is
$\EE(\EE^s|\zeta|^q)^{2m/q}\le\EE\EE^s|\zeta|^{2m}=\EE|\zeta|^{2m}$ by
conditional Jensen's inequality. In particular, $q=4$ is admissible because
$m\ge2$.

\begin{lemma}[Conditional dissipativity bound]
\label{lem:stability-dissipativity}
There exist constants $\eta>0$, $h_0\in(0,1]$, and $C\ge1$ such that, for
every interval $\mathbb J\subset[t_0,T]$ of length $h\le h_0$ and every
$\rho\ge1$,
\begin{align}\label{eq:dissipativity-stability-local}
 \mathcal Q_m(\mathbb J)
 \le Ch^\eta\Bigl[
 \Bigl(\sup_{r\in\mathbb J}\|X_r-Y_r\|_{L_{2m}}\Bigr)^2
 +\mathcal K_m(\mathbb J)^2
 \Bigr]
 +C^\rho\rho^\rho h^{2\chi_\eps\rho}.
\end{align}
Above moment order $2m$, only the marginal $L_{2m+\delta_0}$ bound in
\eqref{eq:dissipativity-fixed-bounds} is used.
\end{lemma}

\begin{proof}
If $\mathcal K_m(\mathbb J)=\infty$, there is nothing to prove, so we assume
that $\mathcal K_m(\mathbb J)<\infty$. Fix $(s,t)\in\Delta_{\mathbb J}$ and,
for $r\in[s,t]$, write
$\mathsf P_{s,r}:=1+|X_r|+|Y_r|+|\phi_r^{s,X_s}|+|\phi_r^{s,Y_s}|$.

\smallskip
\emph{Step 1: energy inequality.}
Since $X_r=X_s+\int_s^rF(X_{r'})\,dr'+K^X_{s,r}+V_{s,r}+B^{t_0}_{s,r}$ and,
by \eqref{eq:restart-RL},
$\phi_r^{s,X_s}=X_s+\int_s^rF(\phi_{r'}^{s,X_s})\,dr'+V_{s,r}+B^{t_0}_{s,r}$,
and similarly for $Y$, the process
$\mathfrak D_r:=\mathsf D_{s,r}-\mathbb K_{s,r}$ satisfies $\mathfrak D_s=0$ and
\begin{align*}
 \mathfrak D_r
 =\int_s^r\bigl(F(X_{r'})-F(\phi_{r'}^{s,X_s})
 -F(Y_{r'})+F(\phi_{r'}^{s,Y_s})\bigr)\,dr'.
\end{align*}
In particular, $\mathfrak D$ is absolutely continuous even though $K^X$ and
$K^Y$ need not be, so the energy calculation can be performed directly for
the given solutions. Omitting the time index $r$, set
$\mathfrak p:=Y+\phi^{s,X_s}-X$, so that
$\phi^{s,Y_s}=\mathfrak p+\mathsf D_{s,r}$, and split
\begin{align*}
 F(X)-F(\phi^{s,X_s})-F(Y)+F(\phi^{s,Y_s})
 =\underbrace{F(X)-F(\phi^{s,X_s})-F(Y)+F(\mathfrak p)}_{=:\mathcal T_1}
 +\underbrace{F(\mathfrak p+\mathsf D_{s,r})-F(\mathfrak p)}_{=:\mathcal T_2}.
\end{align*}
Two applications of the fundamental theorem of calculus and
\eqref{eq:growth} give
$|\mathcal T_1|\lesssim|X-\phi^{s,X_s}|\,|X-Y|\,\mathsf P_{s,r}^{(k-2)\vee0}$,
since all the points involved have norm at most $\mathsf P_{s,r}$. By
\eqref{eq:F-one-sided-global} and \eqref{eq:growth},
$\langle\mathsf D_{s,r},\mathcal T_2\rangle\le\kappa_4|\mathsf D_{s,r}|^2$ and
$|\mathcal T_2|\lesssim\mathsf P_{s,r}^{k-1}|\mathsf D_{s,r}|$. Since
$\mathfrak D_r=\mathsf D_{s,r}-\mathbb K_{s,r}$, it follows that
\begin{align*}
 \frac{d}{dr}|\mathfrak D_r|^2
 &=2\langle\mathsf D_{s,r},\mathcal T_2\rangle
 -2\langle \mathbb K_{s,r},\mathcal T_2\rangle
 +2\langle\mathfrak D_r,\mathcal T_1\rangle\\
 &\lesssim|\mathsf D_{s,r}|^2
 +|\mathbb K_{s,r}|\,|\mathsf D_{s,r}|\,\mathsf P_{s,r}^{k-1}
 +\bigl(|\mathsf D_{s,r}|+|\mathbb K_{s,r}|\bigr)
 |X_r-\phi_r^{s,X_s}|\,|X_r-Y_r|\,\mathsf P_{s,r}^{(k-2)\vee0}.
\end{align*}
Integrating from $s$ to $t$ gives
\begin{align}\label{eq:conditional-dissipativity-energy}
\begin{split}
 |\mathsf D_{s,t}-\mathbb K_{s,t}|^2
 &\lesssim\int_s^t|\mathsf D_{s,r}|^2\,dr
 +\int_s^t
 |\mathsf D_{s,r}|\,|X_r-\phi_r^{s,X_s}|\,
 |X_r-Y_r|\mathsf P_{s,r}^{(k-2)\vee0}\,dr\\
 &\quad+\int_s^t
 |\mathbb K_{s,r}|\,|X_r-\phi_r^{s,X_s}|\,
 |X_r-Y_r|\mathsf P_{s,r}^{(k-2)\vee0}\,dr
 +\int_s^t
 |\mathbb K_{s,r}|\,|\mathsf D_{s,r}|\mathsf P_{s,r}^{k-1}\,dr.
\end{split}
\end{align}

\smallskip
\emph{Step 2: the bad event.}
Let $G_{s,t}:=\{t-s\le\tau_{V+\bar B^{t_0,s};I_s}\}$, an
$\mathcal F_s$-measurable event. On $G_{s,t}^c$ we do not use the energy
estimate. By \Cref{lem:reg-phi},
$\ind_{G_{s,t}^c}|\mathsf D_{s,t}|^2
\lesssim\ind_{G_{s,t}^c}(|X_t-Y_t|^2+|X_s-Y_s|^2)$.
Conditional expectation is a contraction on $L_m$, and H\"older's
inequality with $q_0:=2m+\delta_0$ gives
\begin{align}\label{eq:conditional-bad-holder}
\begin{split}
 \left\|\ind_{G_{s,t}^c}\EE^s|\mathsf D_{s,t}|^2
 \right\|_{L_m}
 \lesssim C_0^2
 \PP(G_{s,t}^c)^{2(\frac1{2m}-\frac1{q_0})}.
\end{split}
\end{align}
By \Cref{lem:bad-event-moments}, applied to
$\mathbf w=V+\bar B^{t_0,s}$ on $I_s$ with $\nu=0$ and a moment order chosen
in terms of $\rho$, and after replacing $\rho$ by $\rho+(2\chi_\eps)^{-1}$
to gain an additional factor $t-s$, we obtain
\begin{align}\label{eq:conditional-bad-final}
\begin{split}
 \left\|\ind_{G_{s,t}^c}\EE^s|\mathsf D_{s,t}|^2
 \right\|_{L_m}
 \le C^\rho\rho^\rho
 (t-s)^{1+2\chi_\eps\rho}.
\end{split}
\end{align}
The constants depend on the fixed gap $q_0-2m$. Since the quantity being
estimated is a square, the sub-Gaussian moment bound produces the growth
$C^\rho\rho^\rho$.

\smallskip
\emph{Step 3: the good event.}
On $G_{s,t}$, for every $r\in[s,t]$ and every $p\in[2,p_k]$, the second bound
in \eqref{eq:dissipativity-fixed-bounds} and the coming-down estimate
\eqref{eq:coming-down} (for the flows, uniformly in the initial point) give
\begin{align}\label{eq:conditional-good-moments}
\begin{split}
 &\ind_{G_{s,t}}
 \bigl(\EE^s|X_r-\phi_r^{s,X_s}|^p\bigr)^{1/p}
 +\ind_{G_{s,t}}
 \bigl(\EE^s|Y_r-\phi_r^{s,Y_s}|^p\bigr)^{1/p}
 \lesssim(r-s)^{\alpha_\gamma},\\
 &\ind_{G_{s,t}}
 \bigl(\EE^s\mathsf P_{s,r}^p\bigr)^{1/p}
 \lesssim1+(r-s)^{-1/(k-1)}.
\end{split}
\end{align}
Conditional H\"older's inequality in \eqref{eq:conditional-good-moments}
gives
\begin{align}\label{eq:conditional-good-weighted}
\begin{split}
 &\ind_{G_{s,t}}
 \bigl(\EE^s[
 |X_r-\phi_r^{s,X_s}|^4\mathsf P_{s,r}^{4((k-2)\vee0)}]
 \bigr)^{1/4}
 \lesssim(r-s)^{\alpha_\gamma-((k-2)\vee0)/(k-1)},\\
 &\ind_{G_{s,t}}
 \bigl(\EE^s\mathsf P_{s,r}^{4(k-1)}\bigr)^{1/4}
 \lesssim(r-s)^{-1},
\end{split}
\end{align}
and all the conditional moments used here have order at most $p_k$.

We now estimate each term in \eqref{eq:conditional-dissipativity-energy}.
Conditional H\"older's inequality with exponents $(2,4,4)$, followed by
\eqref{eq:conditional-good-weighted}, gives
\begin{align}\label{eq:conditional-first-product}
\begin{split}
 &\ind_{G_{s,t}}\EE^s\left[
 |\mathsf D_{s,r}|\,|X_r-\phi_r^{s,X_s}|\,
 |X_r-Y_r|\mathsf P_{s,r}^{(k-2)\vee0}\right]\\
 &\quad\lesssim(r-s)^{\alpha_\gamma-((k-2)\vee0)/(k-1)}
 \bigl(\EE^s|\mathsf D_{s,r}|^2\bigr)^{1/2}
 \bigl(\EE^s|X_r-Y_r|^4\bigr)^{1/4}.
\end{split}
\end{align}
Conditional H\"older's inequality with exponents $(4,4,4)$ similarly yields
\begin{align*}
\begin{split}
 &\ind_{G_{s,t}}\EE^s\left[
 |\mathbb K_{s,r}|\,|X_r-\phi_r^{s,X_s}|\,
 |X_r-Y_r|\mathsf P_{s,r}^{(k-2)\vee0}\right]\\
 &\quad\lesssim(r-s)^{\alpha_\gamma-((k-2)\vee0)/(k-1)}
 \bigl(\EE^s|\mathbb K_{s,r}|^4\bigr)^{1/4}
 \bigl(\EE^s|X_r-Y_r|^4\bigr)^{1/4},
\end{split}
\end{align*}
and conditional H\"older's inequality with exponents $(4,2,4)$ gives
\begin{align}
\label{eq:conditional-third-product}
 \ind_{G_{s,t}}\EE^s\left[
 |\mathbb K_{s,r}|\,|\mathsf D_{s,r}|\mathsf P_{s,r}^{k-1}\right] &\lesssim(r-s)^{-1}
 \bigl(\EE^s|\mathbb K_{s,r}|^4\bigr)^{1/4}
 \bigl(\EE^s|\mathsf D_{s,r}|^2\bigr)^{1/2}.
\end{align}
By \eqref{eq:conditional-mixed-moment},
$\|(\EE^s|X_r-Y_r|^4)^{1/4}\|_{L_{2m}}\le\|X_r-Y_r\|_{L_{2m}}$ and
$\|(\EE^s|\mathbb K_{s,r}|^4)^{1/4}\|_{L_{2m}}\le\|\mathbb K_{s,r}\|_{L_{2m}}
\le\mathcal K_m(\mathbb J)(r-s)^{\beta_\gamma}$. Taking the $L_m$-norm in
\eqref{eq:conditional-first-product}--\eqref{eq:conditional-third-product},
using H\"older's inequality, the definition of $\mathcal Q_m(\mathbb J)$,
and $2ab\le a^2+b^2$, and integrating in $r$, we obtain
\begin{align}\label{eq:conditional-good-integrated}
\begin{split}
 &\left\|\ind_{G_{s,t}}\EE^s|\mathsf D_{s,t}-\mathbb K_{s,t}|^2
 \right\|_{L_m}\\
 &\quad\le C(t-s)\Bigg[
 h\,\mathcal Q_m(\mathbb J)
 +h^{\alpha_\gamma+1/2-((k-2)\vee0)/(k-1)}
 \left(\mathcal Q_m(\mathbb J)
 +\sup_{r\in\mathbb J}\|X_r-Y_r\|_{L_{2m}}^2\right)\\
 &\qquad+
 h^{\alpha_\gamma+\beta_\gamma-((k-2)\vee0)/(k-1)}
 \left(
 \sup_{r\in\mathbb J}\|X_r-Y_r\|_{L_{2m}}^2
 +\mathcal K_m(\mathbb J)^2\right)
 +
 h^{\beta_\gamma-1/2}
 \left(\mathcal Q_m(\mathbb J)+\mathcal K_m(\mathbb J)^2\right)
 \Bigg].
\end{split}
\end{align}
For example, the first product contributes
\begin{align*}
 &\int_s^t
 (r-s)^{\alpha_\gamma-((k-2)\vee0)/(k-1)}
 \|\EE^s|\mathsf D_{s,r}|^2\|_{L_m}^{1/2}
 \|X_r-Y_r\|_{L_{2m}}\,dr\\
 &\quad\le
 C(t-s)h^{\alpha_\gamma+1/2-((k-2)\vee0)/(k-1)}
 \left(\mathcal Q_m(\mathbb J)
 +\sup_{r\in\mathbb J}\|X_r-Y_r\|_{L_{2m}}^2\right),
\end{align*}
and the other terms are treated in the same way.

\smallskip
\emph{Step 4: conclusion.}
The quantity $\mathcal Q_m(\mathbb J)$ is finite: on the good event this
follows from \eqref{eq:conditional-good-moments} and $2\alpha_\gamma>1$,
and on the bad event from \eqref{eq:conditional-bad-final} with one fixed
value of $\rho$. Since
$\|\EE^s|\mathsf D_{s,t}|^2\|_{L_m}\lesssim
\|\EE^s|\mathsf D_{s,t}-\mathbb K_{s,t}|^2\|_{L_m}+\|\mathbb K_{s,t}\|_{L_{2m}}^2$
and $\|\mathbb K_{s,t}\|_{L_{2m}}^2\le\mathcal K_m(\mathbb J)^2(t-s)h^{2\beta_\gamma-1}$,
we set
\begin{align*}
 \eta:=\min\left\{
 1,\,
 \alpha_\gamma+\frac12-\frac{(k-2)\vee0}{k-1},\,
 \alpha_\gamma+\beta_\gamma-\frac{(k-2)\vee0}{k-1},\,
 \beta_\gamma-\frac12,\,
 2\beta_\gamma-1
 \right\},
\end{align*}
which is positive by \eqref{eq:exponent-relations} and since
$(k-2)\vee0<k-1$ and $\alpha_\gamma>1/2$. Dividing
\eqref{eq:conditional-good-integrated} and
\eqref{eq:conditional-bad-final} by $t-s$, taking the supremum over
$(s,t)\in\Delta_{\mathbb J}$, and choosing $h_0$ so that
$Ch_0^{\eta}\le1/2$, we absorb $\mathcal Q_m(\mathbb J)$ into the left-hand
side. This proves \eqref{eq:dissipativity-stability-local}.
\end{proof}

We next bound $\mathcal K_m$ by stochastic sewing.

\begin{lemma}[Sewing bound]\label{lem:stability-sewing}
Let $h_0$ be as in \Cref{lem:stability-dissipativity}. There exists a
constant $C$ such that, for every interval $\mathbb J\subset[t_0,T]$ of
length $h\le h_0$ and every $\rho\ge1$,
\begin{align}\label{eq:K-stability-local}
 \mathcal K_m(\mathbb J)
 \le C\Bigl(
 \sup_{r\in\mathbb J}\|X_r-Y_r\|_{L_{2m}}
 +\|b-g\|_{\mathcal C^{\gamma-1}}
 +h^{1/2}\mathcal Q_m(\mathbb J)^{1/2}
 +C^\rho\rho^{\rho/2}h^{\chi_\eps\rho}\Bigr).
\end{align}
\end{lemma}

\begin{proof}
Let $b^N=P_{1/N}b$ and $g^N=P_{1/N}g$. Uniformly in $N$,
$\|b^N\|_{\mathcal C^\gamma}+\|g^N\|_{\mathcal C^\gamma}\lesssim\Xi$ and
$\|b^N-g^N\|_{\mathcal C^{\gamma-1}}\lesssim\|b-g\|_{\mathcal C^{\gamma-1}}$.
For $u<r$ in $\mathbb J$, define
$\mathfrak G^N_{u,r}:=\EE^u[b^N(\phi_r^{u,X_u})-g^N(\phi_r^{u,Y_u})]$.
The sewing germ is
\begin{align}\label{eq:stability-sewing-germ}
 A_{u,v}:=\int_u^v\mathfrak G^N_{u,r}\,dr.
\end{align}

Split the integrand into
$b^N(\phi_r^{u,X_u})-b^N(\phi_r^{u,Y_u})+(b^N-g^N)(\phi_r^{u,Y_u})$.
With $\mathbf w:=V+\bar B^{t_0,u}$ frozen, on the event
$\{r-u\le\tau_{\mathbf w;I_u}\}$,
\Cref{lem:malliavin-2}, \eqref{eq:odd-one-point-IBP} with
$\sigma=1-\gamma$, and \Cref{prop:flow-estimates} give
\begin{align*}
 &\left|\EE\!\left[
 b^N(\phi_r^{u,x,\mathbf w})-
 b^N(\phi_r^{u,y,\mathbf w})\right]\right|
 \lesssim|x-y|(r-u)^{-(1-\gamma)H},\\
 &\left|\EE\!\left[
 (b^N-g^N)(\phi_r^{u,y,\mathbf w})\right]\right|
 \lesssim\|b-g\|_{\mathcal C^{\gamma-1}}
 (r-u)^{-(1-\gamma)H}.
\end{align*}
Substituting $x=X_u,y=Y_u$ and integrating the time singularity therefore
gives
\begin{align*}
\begin{split}
 &\left\|\int_u^v
 \ind_{\{r-u\le\tau_{V+\bar B^{t_0,u};I_u}\}}
 \mathfrak G^N_{u,r}\,dr
 \right\|_{L_{2m}}\\
 &\quad\lesssim
 \left(\sup_{r\in\mathbb J}\|X_r-Y_r\|_{L_{2m}}
 +\|b-g\|_{\mathcal C^{\gamma-1}}\right)
 (v-u)^{\beta_\gamma}.
\end{split}
\end{align*}
On the complementary event, we return to the original two evaluations
defining $\mathfrak G^N_{u,r}$.  The one-point estimate
\Cref{lem:malliavin-1}, together with
\eqref{eq:Lambda-flow-scale}, gives the two separate bounds
\begin{align*}
 \left|\EE^u b^N(
 \phi_r^{u,X_u})\right|
 &\lesssim
 \bigl(R_{u,r}^{V+\bar B^{t_0,u};I_u}\bigr)^C
 (r-u)^{(\gamma\wedge0)H}\\
 &\le
 \bigl(R_{u,r}^{V+\bar B^{t_0,u};I_u}\bigr)^C
 (r-u)^{-(1-\gamma)H},\\
 \left|\EE^u g^N(
 \phi_r^{u,Y_u})\right|
 &\lesssim
 \bigl(R_{u,r}^{V+\bar B^{t_0,u};I_u}\bigr)^C
 (r-u)^{(\gamma\wedge0)H}\\
 &\le
 \bigl(R_{u,r}^{V+\bar B^{t_0,u};I_u}\bigr)^C
 (r-u)^{-(1-\gamma)H}.
\end{align*}
Here the last inequalities use $r-u\le1$.  The triangle inequality
therefore yields
\begin{align*}
 |\mathfrak G^N_{u,r}|
 \lesssim
 \bigl(R_{u,r}^{V+\bar B^{t_0,u};I_u}\bigr)^C
 (r-u)^{-(1-\gamma)H}.
\end{align*}
Consequently, \Cref{lem:bad-event-moments} gives
\begin{align*}
 \left\|\int_u^v
 \ind_{\{r-u>\tau_{V+\bar B^{t_0,u};I_u}\}}
 \mathfrak G^N_{u,r}\,dr
 \right\|_{L_{2m}}
 \le C^\rho\rho^{\rho/2}
 (v-u)^{\beta_\gamma+\chi_\eps\rho}.
\end{align*}
Combining the two events,
\begin{align}\label{eq:stability-germ-bound}
\begin{split}
 \|A_{u,v}\|_{L_{2m}}
 &\lesssim
 \left(\sup_{r\in\mathbb J}\|X_r-Y_r\|_{L_{2m}}
 +\|b-g\|_{\mathcal C^{\gamma-1}}\right)
 (v-u)^{\beta_\gamma}\\
 &\quad+C^\rho\rho^{\rho/2}
 (v-u)^{\beta_\gamma+\chi_\eps\rho}.
\end{split}
\end{align}
In particular, for $c\in(u,v)$, the triangle inequality and
\eqref{eq:stability-germ-bound} give the second sewing estimate
\begin{align}\label{eq:stability-unconditional-defect}
\begin{split}
 \|\delta A_{u,c,v}\|_{L_{2m}}
 &\le\|A_{u,v}\|_{L_{2m}}
 +\|A_{u,c}\|_{L_{2m}}
 +\|A_{c,v}\|_{L_{2m}}\\
 &\lesssim\Bigg(
 \sup_{r\in\mathbb J}\|X_r-Y_r\|_{L_{2m}}
 +\|b-g\|_{\mathcal C^{\gamma-1}}
 +C^\rho\rho^{\rho/2}h^{\chi_\eps\rho}
 \Bigg)(v-u)^{\beta_\gamma}.
\end{split}
\end{align}
Since $\beta_\gamma>1/2$, this has the exponent required in
\eqref{eq:stochastic-sewing-defect}.

We next estimate the conditional defect at the midpoint $c=(u+v)/2$. The
tower property and the flow property, together with
\eqref{eq:history-cocycle}, give
\begin{align*}
 \EE^u\delta A_{u,c,v}
 =\EE^u\int_c^v\EE^c\Big[
 b^N\big(\phi_r^{c,\phi_c^{u,X_u}}\big)
 -g^N\big(\phi_r^{c,\phi_c^{u,Y_u}}\big)
 -b^N(\phi_r^{c,X_c})
 +g^N(\phi_r^{c,Y_c})
 \Big]\,dr,
\end{align*}
where, following our convention, a flow starting at $c$ carries the input
$V+\bar B^{t_0,c}$ and a flow starting at $u$ carries $V+\bar B^{t_0,u}$.
Define the two integrands
\begin{align*}
 \mathcal J^N_{1;u,c,r}
 &:=b^N\big(\phi_r^{c,\phi_c^{u,X_u}}\big)
 -b^N(\phi_r^{c,X_c})
 -b^N\big(\phi_r^{c,\phi_c^{u,Y_u}}\big)
 +b^N(\phi_r^{c,Y_c}),\\
 \mathcal J^N_{2;u,c,r}
 &:=(b^N-g^N)\big(\phi_r^{c,\phi_c^{u,Y_u}}\big)
 -(b^N-g^N)(\phi_r^{c,Y_c}).
\end{align*}
Then
\begin{align*}
 \EE^u\delta A_{u,c,v}
 =\EE^u\int_c^v\EE^c
 \bigl[\mathcal J^N_{1;u,c,r}+\mathcal J^N_{2;u,c,r}\bigr]\,dr.
\end{align*}
We now decompose these two terms before estimating them. Put
$\xi^X:=\phi_c^{u,X_u}$ and $\xi^Y:=\phi_c^{u,Y_u}$, which are the
restarted flows at time $c$, and, for $z\in\{\xi^X,\xi^Y,X_c,Y_c\}$,
write $\Psi_r^z:=\phi_r^{c,z}$. Then
\begin{align*}
 \mathcal J^N_{1;u,c,r}
 &=b^N(\Psi_r^{\xi^X})-b^N(\Psi_r^{X_c})
   -b^N(\Psi_r^{\xi^Y})+b^N(\Psi_r^{Y_c})\\
 &=\Big[
 b^N(\Psi_r^{\xi^X})-b^N(\Psi_r^{X_c})
 -b^N(\Psi_r^{\xi^Y})
 +b^N(\Psi_r^{X_c}+\Psi_r^{\xi^Y}-\Psi_r^{\xi^X})
 \Big]\\
 &\quad+
 \Big[
 b^N(\Psi_r^{Y_c})
 -b^N(\Psi_r^{X_c}+\Psi_r^{\xi^Y}-\Psi_r^{\xi^X})
 \Big].
\end{align*}
\Cref{lem:malliavin-3}, with $Z_0=\Psi_r^{\xi^X}$, $Z_1=\Psi_r^{X_c}$, and
$Z_2=\Psi_r^{\xi^Y}$, applies to the first bracket; by
\eqref{eq:upperbound-D-x}, its two differences are controlled by
$|\xi^X-X_c|$ and $|\xi^X-\xi^Y|$. For the second bracket, \Cref{lem:malliavin-2} applies and
\eqref{eq:upperbound-4} gives
\begin{align*}
 &\left\|
 \Psi_r^{\xi^X}-\Psi_r^{\xi^Y}-\Psi_r^{X_c}+\Psi_r^{Y_c}
 \right\|_{\mathbb D^{\lfloor-\gamma\rfloor+3,\bar q}}\\
 &\quad\lesssim
 \bigl(R_{c,r}^{V+\bar B^{t_0,c};I_c}\bigr)^C
 \Big(
 |\xi^X-\xi^Y-X_c+Y_c|+|\xi^X-\xi^Y|\,|\xi^X-X_c|
 \Big).
\end{align*}
It follows that
\begin{align*}
 \left|\EE^c\mathcal J^N_{1;u,c,r}\right|
 &\lesssim
 \bigl(R_{c,r}^{V+\bar B^{t_0,c};I_c}\bigr)^C
 |\xi^X-X_c|\,|\xi^X-\xi^Y|(r-c)^{-(2-\gamma)H}\\
 &\quad+
 \bigl(R_{c,r}^{V+\bar B^{t_0,c};I_c}\bigr)^C
 \Big(
 |\xi^X-\xi^Y-X_c+Y_c|+|\xi^X-\xi^Y|\,|\xi^X-X_c|
 \Big)(r-c)^{-(1-\gamma)H}.
\end{align*}
Moreover,
$\mathcal J^N_{2;u,c,r}=(b^N-g^N)(\Psi_r^{\xi^Y})-(b^N-g^N)(\Psi_r^{Y_c})$.
Since
\(\|b^N-g^N\|_{\mathcal C^{\gamma-1}}
\lesssim\|b-g\|_{\mathcal C^{\gamma-1}}\),
\eqref{eq:odd-two-point-IBP} and
\eqref{eq:upperbound-D-x} give
\begin{align*}
 \left|\EE^c\mathcal J^N_{2;u,c,r}\right|
 \lesssim
 \bigl(R_{c,r}^{V+\bar B^{t_0,c};I_c}\bigr)^C
 \|b-g\|_{\mathcal C^{\gamma-1}}
 |\xi^Y-Y_c|(r-c)^{-(2-\gamma)H}.
\end{align*}
Adding the last two inequalities yields
\begin{align}\label{eq:stability-direct-defect-bound}
\begin{split}
 &\left|\EE^c\big[
 \mathcal J^N_{1;u,c,r}+\mathcal J^N_{2;u,c,r}\big]\right|\\
 &\lesssim
 (R_{c,r}^{V+\bar B^{t_0,c};I_c})^C
 |\phi_c^{u,X_u}-X_c|\,
 |\phi_c^{u,X_u}-\phi_c^{u,Y_u}|
 (r-c)^{-(2-\gamma)H}\\
 &\quad+
 (R_{c,r}^{V+\bar B^{t_0,c};I_c})^C
 |\phi_c^{u,X_u}-\phi_c^{u,Y_u}-X_c+Y_c|
 (r-c)^{-(1-\gamma)H}\\
 &\quad+
 (R_{c,r}^{V+\bar B^{t_0,c};I_c})^C
 |\phi_c^{u,X_u}-X_c|\,
 |\phi_c^{u,X_u}-\phi_c^{u,Y_u}|
 (r-c)^{-(1-\gamma)H}\\
 &\quad+
 (R_{c,r}^{V+\bar B^{t_0,c};I_c})^C
 \|b-g\|_{\mathcal C^{\gamma-1}}
 |\phi_c^{u,Y_u}-Y_c|
 (r-c)^{-(2-\gamma)H}.
\end{split}
\end{align}
Again, all the flows in this display carry the common input
$V+\bar B^{t_0,u}$ up to time $c$.

Set $G_1:=\{c-u\le\tau_{V+\bar B^{t_0,u};I_u}\}$ and
$G_2(r):=\{r-c\le\tau_{V+\bar B^{t_0,c};I_c}\}$. Since the comparison
class is defined with the full frozen input $V+\bar B^{t_0,u}$, the local
scale of this sum is used here. The three indicators
$\ind_{G_1\cap G_2(r)}$, $\ind_{G_1^c\cap G_2(r)}$, and
$\ind_{G_2(r)^c}$ form a partition of unity.
For the latter two terms, expand
\begin{align*}
 \mathcal J^N_{1;u,c,r}+\mathcal J^N_{2;u,c,r}
 =b^N(\Psi_r^{\xi^X})-g^N(\Psi_r^{\xi^Y})
 -b^N(\Psi_r^{X_c})+g^N(\Psi_r^{Y_c}).
\end{align*}
Applying \Cref{lem:malliavin-1} separately to these four evaluations
and using \eqref{eq:Lambda-flow-scale} gives
\begin{align}\label{eq:stability-defect-coarse}
 \big|\EE^c[\mathcal J^N_{1;u,c,r}+\mathcal J^N_{2;u,c,r}]\big|
 \lesssim(R_{c,r}^{V+\bar B^{t_0,c};I_c})^C(r-c)^{\beta_\gamma-1}.
\end{align}
Combining this bound with
\eqref{eq:stability-direct-defect-bound}, and using $r-c\le1$, yields
\begin{align}\label{eq:stability-defect-indicators}
\begin{split}
 &\big|\EE^c[\mathcal J^N_{1;u,c,r}+\mathcal J^N_{2;u,c,r}]\big|\\
 &\lesssim\ind_{G_1\cap G_2(r)}
 \Big[|\xi^X-\xi^Y-X_c+Y_c|(r-c)^{\beta_\gamma-1}\\
 &\qquad\quad+
 \big(|\xi^X-X_c|\,|\xi^X-\xi^Y|
 +\|b-g\|_{\mathcal C^{\gamma-1}}|\xi^Y-Y_c|\big)
 (r-c)^{\zeta_\gamma-1}\Big]\\
 &\quad+\Big[\ind_{G_1^c\cap G_2(r)}
 +\ind_{G_2(r)^c}(R_{c,r}^{V+\bar B^{t_0,c};I_c})^C\Big]
 (r-c)^{\beta_\gamma-1}.
\end{split}
\end{align}
For the first line, since
$\xi^X-\xi^Y-X_c+Y_c=-\mathsf D_{u,c}$, conditional Jensen and the definition
of $\mathcal Q_m(\mathbb J)$ give
\begin{align*}
 \big\|\EE^u[\ind_{G_1\cap G_2(r)}
 |\xi^X-\xi^Y-X_c+Y_c|]\big\|_{L_{2m}}
 \le(c-u)^{1/2}\mathcal Q_m(\mathbb J)^{1/2}.
\end{align*}
For the products, $|\xi^X-\xi^Y|\le e^{\kappa_4(c-u)}|X_u-Y_u|$ by
\Cref{lem:reg-phi}, which is $\mathcal F_u$-measurable, while
$|\xi^X-X_c|$ and $|\xi^Y-Y_c|$ are the comparison remainders of $X$ and $Y$
over $[u,c]$. The second bound in \eqref{eq:dissipativity-fixed-bounds} on
$G_1$ therefore gives
\begin{align*}
 &\big\|\EE^u[\ind_{G_1\cap G_2(r)}
 (|\xi^X-X_c|\,|\xi^X-\xi^Y|
 +\|b-g\|_{\mathcal C^{\gamma-1}}|\xi^Y-Y_c|)]\big\|_{L_{2m}}\\
 &\qquad\lesssim(c-u)^{\alpha_\gamma}
 \left(\sup_{r\in\mathbb J}\|X_r-Y_r\|_{L_{2m}}
 +\|b-g\|_{\mathcal C^{\gamma-1}}\right).
\end{align*}
For the last line of \eqref{eq:stability-defect-indicators},
\Cref{lem:bad-event-moments} gives
\begin{align*}
 \|\ind_{G_1^c}\|_{L_{2m}}
 &\le C^\rho\rho^{\rho/2}(c-u)^{\chi_\eps\rho},\\
 \|\ind_{G_2(r)^c}(R_{c,r}^{V+\bar B^{t_0,c};I_c})^C\|_{L_{2m}}
 &\le C^\rho\rho^{\rho/2}(r-c)^{\chi_\eps\rho}.
\end{align*}
Integrating in $r$, we obtain
\begin{align}\label{eq:stability-conditional-defect}
\begin{split}
 \|\EE^u\delta A_{u,c,v}\|_{L_{2m}}
 &\lesssim
 (c-u)^{1/2}(v-c)^{\beta_\gamma}
 \mathcal Q_m(\mathbb J)^{1/2}\\
 &\quad+
 (c-u)^{\alpha_\gamma}(v-c)^{\zeta_\gamma}
 \left(\sup_{r\in\mathbb J}\|X_r-Y_r\|_{L_{2m}}
 +\|b-g\|_{\mathcal C^{\gamma-1}}\right)\\
 &\quad+C^\rho\rho^{\rho/2}
 (v-u)^{\beta_\gamma+\chi_\eps\rho}.
\end{split}
\end{align}
The three sewing powers are admissible because
$\frac12+\beta_\gamma>1$ and
$\alpha_\gamma+\zeta_\gamma>1$ by \eqref{eq:exponent-relations}, and $\rho$ can be chosen so that
$\beta_\gamma+\chi_\eps\rho>1$. To obtain the estimate for every prescribed
$\rho\ge1$, apply the preceding bounds with $\rho+\rho_0$, where the fixed
$\rho_0$ makes the last exponent larger than one, and absorb the shift into
$C$.

Let $S$ be the left endpoint of $\mathbb J$. The process
$\mathcal Y^N_t:=\int_S^t(b^N(X_r)-g^N(Y_r))\,dr$ is the process associated
with the germ \eqref{eq:stability-sewing-germ} by
\Cref{lem:stochastic-sewing}. Since $X$ and $Y$ have singular drifts, we
check the characterisation in that lemma. We first show that there exist
$\varsigma_0>0$ and $C$ such that, for $Z\in\{X,Y\}$ and $u<r$ in $\mathbb J$,
\begin{align}\label{eq:comparison-remainder-L2m}
 \|Z_r-\phi_r^{u,Z_u}\|_{L_{2m}}\le C(r-u)^{\varsigma_0}.
\end{align}
Let $G:=\{r-u\le\tau_{V+\bar B^{t_0,u};I_u}\}$ and $q_1:=2m+\delta_0/2$. By
\eqref{eq:flow-conditional-moments} and the marginal bound in
\eqref{eq:dissipativity-fixed-bounds}, $\|Z_r-\phi_r^{u,Z_u}\|_{L_{q_1}}\le C$.
On $G$, the comparison bound in \eqref{eq:dissipativity-fixed-bounds} with
$p=2$ gives $\|\ind_G(Z_r-\phi_r^{u,Z_u})\|_{L_2}\le C_0(r-u)^{\alpha_\gamma}$,
and interpolation between $L_2$ and $L_{q_1}$ yields
$\|\ind_G(Z_r-\phi_r^{u,Z_u})\|_{L_{2m}}\le C(r-u)^{\vartheta\alpha_\gamma}$,
where $\vartheta\in(0,1]$ is defined by
$\frac1{2m}=\frac\vartheta2+\frac{1-\vartheta}{q_1}$. On $G^c$, H\"older's
inequality and \Cref{lem:bad-event-moments}, with $\rho$ large enough, give
$\|\ind_{G^c}(Z_r-\phi_r^{u,Z_u})\|_{L_{2m}}
\le C\PP(G^c)^{\frac1{2m}-\frac1{q_1}}\le C(r-u)^{\chi_\eps}$. This proves
\eqref{eq:comparison-remainder-L2m} with
$\varsigma_0:=\vartheta\alpha_\gamma\wedge\chi_\eps$. Now write
$b^N(X_r)-\EE^u b^N(\phi_r^{u,X_u})
=[b^N(X_r)-b^N(\phi_r^{u,X_u})]
+[b^N(\phi_r^{u,X_u})-\EE^u b^N(\phi_r^{u,X_u})]$, and similarly for
$g^N$ and $Y$. The first bracket is bounded by
$\|\nabla b^N\|_\infty|X_r-\phi_r^{u,X_u}|$, and the second has zero
conditional mean given $\mathcal F_u$ and is bounded by $2\|b^N\|_\infty$.
Integrating over $[u,v]$ gives
$\|\mathcal Y^N_v-\mathcal Y^N_u-A_{u,v}\|_{L_{2m}}\le C_N(v-u)$ and
$\|\EE^u(\mathcal Y^N_v-\mathcal Y^N_u-A_{u,v})\|_{L_{2m}}
\le C_N(v-u)^{1+\varsigma_0}$ by \eqref{eq:comparison-remainder-L2m}, with
constants depending on $N$, which is admissible in that characterisation.
Hence
\Cref{lem:stochastic-sewing}, together
with \eqref{eq:stability-germ-bound}, \eqref{eq:stability-unconditional-defect},
and \eqref{eq:stability-conditional-defect}, yield
\begin{align}\label{eq:stability-smooth-occupation}
\begin{split}
 &\sup_{(u,v)\in\Delta_{\mathbb J}}
 \frac{\left\|\int_u^v(b^N(X_r)-g^N(Y_r))\,dr\right\|_{L_{2m}}}
 {(v-u)^{\beta_\gamma}}\\
 &\quad\lesssim
 \sup_{r\in\mathbb J}\|X_r-Y_r\|_{L_{2m}}
 +\|b-g\|_{\mathcal C^{\gamma-1}}\\
 &\qquad+
 h^{1/2}\mathcal Q_m(\mathbb J)^{1/2}
 +C^\rho\rho^{\rho/2}h^{\chi_\eps\rho}.
\end{split}
\end{align}
By \Cref{def:solution} and \Cref{ass:singular-drift}, a subsequence of these
smooth occupation integrals converges almost surely, uniformly in time, to
$\mathbb K=K^X-K^Y$.
Fatou's lemma in \eqref{eq:stability-smooth-occupation} proves
\eqref{eq:K-stability-local}.
\end{proof}

Finally, we derive an energy inequality for the difference of the two
solutions, in which the singular drift enters only through a process that
can be estimated by stochastic sewing. Set $Z:=X-Y$,
$w_t:=e^{-2\kappa_4(t-t_0)}$, $\mathfrak F_t:=\int_{t_0}^t(F(X_r)-F(Y_r))\,dr$,
$U^0_t:=Z_{t_0}+\mathfrak F_t$, and $U_t:=w_tU^0_t$, so that $Z=U^0+\mathbb K$. The
paths $\mathfrak F$, $U^0$, and $U$ have finite variation almost surely; see
the beginning of the proof below. We define the continuous scalar process
$\mathcal E$, with $\mathcal E_{t_0}=0$, by
\begin{align}
\label{eq:energy-definition}
 \mathcal E_t
 &:=
 \langle U_t,\mathbb K_t\rangle
 -\int_{t_0}^t\langle \mathbb K_r,dU_r\rangle +\frac12w_t|\mathbb K_t|^2
 -\frac12\int_{t_0}^t|\mathbb K_r|^2\,dw_r,
\end{align}
where both integrals are Lebesgue--Stieltjes integrals against
finite-variation paths. When $\mathbb K$ has finite variation, integration by parts
shows that $\mathcal E_t-\mathcal E_s=\int_s^tw_r\langle Z_r,d\mathbb K_r\rangle$;
in general, \eqref{eq:energy-definition} gives a meaning to this integral.

\begin{lemma}[Energy sewing bound]\label{lem:stability-energy-sewing}
There exists a constant $C$ such that, for every interval
$\mathbb J\subset[t_0,T]$ of length at most $h_0$, every $(s,t)\in\Delta_{\mathbb J}$,
and every $\rho\ge1$,
\begin{align}\label{eq:energy-sewing-local}
\begin{split}
 \|\mathcal E_t-\mathcal E_s\|_{L_m}
 &\le C(t-s)^{\beta_\gamma}
 \left[
 \left(\sup_{r\in\mathbb J}\|Z_r\|_{L_{2m}}\right)^2
 +\|b-g\|_{\mathcal C^{\gamma-1}}^2
 \right]\\
 &\quad+C(t-s)^{\beta_\gamma+1/2}
 \mathcal Q_m(\mathbb J)^{1/2}
 \left(
 \sup_{r\in\mathbb J}\|Z_r\|_{L_{2m}}
 +\|b-g\|_{\mathcal C^{\gamma-1}}\right)
 +C^\rho\rho^{\rho/2}
 (t-s)^{\beta_\gamma+\chi_\eps\rho}.
\end{split}
\end{align}
Moreover, the difference $Z=X-Y$ satisfies the energy inequality
\begin{align}\label{eq:singular-energy-identity}
\begin{split}
 e^{-2\kappa_4(t-t_0)}|Z_t|^2
 \le e^{-2\kappa_4(s-t_0)}|Z_s|^2
 +2(\mathcal E_t-\mathcal E_s),
 \qquad (s,t)\in\Delta_{[t_0,T]}.
\end{split}
\end{align}
\end{lemma}

\begin{proof}
Since $X$ and $Y$ have continuous paths and $F$ is continuous, the path
$r\mapsto F(X_r)-F(Y_r)$ is continuous almost surely. Consequently,
$\mathfrak F$ is absolutely continuous, with
$\operatorname{Var}_{[t_0,T]}(\mathfrak F)
\le\int_{t_0}^T|F(X_r)-F(Y_r)|\,dr<\infty$
almost surely. Since $w$ is continuously differentiable, $U=wU^0$ also has
finite variation.  Thus all the terms in
\eqref{eq:energy-definition} are well defined.

Let $b^N=P_{1/N}b$, $g^N=P_{1/N}g$, and
$\mathbb K_t^N:=\int_{t_0}^t(b^N(X_r)-g^N(Y_r))\,dr$. For $u<r$, abbreviate
$\Psi_r^{u,X}:=\phi_r^{u,X_u}$ and $\Psi_r^{u,Y}:=\phi_r^{u,Y_u}$ (with the
frozen input $V+\bar B^{t_0,u}$). Consider the scalar germ
\begin{align}\label{eq:energy-sewing-germ}
\begin{split}
 A^N_{u,v}
 :=\int_u^v e^{-2\kappa_4(r-t_0)}\EE^u
 \left[
 \left\langle
 \Psi_r^{u,X}-\Psi_r^{u,Y},
 b^N(\Psi_r^{u,X})-g^N(\Psi_r^{u,Y})
 \right\rangle
 \right]dr.
\end{split}
\end{align}
We first bound the germ.  Split its integrand into
\begin{align*}
 \left\langle
 \Psi_r^{u,X}-\Psi_r^{u,Y},
 b^N(\Psi_r^{u,X})-b^N(\Psi_r^{u,Y})
 \right\rangle +
 \left\langle
 \Psi_r^{u,X}-\Psi_r^{u,Y},
 (b^N-g^N)(\Psi_r^{u,Y})
 \right\rangle.
\end{align*}
On
$\{r-u\le\tau_{V+\bar B^{t_0,u};I_u}\}$,
the multiplier $\Psi_r^{u,X}-\Psi_r^{u,Y}$ and the difference between the
two arguments in the first term both have Malliavin--Sobolev norm bounded by
$C|Z_u|$, by \eqref{eq:upperbound-D-x}. Applying \Cref{lem:malliavin-2} to
the first term and \eqref{eq:odd-one-point-IBP} with $\sigma=1-\gamma$ to the
second term, we obtain
\begin{align*}
\begin{split}
 &\left|\EE^u\left[
 \left\langle
 \Psi_r^{u,X}-\Psi_r^{u,Y},
 b^N(\Psi_r^{u,X})-g^N(\Psi_r^{u,Y})
 \right\rangle\right]\right|\\
 &\quad\lesssim
 \left(
 |Z_u|^2+
 \|b-g\|_{\mathcal C^{\gamma-1}}|Z_u|
 \right)(r-u)^{-(1-\gamma)H}.
\end{split}
\end{align*}
Taking the $L_m$-norm, using
$2ab\le a^2+b^2$, and integrating in $r$ yields
\begin{align}\label{eq:energy-germ-good}
\begin{split}
 \left\|\ind_{\{v-u\le\tau_{V+\bar B^{t_0,u};I_u}\}}
 A^N_{u,v}\right\|_{L_m}
 \lesssim
 \left[
 \left(\sup_{r\in\mathbb J}\|Z_r\|_{L_{2m}}\right)^2
 +\|b-g\|_{\mathcal C^{\gamma-1}}^2
 \right](v-u)^{\beta_\gamma}.
\end{split}
\end{align}
On the complementary event, we expand the scalar product before applying
any estimate:
\begin{align*}
\begin{split}
 &\left\langle
 \Psi_r^{u,X}-\Psi_r^{u,Y},
 b^N(\Psi_r^{u,X})-g^N(\Psi_r^{u,Y})
 \right\rangle\\
 &\quad=
 \langle\Psi_r^{u,X},b^N(\Psi_r^{u,X})\rangle
 -\langle\Psi_r^{u,X},g^N(\Psi_r^{u,Y})\rangle\\
 &\qquad-\langle\Psi_r^{u,Y},b^N(\Psi_r^{u,X})\rangle
 +\langle\Psi_r^{u,Y},g^N(\Psi_r^{u,Y})\rangle .
\end{split}
\end{align*}
Each of these four terms has a flow as multiplier, whose size depends on
its initial point. To quantify this, let $\mathbf w:=V+\bar B^{t_0,u}$ and,
for $\xi,\xi'\in\R^d$, $\Psi_r^\xi:=\phi_r^{u,\xi,\mathbf w}$.
Conditionally on $\mathcal F_u$, the path $\mathbf w$ is deterministic, and
the norms in the next three displays are taken only with respect to the fresh
noise $B^u$. The zeroth-order part of the Malliavin--Sobolev norm depends on
the whole frozen input, not only on the local scale: by
\eqref{eq:flow-conditional-moments}, for every finite $q$,
\begin{align*}
 \|\Psi_r^\xi\|_{L_q}
 \le C(1+|\xi|)
 \bigl(1+[\mathbf w]_{C^{H-\eps}(I_u)}\bigr)^k.
\end{align*}
Every Malliavin derivative of positive order is controlled by
\Cref{prop:flow-estimates}; since
$R_{u,r}^{\mathbf w;I_u}
\le1+[\mathbf w]_{C^{H-\eps}(I_u)}$, the complete norm satisfies
\begin{align}\label{eq:flow-as-IBP-multiplier}
\begin{split}
 \|\Psi_r^\xi\|_{\mathbb D^{\lfloor-\gamma\rfloor+3,\bar q}}
 \le
 C(1+|\xi|)
 \bigl(1+[V+\bar B^{t_0,u}]_{C^{H-\eps}(I_u)}\bigr)^C.
\end{split}
\end{align}
On the other hand,
\begin{align*}
 \Lambda_{\lfloor-\gamma\rfloor+3,q}(\Psi_r^{\xi'})^{-(\gamma\wedge0)}
 \le
 C\bigl(R_{u,r}^{V+\bar B^{t_0,u};I_u}\bigr)^C
 (r-u)^{(\gamma\wedge0)H}.
\end{align*}
Consequently, \Cref{lem:malliavin-1} gives
\begin{align}\label{eq:one-point-flow-multiplier}
\begin{split}
 \left|\EE^u\langle\Psi_r^\xi,f(\Psi_r^{\xi'})\rangle\right|
 \le
 C\|f\|_{\mathcal C^\gamma}(1+|\xi|)
 \bigl(1+[V+\bar B^{t_0,u}]_{C^{H-\eps}(I_u)}\bigr)^C
 (r-u)^{(\gamma\wedge0)H}.
\end{split}
\end{align}
The dependence on the initial condition therefore comes only from the flow in
the first, multiplier position; the initial condition of the flow at which
$f$ is evaluated does not create another polynomial factor.  Since
$r-u\le1$, we may weaken $(r-u)^{(\gamma\wedge0)H}$ to
$(r-u)^{-(1-\gamma)H}$.  Applying
\eqref{eq:one-point-flow-multiplier} to the four evaluations gives,
with $G_{u,v}:=\{v-u\le\tau_{V+\bar B^{t_0,u};I_u}\}$,
\begin{align*}
 |A^N_{u,v}|
 &\lesssim\ind_{G_{u,v}}
 \bigl(|Z_u|^2+\|b-g\|_{\mathcal C^{\gamma-1}}|Z_u|\bigr)
 (v-u)^{\beta_\gamma}\\
 &\quad+\ind_{G_{u,v}^c}(1+|X_u|+|Y_u|)
 \bigl(1+[V+\bar B^{t_0,u}]_{C^{H-\eps}(I_u)}\bigr)^C
 (v-u)^{\beta_\gamma}.
\end{align*}
The second term uses the coarse bound on the whole integral; it does not
require that every shorter subinterval exceeds the local scale.
Put $q_0=2m+\delta_0$.  H\"older's inequality, with the factor
$1+|X_u|+|Y_u|$ in $L_{q_0}$ and the remaining factor in the conjugate
space relative to $L_m$, gives
\begin{align*}
 &\left\|
 \ind_{\{v-u>\tau_{V+\bar B^{t_0,u};I_u}\}}
 (1+|X_u|+|Y_u|)
 \bigl(1+[V+\bar B^{t_0,u}]_{C^{H-\eps}(I_u)}\bigr)^C
 \right\|_{L_m}\\
 &\quad\le
 C(1+C_0)
 \left\|
 \ind_{\{v-u>\tau_{V+\bar B^{t_0,u};I_u}\}}
 \bigl(1+[V+\bar B^{t_0,u}]_{C^{H-\eps}(I_u)}\bigr)^C
 \right\|_{L_{\frac{mq_0}{q_0-m}}}.
\end{align*}
\Cref{lem:bad-event-moments} therefore gives, for every $\rho\ge1$,
\begin{align}\label{eq:energy-germ-bad}
\begin{split}
 \left\|\ind_{\{v-u>\tau_{V+\bar B^{t_0,u};I_u}\}}
 A^N_{u,v}\right\|_{L_m}
 \le C^\rho\rho^{\rho/2}
 (v-u)^{\beta_\gamma+\chi_\eps\rho}.
\end{split}
\end{align}
The loss from $q_0$ to $m$ is fixed, so the constants remain of the form
$C^\rho\rho^{\rho/2}$. Equations \eqref{eq:energy-germ-good}--%
\eqref{eq:energy-germ-bad} imply
\begin{align}
\label{eq:energy-germ-total}
 \|A^N_{u,v}\|_{L_m}
 &\lesssim
 \left[
 \left(\sup_{r\in\mathbb J}\|Z_r\|_{L_{2m}}\right)^2
 +\|b-g\|_{\mathcal C^{\gamma-1}}^2
 \right](v-u)^{\beta_\gamma} +C^\rho\rho^{\rho/2}
 (v-u)^{\beta_\gamma+\chi_\eps\rho}.
\end{align}
The triangle inequality therefore gives the unconditional sewing defect
with exponent $\beta_\gamma>1/2$.

We next estimate the conditional defect at the midpoint $c=(u+v)/2$. As in
the proof of \Cref{lem:stability-sewing}, put $\xi^X:=\phi_c^{u,X_u}$ and
$\xi^Y:=\phi_c^{u,Y_u}$, and, for $z,z'\in\{\xi^X,\xi^Y,X_c,Y_c\}$, write
$\Psi_r^z:=\phi_r^{c,z}$ (with input $V+\bar B^{t_0,c}$) and
$\Psi_r^{z,z'}:=\Psi_r^z-\Psi_r^{z'}$. We also set $\mathfrak b^N:=b^N-g^N$.
Define the scalar defect integrand
\begin{align*}
 \mathcal D^N_{u,c,r}
 &:=
 \left\langle
 \Psi_r^{\xi^X,\xi^Y},
 b^N(\Psi_r^{\xi^X})-g^N(\Psi_r^{\xi^Y})
 \right\rangle -
 \left\langle
 \Psi_r^{X_c,Y_c},
 b^N(\Psi_r^{X_c})-g^N(\Psi_r^{Y_c})
 \right\rangle.
\end{align*}
Let $G_1$ and $G_2(r)$, $r\in[c,v]$, be as in the proof of
\Cref{lem:stability-sewing}.
Using the tower property and the flow property, the integrand of
$\EE^u\delta A^N_{u,c,v}$ is $\EE^c\mathcal D^N_{u,c,r}$.
The following algebraic decomposition specifies the four terms that will be
estimated:
\begin{align*}
 \mathcal D^N_{u,c,r}
 &=\mathcal D^N_{1;u,c,r}
 +\mathcal D^N_{2;u,c,r}
 +\mathcal D^N_{3;u,c,r}
 +\mathcal D^N_{4;u,c,r},
\end{align*}
where
\begin{align*}
 \mathcal D^N_{1;u,c,r}
 &:=
 \left\langle
 \Psi_r^{\xi^X,\xi^Y}-\Psi_r^{X_c,Y_c},
 b^N(\Psi_r^{\xi^X})-b^N(\Psi_r^{\xi^Y})
 \right\rangle,\\
 \mathcal D^N_{2;u,c,r}
 &:=
 \left\langle
 \Psi_r^{X_c,Y_c},
 b^N(\Psi_r^{\xi^X})-b^N(\Psi_r^{\xi^Y})
 -b^N(\Psi_r^{X_c})+b^N(\Psi_r^{Y_c})
 \right\rangle,\\
 \mathcal D^N_{3;u,c,r}
 &:=
 \left\langle
 \Psi_r^{\xi^X,\xi^Y}-\Psi_r^{X_c,Y_c},
 \mathfrak b^N(\Psi_r^{\xi^Y})
 \right\rangle,\\
 \mathcal D^N_{4;u,c,r}
 &:=
 \left\langle
 \Psi_r^{X_c,Y_c},
 \mathfrak b^N(\Psi_r^{\xi^Y})-\mathfrak b^N(\Psi_r^{Y_c})
 \right\rangle.
\end{align*}

We bound these terms on $G_1\cap G_2(r)$.  For
$\mathcal D^N_{1;u,c,r}$, apply \Cref{lem:malliavin-2} to the difference
$b^N(\Psi_r^{\xi^X})-b^N(\Psi_r^{\xi^Y})$.  The two flow estimates
\eqref{eq:upperbound-D-x} and \eqref{eq:upperbound-4} then give
\begin{align*}
 \left|\EE^c\mathcal D^N_{1;u,c,r}\right|
 &\lesssim
 (|\xi^X-\xi^Y|+|X_c-Y_c|)|\xi^X-\xi^Y-X_c+Y_c|
 (r-c)^{-(1-\gamma)H}\\
 &\quad+
 (|\xi^X-\xi^Y|+|X_c-Y_c|)^2
 (|\xi^X-X_c|+|\xi^Y-Y_c|)
 (r-c)^{-(2-\gamma)H}.
\end{align*}
For $\mathcal D^N_{2;u,c,r}$, use the exact identity
\begin{align*}
 &b^N(\Psi_r^{\xi^X})-b^N(\Psi_r^{\xi^Y})
 -b^N(\Psi_r^{X_c})+b^N(\Psi_r^{Y_c})\\
 &\quad=
 b^N(\Psi_r^{\xi^X})-b^N(\Psi_r^{\xi^Y})
 -b^N(\Psi_r^{X_c})
 +b^N(\Psi_r^{\xi^Y}+\Psi_r^{X_c}-\Psi_r^{\xi^X})\\
 &\qquad+
 b^N(\Psi_r^{Y_c})
 -b^N(\Psi_r^{\xi^Y}+\Psi_r^{X_c}-\Psi_r^{\xi^X}).
\end{align*}
\Cref{lem:malliavin-3} bounds the first four-point increment, while
\Cref{lem:malliavin-2} bounds the last difference.  Together with
\eqref{eq:upperbound-D-x} and \eqref{eq:upperbound-4}, this yields
\begin{align*}
 \left|\EE^c\mathcal D^N_{2;u,c,r}\right|
 &\lesssim
 (|\xi^X-\xi^Y|+|X_c-Y_c|)|\xi^X-\xi^Y-X_c+Y_c|
 (r-c)^{-(1-\gamma)H}\\
 &\quad+
 (|\xi^X-\xi^Y|+|X_c-Y_c|)^2
 (|\xi^X-X_c|+|\xi^Y-Y_c|)
 (r-c)^{-(2-\gamma)H}.
\end{align*}
Since
$\|\mathfrak b^N\|_{\mathcal C^{\gamma-1}}
\lesssim\|b-g\|_{\mathcal C^{\gamma-1}}$,
\eqref{eq:odd-one-point-IBP} with $\sigma=1-\gamma$ and
\eqref{eq:upperbound-4} give
\begin{align*}
 \left|\EE^c\mathcal D^N_{3;u,c,r}\right|
 &\lesssim
 \|b-g\|_{\mathcal C^{\gamma-1}}
 |\xi^X-\xi^Y-X_c+Y_c|
 (r-c)^{-(1-\gamma)H}\\
 &\quad+
 \|b-g\|_{\mathcal C^{\gamma-1}}
 (|\xi^X-\xi^Y|+|X_c-Y_c|)
 (|\xi^X-X_c|+|\xi^Y-Y_c|)
 (r-c)^{-(2-\gamma)H}.
\end{align*}
Finally, \eqref{eq:odd-two-point-IBP} with $\sigma=1-\gamma$ and
\eqref{eq:upperbound-D-x} give
\begin{align*}
 \left|\EE^c\mathcal D^N_{4;u,c,r}\right|
 &\lesssim
 \|b-g\|_{\mathcal C^{\gamma-1}}
 (|\xi^X-\xi^Y|+|X_c-Y_c|)
 (|\xi^X-X_c|+|\xi^Y-Y_c|)
 (r-c)^{-(2-\gamma)H}.
\end{align*}
Adding the four inequalities gives
\begin{align}\label{eq:energy-defect-pointwise}
\begin{split}
 &\left|\EE^c\mathcal D^N_{u,c,r}\right|\\
 &\quad\lesssim
 (|\xi^X-\xi^Y|+|X_c-Y_c|)\,|\xi^X-\xi^Y-X_c+Y_c|
 (r-c)^{-(1-\gamma)H}\\
 &\qquad+
 (|\xi^X-\xi^Y|+|X_c-Y_c|)^2
 (|\xi^X-X_c|+|\xi^Y-Y_c|)
 (r-c)^{-(2-\gamma)H}\\
 &\qquad+
 \|b-g\|_{\mathcal C^{\gamma-1}}
 |\xi^X-\xi^Y-X_c+Y_c|
 (r-c)^{-(1-\gamma)H}\\
 &\qquad+
 \|b-g\|_{\mathcal C^{\gamma-1}}
 (|\xi^X-\xi^Y|+|X_c-Y_c|)
 (|\xi^X-X_c|+|\xi^Y-Y_c|)
 (r-c)^{-(2-\gamma)H}.
\end{split}
\end{align}
Thus every term on the right contains either the double comparison
remainder $\xi^X-\xi^Y-X_c+Y_c=-\mathsf D_{u,c}$ or one of the individual
comparison remainders $\xi^X-X_c$ and $\xi^Y-Y_c$.

For the factors in \eqref{eq:energy-defect-pointwise}, conditional
Cauchy--Schwarz and the definition of $\mathcal Q_m(\mathbb J)$ give
\begin{align*}
 \left\|
 \bigl(\EE^u|\xi^X-\xi^Y-X_c+Y_c|^2\bigr)^{1/2}
 \right\|_{L_{2m}}
 \le(c-u)^{1/2}\mathcal Q_m(\mathbb J)^{1/2}.
\end{align*}
Moreover,
$|\xi^X-\xi^Y|\le e^{\kappa_4(c-u)}|Z_u|$.  Hence
\eqref{eq:conditional-mixed-moment} with $q=4$ and the triangle inequality
give
\begin{align*}
 &\left\|
 \bigl(\EE^u(|\xi^X-\xi^Y|+|X_c-Y_c|)^4\bigr)^{1/2}
 \right\|_{L_m} \lesssim
 \left(\sup_{r\in\mathbb J}\|Z_r\|_{L_{2m}}\right)^2.
\end{align*}
Finally, the good-event conditional estimates in
\eqref{eq:dissipativity-fixed-bounds} bound each of
$|\xi^X-X_c|$ and $|\xi^Y-Y_c|$ in every conditional moment used here by
$C(c-u)^{\alpha_\gamma}$.  Applying these three estimates term by term in
\eqref{eq:energy-defect-pointwise}, integrating over $r\in[c,v]$, and then
using $2ab\le a^2+b^2$, gives the good-event estimate
\begin{align}\label{eq:energy-defect-on-good}
\begin{split}
 &\left\|\EE^u\int_c^v
 \ind_{G_1\cap G_2(r)}
 \EE^c\mathcal D^N_{u,c,r}\,dr\right\|_{L_m}\\
 &\lesssim
 (c-u)^{1/2}(v-c)^{\beta_\gamma}
 \mathcal Q_m(\mathbb J)^{1/2}
 \left(
 \sup_{r\in\mathbb J}\|Z_r\|_{L_{2m}}
 +\|b-g\|_{\mathcal C^{\gamma-1}}\right)\\
 &\quad+
 (c-u)^{\alpha_\gamma}(v-c)^{\zeta_\gamma}
 \left[
 \left(\sup_{r\in\mathbb J}\|Z_r\|_{L_{2m}}\right)^2
 +\|b-g\|_{\mathcal C^{\gamma-1}}^2
 \right].
\end{split}
\end{align}
For use with the complementary indicators, we derive one coarse bound.  We first expand
$\mathcal D^N_{u,c,r}$ into its eight individual evaluations:
\begin{align}\label{eq:energy-defect-eight-evaluations}
\begin{split}
 \mathcal D^N_{u,c,r}
 &=
 \langle\Psi_r^{\xi^X},b^N(\Psi_r^{\xi^X})\rangle
 -\langle\Psi_r^{\xi^X},g^N(\Psi_r^{\xi^Y})\rangle
 -\langle\Psi_r^{\xi^Y},b^N(\Psi_r^{\xi^X})\rangle\\
 &\quad+\langle\Psi_r^{\xi^Y},g^N(\Psi_r^{\xi^Y})\rangle
 -\langle\Psi_r^{X_c},b^N(\Psi_r^{X_c})\rangle
 +\langle\Psi_r^{X_c},g^N(\Psi_r^{Y_c})\rangle\\
 &\quad+\langle\Psi_r^{Y_c},b^N(\Psi_r^{X_c})\rangle
 -\langle\Psi_r^{Y_c},g^N(\Psi_r^{Y_c})\rangle .
\end{split}
\end{align}
Let $\mathfrak M_{u,c}:=1+|\xi^X|+|\xi^Y|+|X_c|+|Y_c|$.
The one-point estimate \Cref{lem:malliavin-1} and the flow estimates
apply separately to the eight terms in
\eqref{eq:energy-defect-eight-evaluations}.  As in
\eqref{eq:one-point-flow-multiplier}, the only initial-condition factor comes
from the flow in the multiplier position.  More precisely,
\begin{align*}
 \left|\EE^c\langle\Psi_r^{z},b^N(\Psi_r^{z'})\rangle\right|
 &\lesssim (1+|z|)
 \bigl(1+[V+\bar B^{t_0,c}]_{C^{H-\eps}(I_c)}\bigr)^C
 (r-c)^{-(1-\gamma)H},
\end{align*}
where $z,z'\in\{\xi^X,\xi^Y,X_c,Y_c\}$.
The sum of the above inequalities gives the coarse
bound
\begin{align}\label{eq:energy-defect-coarse}
\begin{split}
 &\left|\EE^c\mathcal D^N_{u,c,r}\right| \lesssim
 \bigl(1+[V+\bar B^{t_0,c}]_{C^{H-\eps}(I_c)}\bigr)^C
 \mathfrak M_{u,c}
 (r-c)^{-(1-\gamma)H}.
\end{split}
\end{align}
We now justify the moment bound for $\mathfrak M_{u,c}$. Since
$\xi^X=\phi_c^{u,X_u}$ and $\xi^Y=\phi_c^{u,Y_u}$,
\eqref{eq:flow-conditional-moments} gives, for every finite $q$,
\begin{align}\label{eq:energy-starting-point-moments}
\begin{split}
 &\left(\EE^u\bigl[(1+|\xi^X|+|\xi^Y|)^q\bigr]\right)^{1/q} \le
 C(1+|X_u|+|Y_u|)
 \bigl(1+[V+\bar B^{t_0,u}]_{C^{H-\eps}(I_u)}\bigr)^C.
\end{split}
\end{align}
The remaining terms $X_c$ and $Y_c$ in $\mathfrak M_{u,c}$ are controlled
directly by the assumed marginal moments.

Put $q_0=2m+\delta_0$.  Apply H\"older's inequality with the state factors
from \eqref{eq:energy-starting-point-moments}, $X_c$, and $Y_c$ in
$L_{q_0}$ and the input factors in the conjugate space relative to $L_m$.
Then \Cref{lem:bad-event-moments} gives
\begin{align*}
 \left\|\ind_{G_1^c}\mathfrak M_{u,c}
 \bigl(1+[V+\bar B^{t_0,c}]_{C^{H-\eps}(I_c)}\bigr)^C
 \right\|_{L_m}
 &\le C^\rho\rho^{\rho/2}(c-u)^{\chi_\eps\rho},\\
 \left\|\ind_{G_2(r)^c}
 \mathfrak M_{u,c}
 \bigl(1+[V+\bar B^{t_0,c}]_{C^{H-\eps}(I_c)}\bigr)^C
 \right\|_{L_m}
 &\le C^\rho\rho^{\rho/2}(r-c)^{\chi_\eps\rho}.
\end{align*}
In the first inequality, the failed local scale is that of
$V+\bar B^{t_0,u}$ on $I_u$; in the second, that of $V+\bar B^{t_0,c}$ on
$I_c$. In both cases, all other input factors, including the one from
\eqref{eq:energy-starting-point-moments}, have fixed finite moments.

We now combine the estimates using the partition of unity
\begin{align*}
 1=\ind_{G_1\cap G_2(r)}
 +\ind_{G_1^c\cap G_2(r)}+\ind_{G_2(r)^c}.
\end{align*}
On the first term use \eqref{eq:energy-defect-pointwise}; on the other
two use \eqref{eq:energy-defect-coarse}. Taking conditional expectations, applying the two input-moment estimates
just proved, and integrating in $r$ gives
\begin{align}\label{eq:energy-defect-good}
\begin{split}
 \|\EE^u\delta A^N_{u,c,v}\|_{L_m}
 &\lesssim
 (c-u)^{1/2}(v-c)^{\beta_\gamma}
 \mathcal Q_m(\mathbb J)^{1/2}
 \left(\sup_{r\in\mathbb J}\|Z_r\|_{L_{2m}}
 +\|b-g\|_{\mathcal C^{\gamma-1}}\right)\\
 &\quad+(c-u)^{\alpha_\gamma}(v-c)^{\zeta_\gamma}
 \left[\left(\sup_{r\in\mathbb J}\|Z_r\|_{L_{2m}}\right)^2
 +\|b-g\|_{\mathcal C^{\gamma-1}}^2\right]\\
 &\quad+C^\rho\rho^{\rho/2}
 \Big[(c-u)^{\chi_\eps\rho}(v-c)^{\beta_\gamma}
 +(v-c)^{\beta_\gamma+\chi_\eps\rho}\Big].
\end{split}
\end{align}
The last line is bounded by
$C^\rho\rho^{\rho/2}(v-u)^{\beta_\gamma+\chi_\eps\rho}$,
since $c$ is the midpoint of $[u,v]$. The three conditional sewing
exponents are admissible because $\frac12+\beta_\gamma>1$,
$\alpha_\gamma+\zeta_\gamma>1$, and $\beta_\gamma+\chi_\eps\rho>1$ after
increasing $\rho$ if necessary. As in the preceding sewing argument,
the estimate for every prescribed $\rho\ge1$ follows by applying the bounds
with $\rho+\rho_0$ and absorbing the fixed shift into $C$.

The process
$\mathcal E^N_t:=\int_{t_0}^t w_r\langle Z_r,b^N(X_r)-g^N(Y_r)\rangle\,dr$ is
the process associated with the germ \eqref{eq:energy-sewing-germ} by
\Cref{lem:stochastic-sewing}, applied in $L_m$ (recall $m\ge2$). To check
the characterisation in that lemma, write, for $u<r$,
\begin{align*}
 &w_r\langle Z_r,b^N(X_r)-g^N(Y_r)\rangle
 -\EE^u\bigl[w_r\langle\Psi_r^{u,X}-\Psi_r^{u,Y},
 b^N(\Psi_r^{u,X})-g^N(\Psi_r^{u,Y})\rangle\bigr]\\
 &\quad=w_r\Bigl(\langle Z_r,b^N(X_r)-g^N(Y_r)\rangle
 -\langle\Psi_r^{u,X}-\Psi_r^{u,Y},
 b^N(\Psi_r^{u,X})-g^N(\Psi_r^{u,Y})\rangle\Bigr)+\mathfrak c_{u,r},
\end{align*}
where $\mathfrak c_{u,r}$ has zero conditional mean given $\mathcal F_u$.
Since $|\Psi_r^{u,X}-\Psi_r^{u,Y}|\le e^{\kappa_4}|Z_u|$ by
\Cref{lem:reg-phi}, the first term on the right is bounded by
$C_N\bigl(1+|Z_u|\bigr)\bigl(|X_r-\Psi_r^{u,X}|+|Y_r-\Psi_r^{u,Y}|\bigr)$,
with $C_N$ depending on $\|b^N\|_{C^1}+\|g^N\|_{C^1}$, and
$|\mathfrak c_{u,r}|\le C_N|Z_u|$. H\"older's inequality and
\eqref{eq:comparison-remainder-L2m} bound the $L_m$-norm of the first term by
$C_N(r-u)^{\varsigma_0}$. Integrating over $[u,v]$ gives
$\|\mathcal E^N_v-\mathcal E^N_u-A^N_{u,v}\|_{L_m}\le C_N(v-u)$ and
$\|\EE^u(\mathcal E^N_v-\mathcal E^N_u-A^N_{u,v})\|_{L_m}
\le C_N(v-u)^{1+\varsigma_0}$, which is admissible in that characterisation.
Hence
\Cref{lem:stochastic-sewing}, together
with \eqref{eq:energy-germ-total} and \eqref{eq:energy-defect-good}, give
\eqref{eq:energy-sewing-local} uniformly in $N$, with $\mathcal E$ replaced
by $\mathcal E^N$.

We next identify the limit. Let $m':=m+\delta_0/4$ and $\delta_0':=\delta_0/2$,
so that $2m'+\delta_0'=2m+\delta_0$. The standing assumptions remain valid
with $(m,\delta_0)$ replaced by $(m',\delta_0')$, and
$\mathcal Q_{m'}(\mathbb J)<\infty$. Indeed, the comparison estimates in
\eqref{eq:dissipativity-fixed-bounds} give
\begin{align*}
 \left\|\ind_{\{t-s\le\tau_{V+\bar B^{t_0,s};I_s}\}}
 \EE^s|\mathsf D_{s,t}|^2\right\|_{L_{m'}}
 \le C(t-s)^{2\alpha_\gamma}.
\end{align*}
On the complementary event, the proof of
\eqref{eq:conditional-bad-final}, with
$(m,\delta_0)$ replaced by $(m',\delta_0')$, gives, for every
$\rho\ge1$,
\begin{align*}
 \left\|\ind_{\{t-s>\tau_{V+\bar B^{t_0,s};I_s}\}}
 \EE^s|\mathsf D_{s,t}|^2\right\|_{L_{m'}}
 \le C^\rho\rho^\rho(t-s)^{1+2\chi_\eps\rho}.
\end{align*}
Since $2\alpha_\gamma>1$, these two bounds show that
$\mathcal Q_{m'}(\mathbb J)<\infty$. Hence
\eqref{eq:stability-smooth-occupation}, applied on a finite cover of
$[t_0,T]$, gives
\begin{align}\label{eq:KN-uniform-holder}
 \sup_{N\ge1}
 [\mathbb K^N]_{C^{\beta_\gamma}([t_0,T];L_{2m'})}<\infty.
\end{align}
By \Cref{def:solution}, $\mathbb K^N\to \mathbb K$ in probability in $C([t_0,T];\R^d)$.
Since
$2m'\beta_\gamma>1$, the
Garsia--Rodemich--Rumsey inequality applied to
\eqref{eq:KN-uniform-holder}, together with $\mathbb K^N_{t_0}=0$, yields
\begin{align*}
 \sup_{N\ge1}
 \left\|\sup_{t\in[t_0,T]}|\mathbb K_t^N|\right\|_{L_{2m'}}<\infty.
\end{align*}
After extracting an almost surely uniformly convergent subsequence and
using Fatou's lemma, the same bound holds for $\mathbb K$.  Therefore
$\sup_t|\mathbb K_t^N-\mathbb K_t|^{2m}$ is uniformly integrable.  The convergence in
probability from \Cref{def:solution} consequently improves to
\begin{align}\label{eq:KN-uniform-L2m}
 \left\|\sup_{t\in[t_0,T]}|\mathbb K_t^N-\mathbb K_t|\right\|_{L_{2m}}
 \longrightarrow0.
\end{align}
Fatou's lemma applied to the increments in
\eqref{eq:KN-uniform-holder} also gives
$\mathbb K\in C^{\beta_\gamma}([t_0,T];L_{2m'})$.

Let $\beta'\in(1/2,\beta_\gamma)$. The interpolation inequality
\eqref{eq:holder-interpolation}, applied to $\mathbb K^N-\mathbb K$ with
$(\beta,\alpha)=(\beta',\beta_\gamma)$, bounds
$[\mathbb K^N-\mathbb K]_{C^{\beta'}([t_0,T];L_{2m})}$ by
$2(\sup_t\|\mathbb K_t^N-\mathbb K_t\|_{L_{2m}})^{1-\beta'/\beta_\gamma}
[\mathbb K^N-\mathbb K]_{C^{\beta_\gamma}([t_0,T];L_{2m})}^{\beta'/\beta_\gamma}$.
The second factor is uniformly bounded by
\eqref{eq:KN-uniform-holder} and its Fatou limit, whereas the first tends
to zero by \eqref{eq:KN-uniform-L2m}.  Thus
\begin{align}\label{eq:KN-holder-convergence}
 \mathbb K^N\longrightarrow \mathbb K
 \quad\text{in }C^{\beta'}([t_0,T];L_{2m}).
\end{align}

We now pass to the limit in $\mathcal E^N$. Since $Z=U^0+\mathbb K$ and $U=wU^0$,
\begin{align*}
 \mathcal E_t^N-\mathcal E_s^N
 =
 \int_s^t\langle U_r,d\mathbb K_r^N\rangle
 +\int_s^t w_r\langle \mathbb K_r,d\mathbb K_r^N\rangle.
\end{align*}
For the first integral, ordinary integration by parts gives
\begin{align*}
 \int_s^t\langle U_r,d\mathbb K_r^N\rangle
 =
 \langle U_t,\mathbb K_t^N\rangle-\langle U_s,\mathbb K_s^N\rangle
 -\int_s^t\langle \mathbb K_r^N,dU_r\rangle.
\end{align*}
Consequently,
\begin{align*}
 &\left|
 \int_s^t\langle U_r,d\mathbb K_r^N\rangle
 -\left(
 \langle U_t,\mathbb K_t\rangle-\langle U_s,\mathbb K_s\rangle
 -\int_s^t\langle \mathbb K_r,dU_r\rangle
 \right)\right|\\
 &\quad\le
 \left(|U_s|+|U_t|
 +\operatorname{Var}_{[s,t]}(U)\right)
 \sup_{r\in[s,t]}|\mathbb K_r^N-\mathbb K_r|.
\end{align*}
The right-hand side converges to zero in probability, uniformly in
$(s,t)$, by \eqref{eq:KN-uniform-L2m}.

For the second integral, use the continuous bilinear pairing
$L_{2m}\times L_{2m}\to L_m$.  The Banach-valued Young estimate and
\eqref{eq:KN-holder-convergence} give
\begin{align*}
 &\left\|
 \int_s^t w_r\langle \mathbb K_r,d\mathbb K_r^N\rangle
 -\int_s^t w_r\langle \mathbb K_r,d\mathbb K_r\rangle
 \right\|_{L_m} \le
 C\|w\mathbb K\|_{C^{\beta'}([s,t];L_{2m})}
 \|\mathbb K^N-\mathbb K\|_{C^{\beta'}([s,t];L_{2m})}
 \longrightarrow0,
\end{align*}
where the Young integral exists because $2\beta'>1$.  The Young product
formula, first for smooth paths and then by continuity in
$C^{\beta'}$, yields
\begin{align*}
 \int_s^t w_r\langle \mathbb K_r,d\mathbb K_r\rangle
 &=
 \frac12w_t|\mathbb K_t|^2-\frac12w_s|\mathbb K_s|^2
 -\frac12\int_s^t|\mathbb K_r|^2\,dw_r.
\end{align*}
The two preceding limits show that
$\mathcal E_t^N-\mathcal E_s^N\to\mathcal E_t-\mathcal E_s$ in probability
for every $(s,t)\in\Delta_{\mathbb J}$.  Passing to an almost surely
convergent subsequence and applying Fatou's lemma to the sewing estimate for
$\mathcal E^N$, which is uniform in $N$, proves \eqref{eq:energy-sewing-local}.

We finally prove the energy inequality. Since $U^0$ has finite variation
and $dU^0=d\mathfrak F$,
\begin{align*}
 (w|U^0|^2)_{s,t}
 =
 2\int_s^t w_r\langle U^0_r,d\mathfrak F_r\rangle
 +\int_s^t|U^0_r|^2\,dw_r.
\end{align*}
On the other hand, the definition of $\mathcal E$ gives
\begin{align*}
 2(\mathcal E_t-\mathcal E_s)
 &=
 2\langle U,\mathbb K\rangle_{s,t}
 -2\int_s^t\langle \mathbb K_r,dU_r\rangle +(w|\mathbb K|^2)_{s,t}
 -\int_s^t|\mathbb K_r|^2\,dw_r .
\end{align*}
Expanding $w|Z|^2=w|U^0|^2+2\langle U,\mathbb K\rangle+w|\mathbb K|^2$, using
$dU_r=w_r\,d\mathfrak F_r+U^0_r\,dw_r$, and substituting the two preceding identities, we obtain
\begin{align*}
 (w|Z|^2)_{s,t}
 &=
 2(\mathcal E_t-\mathcal E_s)
 +2\int_s^t w_r
 \left(
 \langle Z_r,F(X_r)-F(Y_r)\rangle
 -\kappa_4|Z_r|^2
 \right)dr .
\end{align*}
By \eqref{eq:F-one-sided-global}, the last integral is nonpositive.
This proves \eqref{eq:singular-energy-identity}.
\end{proof}

We can now state the final stability result, which is where the
compatibility condition \eqref{eq:H-k-iteration-condition} enters. Under
\eqref{eq:H-k-iteration-condition}, the range
\begin{align}\label{eq:stability-epsilon-range}
0<\eps<H\wedge\left(H+\frac1{k-1}-\frac{(1-\gamma)H}{2} \right)
\end{align}
is non-empty, and for every $\eps$ in this range,
$\chi_\eps>(1-\gamma)H/2$, that is, $2\chi_\eps>1-\beta_\gamma$.

\begin{proposition}[Stability]\label{prop:stab}
Assume \eqref{eq:H-k-iteration-condition}, let $\theta\in(0,1)$, and place
ourselves in the standing setting of this subsection with $\eps$ satisfying
\eqref{eq:stability-epsilon-range}. Let $X_{t_0}$ and $Y_{t_0}$ be the initial
conditions of $X$ and $Y$. There exists a constant $C$, depending only on
$d,H,\gamma,k,\boldsymbol\kappa,\eps,T-t_0,\Lambda_{V,\eps},
m,\delta_0,\theta,\Xi$, and $C_0$, such that
\begin{align}\label{eq:stability-final}
 \sup_{t\in[t_0,T]}\|X_t-Y_t\|_{L_{2m}}
 \le
 C\left(
 \|X_{t_0}-Y_{t_0}\|_{L_{2m}}
 +\|b-g\|_{\mathcal C^{\gamma-1}}
 \right)^{1-\theta}.
\end{align}
\end{proposition}

\begin{proof}
Fix an interval $\mathbb J=[S,S']\subset[t_0,T]$ of length $h\le h_0$. By
Step 4 of the proof of \Cref{lem:stability-dissipativity},
$\mathcal Q_m(\mathbb J)<\infty$, and hence
$\mathcal K_m(\mathbb J)<\infty$ by \Cref{lem:stability-sewing}.
Combining
\Cref{lem:stability-dissipativity,lem:stability-sewing}, squaring
\eqref{eq:K-stability-local}, and using
$2ab\le a^2+b^2$, we obtain, after decreasing $h_0$ if necessary,
\begin{align}\label{eq:conditional-square-closed}
\begin{split}
 \mathcal Q_m(\mathbb J)
 \le Ch^\eta\left[
 \left(\sup_{r\in\mathbb J}\|X_r-Y_r\|_{L_{2m}}\right)^2
 +\|b-g\|_{\mathcal C^{\gamma-1}}^2
 \right]
 +C^\rho\rho^\rho h^{2\chi_\eps\rho}.
\end{split}
\end{align}
Indeed, squaring the last term in
\eqref{eq:K-stability-local} gives
$C^\rho\rho^\rho h^{2\chi_\eps\rho}$, the same moment profile as
the complementary-event term in
\eqref{eq:dissipativity-stability-local}. The only term to be reabsorbed
after the substitution is $Ch^{\eta+1}\mathcal Q_m(\mathbb J)$.

We now use the energy estimate.  Taking the $L_m$-norm in
\eqref{eq:singular-energy-identity}, applying
\eqref{eq:energy-sewing-local}, and then using
\eqref{eq:conditional-square-closed} gives, uniformly for $t\in\mathbb J$,
\begin{align}\label{eq:local-energy-recursion}
\begin{split}
 \|X_t-Y_t\|_{L_{2m}}^2
 &\le e^{2\kappa_4h}\|X_S-Y_S\|_{L_{2m}}^2\\
 &\quad+Ch^{\beta_\gamma}
 \left[
 \left(\sup_{r\in\mathbb J}\|X_r-Y_r\|_{L_{2m}}\right)^2
 +\|b-g\|_{\mathcal C^{\gamma-1}}^2
 \right]\\
 &\quad+C^\rho\rho^{\rho/2}
 h^{\beta_\gamma+\chi_\eps\rho}.
\end{split}
\end{align}
The mixed term in \eqref{eq:energy-sewing-local} is bounded using
$2ab\le a^2+b^2$; because
\eqref{eq:conditional-square-closed} contains the additional factor
$h^\eta$, it is of strictly higher order than the second line of
\eqref{eq:local-energy-recursion}. For its complementary-event part,
taking $\mathcal Q_m(\mathbb J)^{1/2}$ turns
$C^\rho\rho^\rho h^{2\chi_\eps\rho}$ into
$C^\rho\rho^{\rho/2}h^{\chi_\eps\rho}$ after a change of $C$, and the
factor $\sup_{r\in\mathbb J}\|Z_r\|_{L_{2m}}+\|b-g\|_{\mathcal C^{\gamma-1}}$
multiplying it is bounded by $C(C_0+\Xi)$; this gives the last line of
\eqref{eq:local-energy-recursion}.

Choose $h_0$ smaller once more so that the first term containing the
supremum over $\mathbb J$ in \eqref{eq:local-energy-recursion} can be absorbed.
Let $\mathfrak d_{b,g}:=\|b-g\|_{\mathcal C^{\gamma-1}}$.
For a uniform mesh $t_j=t_0+jh$, with $h=(T-t_0)/N\le h_0$, let
$e_j:=\|X_{t_j}-Y_{t_j}\|_{L_{2m}}^2$ and
$E_j:=\sup_{r\in[t_j,t_{j+1}]}\|X_r-Y_r\|_{L_{2m}}^2$.
Taking the supremum over each mesh interval in
\eqref{eq:local-energy-recursion}, and using
$e^{2|\kappa_4|h}\le1+Ch^{\beta_\gamma}$, gives
\begin{align}\label{eq:energy-scalar-recurrence}
\begin{split}
 E_j
 &\le(1+Ch^{\beta_\gamma})e_j
 +Ch^{\beta_\gamma}\mathfrak d_{b,g}^2
 +C^\rho\rho^{\rho/2}
 h^{\beta_\gamma+\chi_\eps\rho}.
\end{split}
\end{align}
Since $e_{j+1}\le E_j$, the discrete Gr\"onwall inequality applied to
\eqref{eq:energy-scalar-recurrence} yields
\begin{align}\label{eq:energy-discrete-gronwall}
\begin{split}
 \max_{0\le j\le N}e_j+\max_{0\le j<N}E_j
 \le
 C\exp\!\left(Ch^{-(1-\beta_\gamma)}\right)
 \left[
 s_*+
 \bigl(1+C^\rho\rho^{\rho/2}\bigr)
 h^{\chi_\eps\rho-(1-\beta_\gamma)}
 \right],
\end{split}
\end{align}
where $s_*:=\|X_{t_0}-Y_{t_0}\|_{L_{2m}}^2+\mathfrak d_{b,g}^2$.
The loss $h^{-(1-\beta_\gamma)}$ in the last term is the cost of
summing the $N\asymp h^{-1}$ local remainders of size
$h^{\beta_\gamma+\chi_\eps\rho}$.

It remains to choose $h$ and $\rho$. The relation
$2\chi_\eps>1-\beta_\gamma$, which holds by
\eqref{eq:stability-epsilon-range}, is used at the end of this step, to
bound the Gr\"onwall factor $\exp(Ch^{-(1-\beta_\gamma)})$ in
\eqref{eq:energy-discrete-gronwall} by an arbitrarily small negative power
of $s_*$. For a sufficiently small $h$,
choose $\rho=c\,h^{-2\chi_\eps}$ with $c>0$ small enough (and $h$ small enough that
$\rho\ge1$). This choice produces the cancellation
\begin{align*}
 \log\!\left(
 C^\rho\rho^{\rho/2}h^{\chi_\eps\rho}
 \right)
 &=
 \rho\log C+\frac{\rho}{2}\log\rho
 +\chi_\eps\rho\log h =
 \rho\left(\log C+\frac12\log c\right).
\end{align*}
Choose $c$ sufficiently small that the last parenthesis is negative. It
follows that, for some $a_0>0$,
\begin{align*}
 C^\rho\rho^{\rho/2}h^{\chi_\eps\rho}
 \le e^{-a_0h^{-2\chi_\eps}}.
\end{align*}
The term $h^{\chi_\eps\rho}$ has at least the same decay, and the polynomial
factor $h^{-(1-\beta_\gamma)}$ can be absorbed by decreasing
$a_0$. Consequently, there exists $a_1>0$, depending only on the
fixed parameters, such that
\begin{align*}
 \bigl(1+C^\rho\rho^{\rho/2}\bigr)
 h^{\chi_\eps\rho-(1-\beta_\gamma)}
 \le2e^{-a_1h^{-2\chi_\eps}}.
\end{align*}
If $0<s_*<e^{-1}$, choose the integer $N$ so that
$a_1h^{-2\chi_\eps}$ is comparable to $2\log(1/s_*)$; all the estimates
above are unchanged up to constants. Since
$(1-\beta_\gamma)/(2\chi_\eps)<1$, for every $\theta\in(0,1)$,
\begin{align*}
 \exp\!\left(Ch^{-(1-\beta_\gamma)}\right)
 \le C_\theta s_*^{-\theta}.
\end{align*}
Inserting the last two estimates into
\eqref{eq:energy-discrete-gronwall} gives
\begin{align*}
 \sup_{r\in[t_0,T]}\|X_r-Y_r\|_{L_{2m}}^2
 \le C_\theta s_*^{1-\theta}.
\end{align*}
Since $s_*^{1/2}\le\|X_{t_0}-Y_{t_0}\|_{L_{2m}}+\mathfrak d_{b,g}$, taking square
roots proves \eqref{eq:stability-final}.  If $s_*$ is bounded away from
zero, the conclusion follows by fixing one admissible mesh and enlarging
the constant.  If $s_*=0$, apply the preceding estimate with $s_*$ replaced
by an arbitrary positive number and then let that number tend to zero.
\end{proof}

\section{Proof of the main theorems}\label{sec:proof-main}

We construct the solution of \Cref{thm:main} as the limit of the solutions of
equations in which $b$ is replaced by its heat regularisations. The condition
\eqref{eq:H-k-iteration-condition} permits us to use \Cref{prop:stab} both
to prove convergence of the approximating equations and to obtain
uniqueness; \Cref{thm:stability} then follows from the same proposition.
Throughout this section, we fix
\begin{align*}
 \eps_*:=\frac12\min\left\{H,\,
 H+\frac1{k-1}-\frac{(1-\gamma)H}{2}\right\}>0,
\end{align*}
which satisfies \eqref{eq:stability-epsilon-range} and depends only on
$H,k,\gamma$. All equations below are constructed on the probability space
fixed in \Cref{sec:framework}, with the same input $V$ and the same
Riemann--Liouville process $B^{t_0}$.

\subsection{Existence}

We first solve equations with smooth drifts and pass to a limit in the uniform
a priori estimates. Set $b^N:=P_{1/N}b$, $N\ge1$.
By \Cref{ass:singular-drift},
\begin{align}\label{eq:canonical-drift-properties}
\begin{split}
 b^N\in C_b^\infty(\R^d;\R^d),
 \qquad
 \sup_{N\ge1}\|b^N\|_{\mathcal C^\gamma}\lesssim\Xi,
 \qquad
 b^N\longrightarrow b\quad\text{in }\mathcal C^{\gamma-1}.
\end{split}
\end{align}
For every $N$, let $X^N=X^{t_0,X_{t_0},V;b^N}$ be the unique global solution of
the smooth equation with drift $F+b^N$, that is,
\begin{align*}
 X_t^N
 =X_{t_0}+\int_{t_0}^tF(X_r^N)\,dr
 +\widehat K_t^N
 +V_{t_0,t}+B_t^{t_0},
 \qquad
 \widehat K_t^N:=\int_{t_0}^tb^N(X_r^N)\,dr.
\end{align*}

\paragraph{Uniform estimates.}
Fix $m\ge2$ and $\eps\in(0,H)$. By
\Cref{prop:regsol-inf,prop:regsol,lem:smooth-solution-moments,prop:regsol-usual},
the constants in the a priori estimates are uniform in $N$. In particular,
for every sufficiently large moment index $p$,
\begin{align}\label{eq:uniform-approximation-estimates}
\begin{split}
 \sup_{N\ge1}\Bigg(
 \sup_{r\in[t_0,T]}\|X_r^N\|_{L_p}
 +[X^N-V-B^{t_0}]_{C^{\alpha_\gamma}([t_0,T];L_p)}
 \Bigg)<\infty,
\end{split}
\end{align}
and the conditional comparison estimate defining
$\mathcal V_\eps(\infty)$ holds uniformly in $N$.

\paragraph{Cauchy property of the smooth solutions.}
The uniform estimates above allow us to compare two approximating equations in
a topology strong enough to pass to the limit. By
\eqref{eq:canonical-drift-properties},
$d_{N,M}:=\|b^N-b^M\|_{\mathcal C^{\gamma-1}}\to0$ as $N,M\to\infty$.
Let $Z^{N,M}:=X^N-X^M$. The bounds in \eqref{eq:uniform-approximation-estimates} and the uniform
conditional comparison estimates imply that the quantitative hypotheses of
\Cref{prop:stab}, with $\eps=\eps_*$, hold uniformly in $N,M$ at every
fixed moment order.
Indeed, for a given $q\ge2$, one may take $\delta_0=1$ and use the uniform
$L_{2q+1}$ marginal bound. Applying the proposition with moment parameter
$q$, $X_{t_0}=Y_{t_0}$, $b=b^N$, and $g=b^M$ gives, for every
$\theta\in(0,1)$,
\begin{align*}
 \sup_{r\in[t_0,T]}\|Z_r^{N,M}\|_{L_{2q}}
 \le C_{\theta,q}d_{N,M}^{\,1-\theta}.
\end{align*}
Since $q$ is arbitrary, $d_{N,M}\to0$ and the monotonicity of
$L_p$-norms on a probability space imply directly that, for every
$q\in[2,\infty)$,
\begin{align}\label{eq:XNM-all-moments}
 \sup_{r\in[t_0,T]}\|Z_r^{N,M}\|_{L_q}\longrightarrow0.
\end{align}
Moreover, the triangle inequality applied to the two uniform estimates in
\eqref{eq:uniform-approximation-estimates}, followed if necessary by
monotonicity of $L_p$-norms, gives, for every $q\in[2,\infty)$,
\begin{align*}
 \sup_{N,M\ge1}
 [Z^{N,M}]_{C^{\alpha_\gamma}([t_0,T];L_q)}<\infty.
\end{align*}

Since $\alpha_\gamma>1/2$, the interpolation inequality
\eqref{eq:holder-interpolation} with $(\beta,\alpha)=(1/2,\alpha_\gamma)$
applies. Applying it to $Z^{N,M}$ and using
\eqref{eq:XNM-all-moments} gives, for every $q\in[2,\infty)$,
\begin{align}\label{eq:XNM-global-control}
 [X^N-X^M]_{C^{1/2}([t_0,T];L_q)}\longrightarrow0.
\end{align}

Fix $m\ge2$, choose $q>m$, and choose
$\iota\in(1/q,1/2)$. The quantitative
Kolmogorov and Garsia--Rodemich--Rumsey estimate applied to
\eqref{eq:XNM-global-control} yields
\begin{align*}
 &\left\|
 \sup_{(s,t)\in\Delta_{[t_0,T]}}
 \frac{|Z_t^{N,M}-Z_s^{N,M}|}
 {(t-s)^{1/2-\iota}}
 \right\|_{L_q} \le
 C
 [Z^{N,M}]_{C^{1/2}([t_0,T];L_q)}
 \longrightarrow0.
\end{align*}
The restriction $\iota>1/q$ is the usual loss when passing from an $L_q$
increment estimate to a pathwise H\"older estimate. Since $Z_{t_0}^{N,M}=0$ and $q>m$, it follows that
\begin{align*}
 \left\|\sup_{t\in[t_0,T]}|Z_t^{N,M}|\right\|_{L_m}
 &\le
 (T-t_0)^{1/2-\iota}
 \left\|
 \sup_{(s,t)\in\Delta_{[t_0,T]}}
 \frac{|Z_t^{N,M}-Z_s^{N,M}|}
 {(t-s)^{1/2-\iota}}
 \right\|_{L_q} \longrightarrow0.
\end{align*}
Thus $(X^N)_{N\ge1}$ is Cauchy in $L_m(\Omega;C([t_0,T];\R^d))$, and we
denote its limit by $X$:
\begin{align}\label{eq:XN-path-convergence}
 X^N\longrightarrow X
 \quad\text{in }L_m\bigl(\Omega;C([t_0,T];\R^d)\bigr).
\end{align}
The limit does not depend on the chosen value of $m$: for any two admissible
moments, both constructions are limits in probability in
$C([t_0,T];\R^d)$ of the same sequence $(X^N)_{N\ge1}$, and limits in
probability are unique.
Moreover, letting $M\to\infty$ in \eqref{eq:XNM-global-control} and using
lower semicontinuity of the H\"older seminorm gives
\begin{align}\label{eq:XN-controlled-convergence}
\begin{split}
 \sup_{r\in[t_0,T]}\|X_r^N-X_r\|_{L_m}
 +[X^N-X]_{C^{1/2}([t_0,T];L_m)}
 \longrightarrow0.
\end{split}
\end{align}
The process $X$ is continuous because the convergence takes place in
$L_m(\Omega;C([t_0,T];\R^d))$. For each fixed $t$, the space of
$\mathcal F_t$-measurable elements of $L_m$ is closed. Since
$X_t^N\to X_t$ in $L_m$, the process $X$ is also adapted.

\paragraph{Passage of the a priori estimates to the limit.}
Having identified the limit process, we next transfer the uniform temporal and
conditional comparison estimates. The temporal estimate in
\Cref{prop:regsol-usual} passes to the limit by
lower semicontinuity and gives, with $q_m\ge m$ the moment index provided
by that proposition,
\begin{align}\label{eq:limit-temporal-regularity}
\begin{split}
 [X-V-B^{t_0}]_{C^{\alpha_\gamma}([t_0,T];L_m)}
 \le C\left(
 1+\|X_{t_0}\|_{L_{q_m}}
 +\left\|[V]_{C^{H-\eps}([t_0,T])}\right\|_{L_{q_m}}
 \right)^C.
\end{split}
\end{align}
Since $\alpha_\gamma>1/2$, this proves $X\in\mathcal V(m)$ for every $m\ge2$.

Fix now any $\eps\in(0,H)$, $m\ge2$, and $(s,t)\in\Delta_{[t_0,T]}$.
By \eqref{eq:XN-path-convergence} and the pathwise Lipschitz estimate for the
flow,
the comparison remainders
$\mathfrak R_{s,t}^{X^N}:=X_t^N-\phi_t^{s,X_s^N,V+\bar B^{t_0,s}}$ converge
to $\mathfrak R_{s,t}^{X}:=X_t-\phi_t^{s,X_s,V+\bar B^{t_0,s}}$ in $L_m$. After taking an almost surely convergent subsequence, conditional
Fatou's lemma gives
\begin{align*}
 \EE^s|\mathfrak R_{s,t}^{X}|^m
 \le\liminf_{N\to\infty}\EE^s|\mathfrak R_{s,t}^{X^N}|^m
 \qquad\text{almost surely}.
\end{align*}
The event $\{t-s\le\tau_{V+\bar B^{t_0,s};I_s}\}$ does not depend on $N$. Therefore the uniform estimate of
\Cref{prop:regsol} passes to the limit and yields
\begin{align*}
 \left\|
 \ind_{\{t-s\le\tau_{V+\bar B^{t_0,s};I_s}\}}
 \bigl(\EE^s|\mathfrak R_{s,t}^{X}|^m\bigr)^{1/m}
 \right\|_{L_\infty}
 \le C(t-s)^{\alpha_\gamma}.
\end{align*}
Taking the supremum over $(s,t)$ proves that
$X\in\mathcal V_\eps(\infty)$. Since $\eps\in(0,H)$ was arbitrary,
\begin{align}\label{eq:constructed-all-eps}
 X\in\bigcap_{\eps\in(0,H)}\mathcal V_\eps(\infty).
\end{align}

\paragraph{Construction of the drift process and independence of the
regularisation.}
It remains to identify the limiting distributional drift and to show that it
does not depend on the chosen approximation. Define directly from the limiting
process
\begin{align}\label{eq:limit-drift-definition}
 K_t
 :=X_t-X_{t_0}-\int_{t_0}^tF(X_r)\,dr-V_{t_0,t}-B_t^{t_0},
 \qquad t\in[t_0,T].
\end{align}
Then $K$ is adapted and continuous, $K_{t_0}=0$, and item (i) of
\Cref{def:solution} holds by construction. Using the smooth equation for $X^M$ and the definition
\eqref{eq:limit-drift-definition}, we have
\begin{align*}
 K_t-\widehat K_t^M
 =X_t-X_t^M+\int_{t_0}^t\bigl(F(X_r^M)-F(X_r)\bigr)\,dr.
\end{align*}
The first term tends to zero in $L_m(\Omega;C([t_0,T]))$ by
\eqref{eq:XN-path-convergence}. For the second one, choose a sufficiently
large moment index $q$ and use the polynomial Lipschitz estimate for $F$,
H\"older's inequality, the uniform moments, and
\eqref{eq:XN-controlled-convergence}. It follows that
\begin{align}\label{eq:K-canonical-convergence}
\begin{split}
 \widehat K^M\longrightarrow K
 \quad\text{in }L_m\bigl(\Omega;C([t_0,T];\R^d)\bigr).
\end{split}
\end{align}

Suppose first that $\gamma>0$. By interpolation,
\eqref{eq:canonical-drift-properties} implies $\|b^M-b\|_\infty\to0$.
Since $b\in\mathcal C^\gamma$ and $X^M\to X$ in
$L_q(\Omega;C([t_0,T]))$ for every finite $q$, we also have
\begin{align*}
 \int_{t_0}^{\cdot}b^M(X_r^M)\,dr
 \longrightarrow
 \int_{t_0}^{\cdot}b(X_r)\,dr
 \quad\text{in }L_m\bigl(\Omega;C([t_0,T];\R^d)\bigr).
\end{align*}
Comparing this convergence with \eqref{eq:K-canonical-convergence} gives
$K_t=\int_{t_0}^tb(X_r)\,dr$ for all $t\in[t_0,T]$. If $(c^j)_{j\ge1}$ is any admissible approximation, interpolation between
its uniform $\mathcal C^\gamma$ bound and its convergence in
$\mathcal C^{\gamma-1}$ gives $\|c^j-b\|_\infty\to0$.  Therefore
\begin{align*}
 \sup_{t\in[t_0,T]}
 \left|\int_{t_0}^tc^j(X_r)\,dr-K_t\right|
 \le(T-t_0)\|c^j-b\|_\infty
 \longrightarrow0.
\end{align*}
This proves the approximation requirement when $\gamma>0$.

It remains to consider $\gamma\le0$. We first record an occupation bound
for the smooth solutions, in which the integrand is measured in
$\mathcal C^{\gamma-1}$. Let $f\in C_b^\infty(\R^d;\R^d)$ and $M\ge1$. We
apply \cite[Proposition~3.3]{DHL25}, which holds for an arbitrary adapted
process $\psi$, with $\psi:=X^M-B^{t_0}$, $q=m$, the regularity index
$\gamma-1$ in place of $\gamma$, and the H\"older exponent $\alpha_\gamma$
for $\psi$. 
Its hypotheses hold:
$\gamma-1>-1/(2H)$ is the condition $(1-\gamma)H<1/2$, and
$H(\gamma-2)+\alpha_\gamma>0$ because
$\alpha_\gamma-(2-\gamma)H=1-2H+2\gamma H>0$ when $\gamma<0$, while
$1-(2-\gamma)H>0$ when $\gamma=0$. Since $V$ is
$\mathcal F_{t_0}$-measurable, the centred conditional seminorm of $\psi$
coincides with that of $X^M-V-B^{t_0}$, which, as in the proof of
\Cref{prop:regsol-usual}, is bounded by twice its
$C^{\alpha_\gamma}(L_m)$-seminorm, hence uniformly in $M$ by
\eqref{eq:uniform-approximation-estimates}. Since
$1+H(\gamma-1)=\beta_\gamma$ and
$1+H(\gamma-2)+\alpha_\gamma\ge\beta_\gamma$, we obtain, for
$(s,t)\in\Delta_{[t_0,T]}$ with $t-s\le1$,
\begin{align*}
 \left\|\int_s^tf(X_r^M)\,dr\right\|_{L_m}
 \le C\|f\|_{\mathcal C^{\gamma-1}}(t-s)^{\beta_\gamma},
\end{align*}
with $C$ independent of $f$ and $M$. As $m\ge2$ and $\beta_\gamma>1/2$, we
have $m\beta_\gamma>1$, and the Kolmogorov continuity theorem gives
\begin{align}\label{eq:occupation-C-gamma-minus-one}
\begin{split}
 \left\|\sup_{t\in[t_0,T]}
 \left|\int_{t_0}^tf(X_r^M)\,dr\right|\right\|_{L_m}
 \le C\|f\|_{\mathcal C^{\gamma-1}},
 \qquad f\in C_b^\infty(\R^d;\R^d),\ M\ge1.
\end{split}
\end{align}

Let $(c^j)_{j\ge1}$ be any admissible approximation of $b$, and set
$\widetilde K_t^j:=\int_{t_0}^tc^j(X_r)\,dr$. For every $j,M$ and
$t\in[t_0,T]$,
\begin{align}\label{eq:regularisation-independence-decomposition}
\begin{split}
 K_t-\widetilde K_t^j
 &=\bigl(K_t-\widehat K_t^M\bigr)+\int_{t_0}^t\bigl(b^M-c^j\bigr)(X_r^M)\,dr +\int_{t_0}^t\bigl(c^j(X_r^M)-c^j(X_r)\bigr)\,dr.
\end{split}
\end{align}
Fix $j$ and let $M\to\infty$. The first term tends to zero in
$L_m(\Omega;C([t_0,T]))$ by \eqref{eq:K-canonical-convergence}. Since $c^j$
is smooth with bounded derivative, the third term is bounded by
$(T-t_0)\|\nabla c^j\|_\infty\sup_{r\in[t_0,T]}|X_r^M-X_r|$, which also
tends to zero in $L_m$ by \eqref{eq:XN-path-convergence}. By
\eqref{eq:occupation-C-gamma-minus-one} with $f=b^M-c^j$, the second term is
bounded in $L_m(\Omega;C([t_0,T]))$ by
$C\|b^M-c^j\|_{\mathcal C^{\gamma-1}}
\le C\|b^M-b\|_{\mathcal C^{\gamma-1}}+C\|b-c^j\|_{\mathcal C^{\gamma-1}}$,
and the first part vanishes as $M\to\infty$. Therefore
\begin{align*}
 \left\|\sup_{t\in[t_0,T]}|K_t-\widetilde K_t^j|\right\|_{L_m}
 \le C\|b-c^j\|_{\mathcal C^{\gamma-1}}
 \longrightarrow0
 \qquad\text{as }j\to\infty,
\end{align*}
since $c^j\to b$ in $\mathcal C^{\gamma-1}$, the constant being independent
of $j$.
Thus $K$ satisfies the approximation requirement for every admissible
regularisation also when $\gamma\le0$. The pair $(X,K)$ is a solution in the sense of
\Cref{def:solution}. Estimate \eqref{eq:thm-reg} is
\eqref{eq:limit-temporal-regularity}, and the existence part of
\Cref{thm:main} follows.

\subsection{Uniqueness}

Let $X$ be the solution constructed above and let $Y$ be another solution
with the same initial condition, input, and Riemann--Liouville process.
Assume that $Y\in\mathcal V_{\eps_Y}(\infty)$ and
$\sup_{t\in[t_0,T]}\|Y_t\|_{L_{4+\delta_0}}<\infty$ for some $\delta_0>0$
and some $\eps_Y$ in the range \eqref{eq:stability-epsilon-range}.
By \eqref{eq:constructed-all-eps} and
\eqref{eq:limit-temporal-regularity}, the constructed solution belongs to
$\mathcal V_{\eps_Y}(\infty)$ and has moments of every finite order. The
comparison estimate at the largest required order $p_k$, together with
conditional Jensen's inequality, supplies one constant for every
$p\in[2,p_k]$. Thus the two bounds in
\eqref{eq:dissipativity-fixed-bounds} hold for $X$ and $Y$ with a common
finite $C_0$, and \Cref{prop:stab} applies with $m=2$, the above
$\delta_0$, $\eps=\eps_Y$, $g=b$, and $X_{t_0}=Y_{t_0}$. Its right-hand side is
zero, so that $\sup_{t\in[t_0,T]}\|X_t-Y_t\|_{L_4}=0$. Thus $X_t=Y_t$ almost surely for every fixed $t$. Since both processes have
continuous trajectories, they are indistinguishable. This completes the
proof of \Cref{thm:main}.

\subsection{Proof of \texorpdfstring{\Cref{thm:stability}}{the stability theorem}}

Let $X$ and $Y$ be the solutions given by \Cref{thm:main} with initial
conditions $X_{t_0},Y_{t_0}$ and singular drifts $b,g$. By
\eqref{eq:constructed-all-eps} and the uniform estimates of
\Cref{prop:regsol,lem:smooth-solution-moments}, passed to the limit as
above, both solutions have finite marginal moments of every order and
satisfy the comparison estimate \eqref{eq:classVinf} at $\eps_*$, with
bounds controlled by sufficiently high moments of $X_{t_0}$, $Y_{t_0}$ and by the
input moment bounds in \eqref{eq:assumpV}. Hence the standing setting of
\Cref{sec:stab} holds with $\eps=\eps_*$, $\delta_0=1$, and a constant $C_0$
depending only on these quantities, and \Cref{prop:stab} gives
\eqref{eq:main-stability}. Since $\eps_*$ depends only on $H,k,\gamma$, no
additional free parameter $\eps$ appears in \Cref{thm:stability}.

\section{Concluding remarks}\label{sec:conclusion}

We have shown that the regularisation by a fractional noise and a
superlinear dissipation can be combined within one framework: the
regularising flow of the dissipative part replaces the free noise in the
regularisation estimates, and the local scale controls the inputs for which
the flow estimates degenerate. Several questions remain open.
\begin{enumerate}[label=(\roman*)]
  \item The compatibility condition \eqref{eq:H-k-iteration-condition},
 which is automatic for $k\le5$ and restricts the parameters only when
 $H<\frac14-\frac1{k-1}$, enters only through
 the final closure in \Cref{prop:stab}; we do not know whether it can be
 removed.
 \item The critical case $(1-\gamma)H=1/2$, which for $F=0$ requires a
 critical version of the stochastic sewing lemma with a logarithmic
 correction, is not covered. For $H>1/2$, the second condition in
 \eqref{eq:subcritical} is the more restrictive one, and the range
 $1-\frac1{2H}<\gamma\le2-\frac1H$ is not covered either.
 \item For Lipschitz dissipative drifts, the long-time behaviour and the
 invariant measures were studied in \cite{DHL25}. Extending these results
 to the superlinear setting considered here would require estimates that
 are uniform in time, rather than on bounded intervals.
\end{enumerate}

\end{document}